\documentclass[a4paper,11pt]{article}
\usepackage{amsmath,latexsym,amssymb,amsfonts,natbib}
\usepackage[a4paper,hmargin=1.87 cm]{geometry}
\usepackage[dvips]{graphicx}
\usepackage{pifont}
\usepackage{pstcol,pst-fill,pst-grad}
\usepackage{rotating}
\usepackage{subfigure}
\usepackage{graphics}
\usepackage{hyperref}
\usepackage{multirow}
\usepackage{bigints}
\usepackage{tikz} 
\usetikzlibrary{positioning, fit, shapes.geometric}
\usepackage{multicol}
\usepackage{courier}
\usepackage[latin1]{inputenc}
\usepackage{indentfirst}
\usepackage{ccaption}
\usepackage{verbatim}
\usepackage{makeidx}
\usepackage{listings}
\usepackage{dsfont}
\usepackage[toc,page]{appendix}

\newcommand{\proof}     {\paragraph{Proof}}
\newcommand{\carre}     {\hfill$\Box$}
\numberwithin{equation}{section}
\newtheorem{defi}{Definition}[section]
\newtheorem{lem}{Lemma}[section]
\newtheorem{theo}{Theorem}[section]
\newtheorem{coro}{Corollary}[section]

\newtheorem{rem}{Remark}

\newtheorem{ass}{Assumption}
\newtheorem{ex}{Example}[section]

\title{Large-time behavior of finite-state mean-field systems with multi-classes}
\author{ Donald A. Dawson\thanks{School of Mathematics and Statistics, Carleton University, 1125 Colonel By Drive
Ottawa, Ontario K1S 5B6, Canada. ddawson@math.carleton.ca}, Ahmed Sid-Ali\thanks{School of Mathematics and Statistics,
Carleton University, 1125 Colonel By Drive Ottawa, Ontario K1S 5B6, Canada. ahmedsidali@cunet.carleton.ca}, Yiqiang Q. Zhao\thanks{School of Mathematics and Statistics, Carleton University, 1125 Colonel By Drive Ottawa, Ontario K1S 5B6, Canada. zhao@math.carleton.ca} \\
{\small }\\
{\small }\\
}
\date{}
\begin{document}
\maketitle

\begin{abstract}
We study in this paper the large-time asymptotics of the empirical vector associated with a family of finite-state mean-field systems with multi-classes. The empirical vector is composed of local empirical measures characterizing the different classes within the system. As the number of particles in the system goes to infinity, the empirical vector process converges towards the solution to a McKean-Vlasov system. First, we investigate the large deviations principles of the invariant distribution from the limiting McKean-Vlasov system. Then, we examine the metastable phenomena arising at a large scale and large time. Finally, we estimate the rate of convergence of the empirical vector process to its invariant measure. Given the local homogeneity in the system, our results are established in a product space. 
\end{abstract}

{\it 2020 Mathematics Subject Classification:} Primary 60K35 60F10 60J74; Secondary 60G10.
\\
\\
{\it Keywords:} Large deviations; Mean-field; Invariant measure; Metastability; McKean-Vlasov; Jump processes. 

\section*{Introduction}
Interacting particle systems with multi-classes, widely encountered in a variety of domains going from statistical physics, chemistry, communication networks, biology, to finance, have recently attracted the interest of many researchers and several models were proposed to understand their large-scale behavior. See, e.g.,\cite{Chong+Klu2019,Collet2014,Collet+al2016,Graha2008,Kno+Lo+schu2020,Mey2020,Nguyen+al2020} and the references therein for an overview of recent advances on the subject. In multi-classes systems, the particles come from different sub-populations within which they are homogeneous and thus, the entire system is heterogeneous but composed of homogeneous sub-populations. Thence, one can average over the local symmetries within the different classes to describe mean-field interactions through local empirical measures. Thus, the entire system is characterized by an empirical vector composed of the local empirical measures. 

The focus in the current article is on a particular family of mean-field multi-class models describing the evolution of block-structured networks with dynamically changing multi-color nodes. The state-space here is the finite set of colors and the particles are identified as the nodes of the network. This class of models was proposed in \cite{Daw+Sid+zha2020} to describe the dynamic of various physical phenomena and the large-scale asymptotics were established. In particular, a multi-class propagation of chaos was proved to hold together with a law of large numbers, implying the convergence of the empirical vector towards the solution to a McKean-Vlasov system of equations as the total number of particles $N$ in the system goes to infinity. Moreover, the authors studied the large deviations principle for the empirical vector process over finite time intervals.    

We propose in this article to study the large-time behavior of the family of models introduced in \cite{Daw+Sid+zha2020}. Our motivation comes from the interesting characteristics of these systems as well as the importance of their large-time behavior for many applications. Notice that numerous works on the large-time behavior of various systems of interacting particles exist in the literature. See, for example, \cite{Cox+Grev90,Daw+Grev93, Daw+Grev99, Daw+Grev+Vai95,  Dor+Myt2013, Grev+den2007, Kuehn2015} and references therein. Therefore, our current contribution aims to be in the continuation of the aforementioned references. The approach taken in this document is summarized as follows.

First of all, let us mention that the propagation of chaos and the law of large numbers established in \cite{Daw+Sid+zha2020} are valid over finite time intervals. Thence, on any finite time horizon (not too large), one could use the McKean-Vlasov limit system as an approximation for the large $N$-particle system as it is classical in the mean-field literature. However, when time tends to infinity, these results are no longer necessarily true, and care should be taken in using this approximation. Indeed, as we detail throughout this article, the validity of the approximation is intimately linked to the critical points of the McKean-Vlasov limit system and its stability. Intuitively, if the McKean-Vlasov limit system has several $\omega$-limit sets, one can wonder which of these sets characterizes the large-time behavior of the large $N$-particle system?

As a starting point to tackle this question, we study in Section \ref{LDP-IM-sec} the asymptotics of the invariant measure of the empirical vector. In particular, we establish, under mild conditions, the large deviations principles for the invariant measure in two different cases. First, when there exists a unique asymptotically stable equilibrium to the McKean-Vlasov system (Theorem \ref{LDP-invariant-GAS}), and then when the McKean-Vlasov system has multiple $\omega$-limit sets (Theorem \ref{LDP-invariant-MEP}). The latter case is studied under the classical hypothesis of Friedlin and Wentzell. Note that the large deviations principles of the invariant distribution have been established for interacting diffusions in \cite{Daw+Gart89}, and for finite-state mean-field systems on complete graphs in \cite{Bork+Sund2012}. We adopt here the control theory approach introduced in \cite{Bis+Bork2011} for small noise diffusions and extended in \cite{Bork+Sund2012} to interacting jump processes on complete graphs. We therefore further extend this approach to the heterogeneous case with block-structured interaction graphs.

Next, we investigate the metastable phenomena, in the sense of Friedlin and Wentzell small noise stochastic systems, emerging when the McKean-Vlasov limiting system has several $\omega$-sets. Namely, in the case of multiple attractors, the empirical vector process, associated with a large but finite number of particles $N$, would transit between these attractors when the time is large. One then aims to estimate the most probable order in which these transitions occur, and the average time spent in the neighborhood of each of the $\omega$-limit sets. Understanding these phenomena is of great practical interest and hence our motivation. We describe in Section \ref{Meta-sect} the metastable phenomena in more detail and give important estimates. For this, we adopt the classical approaches of \cite[Chapter 6]{Freid+Wentz2012} and \cite {Hwang+Sheu90}. The main ingredients are the large deviation properties of the empirical vector process over finite time periods established in Section \ref{LDP-EP-sec}, and the main tool is the Freidlin-Wentzell quasipotential derived from the rate function characterizing the large deviations principle. Subsequently, we develop in Theorem \ref{rate-conv-invar-meas} an estimate of the time required for the empirical vector process to converge to its invariant measure. We find out that when the time is of the order $ \exp\{N (\Lambda + \delta) \} $, for all $\delta> 0$ and $ \Lambda$ being an appropriate constant, the empirical vector process is very close to its invariant measure. Interestingly, this coincides with the timescale found in \cite{Hwang+Sheu90} for diffusion processes and in \cite{Bork+Sund2012} for finite-state mean-field systems on complete graphs. We underline that the metastable phenomena in the section \ref{Meta-sect} are studied under the Hypothesis \ref{ass-mult-omeg}. Note, however, that if these assumptions are not imposed, things are more complicated and this is not pursued here. The reader may consult \cite{Bouc+Gaw2016, Tang+Yuan2017,zho+Ali2012} and the references therein for discussions and examples. Also, we emphasize that the results obtained in this study are for empirical vectors and thus, on product spaces.

The rest of this paper is organized as follows. In Section \ref{model} we revisit the family of models introduced in \cite{Daw+Sid+zha2020}. One may also consult \cite[Sec. 2]{Daw+Sid+zha2020} for a full detailed description together with some examples. Section \ref{LDP-EP-sec} is dedicated to the large deviations properties of the empirical vector process over finite time intervals. We first recall the principal results of \cite{Daw+Sid+zha2020} and then introduce some additional results used in the following sections. Then, we prove in Section \ref{LDP-IM-sec} our first set of the main results. Namely, Theorem \ref{LDP-invariant-GAS} gives the large deviations principle for the empirical measures when the McKean-Vlasov system has a unique globally asymptotically stable equilibrium, and Theorem \ref{LDP-invariant-MEP} gives the large deviations principle of the invariant measure when there are multiple $\omega$-limit sets. In addition, we recall some important concepts from the Friedlin-Wentzell large deviations theory. We then investigate in Section \ref{Meta-sect} the metastability of the $N$-particles system. Based on a set of results we provide estimates of the metastable transitions. Finally, Theorem \ref{rate-conv-invar-meas} gives the time required for the empirical vector process to converge towards its invariant measure.  

Since the results of Sections \ref{LDP-IM-sec} and \ref{Meta-sect} rely on the Freidlin-Wentzell program described in \cite[Ch. 6]{Freid+Wentz2012}, we give in Appendix \ref{Freid-Went-appe} the generalization of some of these results to our current setting. Finally, to facilitate the reading, we leave the lengthy proofs of some technical results in Appendix \ref{proof-appen}.

\section{The setting}
\label{model}

 Consider a graph $\mathcal{G}=(\mathcal{V},\Xi)$ composed of $r$ blocks $C_1,\ldots,C_r$ of sizes $N_1,\ldots,N_r$, respectively, where $\mathcal{V}$ is the set of the nodes and $\Xi$ is the set of the edges. Denote by $|\mathcal{V}|=N_1+\cdots+N_r=N$ the total number of the nodes in the network. Moreover, suppose that each block $C_j$ is a clique, i.e., all the $N_j$ nodes of the same block are connected to each other. We divide the nodes of each block $C_j$ into two sets: 
      \begin{itemize}
       \item \textbf{The central nodes $C^c_j$:} connected to all the other nodes of the same block but not to any node from the other blocks. We set $|C_j^c|=N_j^c$.
       \item \textbf{The peripheral nodes $C^p_j$:} connected to all the other nodes of the same block and to all the peripheral nodes of the other blocks. We set $|C^p_j|=N^p_j$.  
      \end{itemize} 
         
 Let $\mathcal{Z}=\{1,2,\ldots,K\}\subset\mathbb{N}$ be a finite set of $K$ colors. Suppose that each node of the graph $\mathcal{G}=(\mathcal{V},\Xi)$ is colored by one of the $K$ colors at each time. For each $1\leq j\leq r$ and $n\in C_j^c$ (resp. $n\in C_j^p$), denote by $(X_{n}(t),t\geq 0)$ the stochastic jump process that describes the evolution of the color of the node $n$ through time. Let $(\mathcal{Z},\mathcal{E})$ be the directed graph where $\mathcal{E}\subset\mathcal{Z}\times\mathcal{Z}\backslash \{(z,z)| z \in\mathcal{Z}\}$ describes the set of admissible jumps. In addition, whenever $(z,z')\in\mathcal{E}$, a node colored by $z$ is allowed to jump from $z$ to $z'$ at a rate that depends on the current state of the node and the state of its neighbors (adjacent nodes). To characterize these neighborhoods, we introduce, for each block $1\leq j\leq r$, the following local empirical measures describing respectively the state of the central and the peripheral nodes of the $j$-th block at time $t$
\begin{equation}
\begin{split}
\mu_j^{c,N}(t)=\frac{1}{N_j^c}\sum_{n\in C^c_j}\delta_{X_{n}(t)}\quad\mbox{and}\quad \mu_j^{p,N}(t)=\frac{1}{N_j^p}\sum_{n\in C^p_j}\delta_{X_{n}(t)},
\label{local-emp-meas-def}
\end{split}
\end{equation}   
and taking values in the set $\mathcal{M}_1(\mathcal{Z})$ of probability measures over $\mathcal{Z}$, endowed with the topology of weak convergence. The random dynamic in each block $1\leq j\leq r$ is thus summarized as follows 
\begin{itemize}
\item \textbf{The central nodes dynamic.} Each central node $n\in C^c_j$ jumps from color $z$ to $z'$, with $(z,z')\in (\mathcal{Z},\mathcal{E})$, at rate
\begin{align}
 \lambda_{z,z'}^c\big(\mu^{c,N}_j(t),\mu^{p,N}_j(t)\big),
\label{lamb-c}
\end{align}
 which depends on its current state and the states of its neighbors through the empirical measures $\mu^{c,N}_j(t)$ and $\mu^{p,N}_j(t)$.

\item \textbf{The peripheral nodes dynamic.} Each peripheral node $n\in C^p_j$ jumps from color $z$ to $z'$, with $(z,z')\in (\mathcal{Z},\mathcal{E})$, at rate  
\begin{equation}
\begin{split}
\lambda^p_{z,z'}\big(\mu_j^{c,N}(t),\mu_1^{p,N}(t),\ldots,\mu_r^{p,N}(t)\big),
\label{lamb-p}
\end{split}
\end{equation}
which depends on its state and the states of its neighbors through the local empirical measures $\mu_j^{c,N}(t),\mu_1^{p,N}(t),\ldots,\mu_r^{p,N}(t)$. 
\end{itemize}

Notice that the rate functions $\lambda^c_{z,z'}$ and $\lambda^p_{z,z'}$ depend on the number of nodes $N_j^c$ and $N_j^p$ within each category, but we omit this dependency to not overload the notations. See \cite[Sec. 2]{Daw+Sid+zha2020} for a detailed description. 

\begin{figure}[h]
\centering
\includegraphics[scale=0.15]{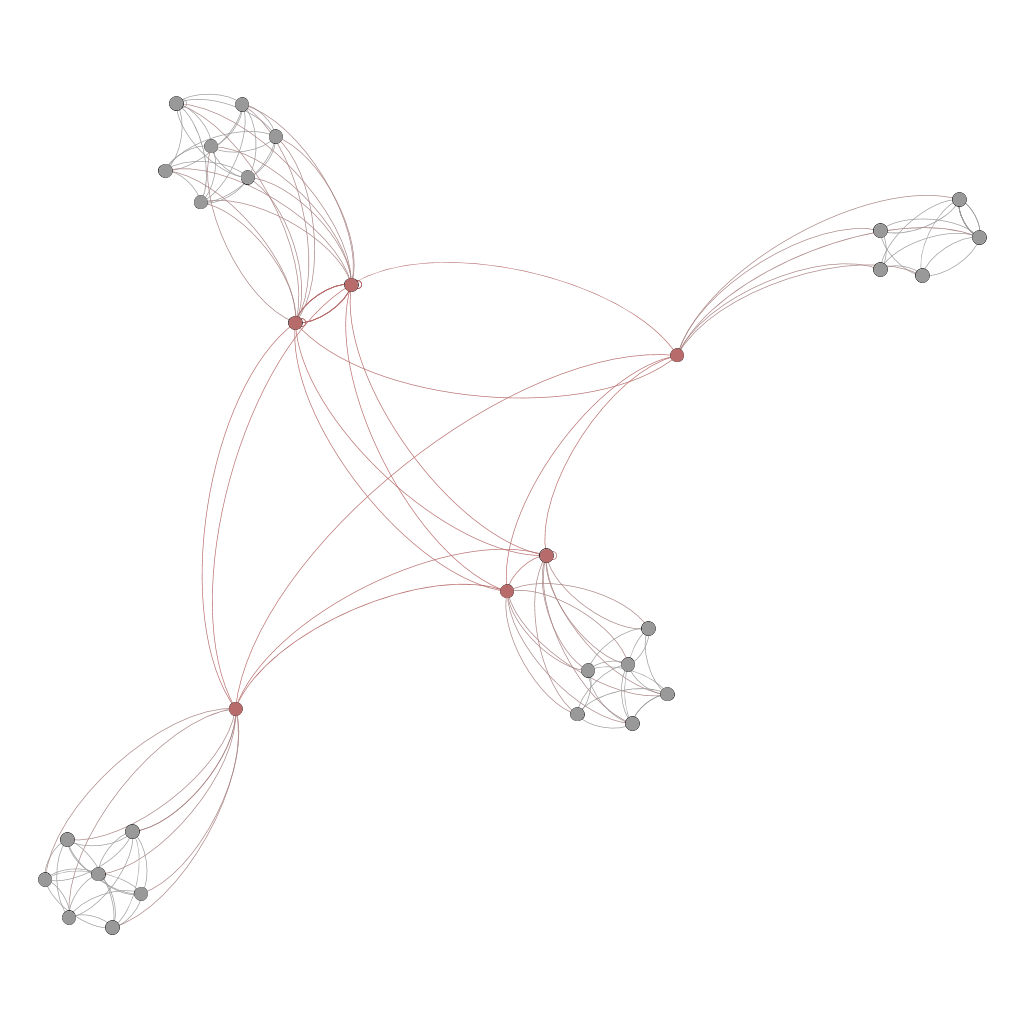}
\caption{A graph composed of four blocks with different numbers of central (grey) and peripheral (red) nodes. }
\end{figure}
We make the following assumptions throughout the paper.
\begin{ass}
\begin{enumerate}
\item The directed graph $(\mathcal{Z},\mathcal{E})$ is irreducible.
\item The rate functions $\lambda^c_{z,z'}$ and $\lambda^p_{z,z'}$ are Lipschitz and uniformly bounded away from zero, that is, there exists $c > 0$ such that for all $(z,z')\in\mathcal{E}$, $\lambda^c_{z,z'}(\cdot)\geq c$ and $\lambda_{z,z'}^p(\cdot)\geq c$.
\item For each block $1\leq j\leq r$, there exist $p_j^c,p_j^p,\alpha_j\in(0,1)$ such that, as $N\rightarrow\infty$,
\begin{align}
\frac{N_j}{N}\rightarrow \alpha_j,\quad \frac{N_j^p}{N_j}\rightarrow p_j^p,\quad \frac{N_j^c}{N_j}\rightarrow p_j^c,\quad p^p_j+p^c_j=1,\quad\text{and}\quad\sum_j\alpha_j=1.
\label{p-regul}
\end{align}
\end{enumerate}
\label{ass-prin}
\end{ass}

\begin{rem}
Notice that since $\mathcal{M}_1(\mathcal{Z})$ is compact and the rate functions $\lambda^c_{z,z'}$ and $\lambda^p_{z,z'}$ are continuous since Lipschitz, the rates are uniformly bounded from above, that is, there exists a constant $C <\infty$ such that for all $(z,z')\in\mathcal{E}$, we have $\lambda^c_{z,z'}(\cdot)\leq C$ and $\lambda^p_{z,z'}(\cdot)\leq C$. 
\end{rem}

\section{Large deviations of the empirical measures}
\label{LDP-EP-sec}
This section gathers large deviations results of the empirical vector process over finite time intervals. We start by introducing some additional notations. 

Denote by $\mathcal{D}([0,T],\mathcal{Z})$ the Skorokhod space of \textit{c\`adl\`ag} functions defined on $[0,T]$ with values in $\mathcal{Z}$, and by $\mathcal{M}_1(\mathcal{D}([0,T],\mathcal{Z}))$ the set of the probability measures over it. Denoting by $X^N=(X_n,X_m,n\in C_j^c,m\in C_j^p, 1\leq j\leq N)\in\mathcal{D}([0,T],\mathcal{Z}^N)$ the full description of the $N$ particles over the finite time interval $[0,T]$, we let $M^N\in\left(\mathcal{M}_1(\mathcal{D}([0,T],\mathcal{Z}))\right)^{2r}$  be the vector of "historical" empirical measures defined by
\begin{equation}
\begin{split}
M^N&=\bigg(M_1^{c,N},M_1^{p.N},\cdots,M_r^{c,N},M_r^{p,N}\bigg)=\bigg(\frac{1}{N_1^c}\sum_{n\in C^c_1}\delta_{X_n},\frac{1}{N_1^p}\sum_{n\in C^p_1}\delta_{X_n},\ldots, \frac{1}{N^c_r}\sum_{n\in C^c_r}\delta_{X_n}, \frac{1}{N^p_r}\sum_{n\in C^p_r}\delta_{X_n}  \bigg),
\label{bar-vector}
\end{split}
\end{equation}
where  $M_j^{c,N}$ (resp.  $M_j^{p,N}$) is the historical empirical measure of the central (resp. peripheral) nodes of the $j$-th block. With a slight abuse of notations, denote by $G_N$ the mapping that takes the full description $X^N$ to the empirical measures vector $M^N$, that is,
\begin{align*}
G_N: (X_n, 1\leq n\leq N)\in\mathcal{D}([0,T],\mathcal{Z}^N)\rightarrow M^N.
\end{align*}
Thus, $M^N=G_N(X^N)$. Denote by $\mathbb{P}_{z^N}^N$ the law of $X^N$ with initial condition $z^N=(z_n,z_m,n\in C_j^c,m\in C_j^p, 1\leq j\leq N)$. Note that the distribution of  $M^N$ depends on the initial condition only through its empirical vector defined by
\begin{align}
\nu^N=\left(\nu_N^{1,c},\nu_N^{1,p},\ldots,\nu_N^{r,c},\nu_N^{r,p} \right)=\bigg(\frac{1}{N_1^c}\sum_{n\in C_1^c}\delta_{z_n},\frac{1}{N_1^c}\sum_{n\in C_1^p}\delta_{z_n},\ldots,\frac{1}{N_r^c}\sum_{n\in C_r^c}\delta_{z_n},\frac{1}{N_r^p}\sum_{n\in C_r^p}\delta_{z_n}\bigg).
\label{init-emp-vect}
\end{align}

Define by $P_{\nu^N}^N=\mathbb{P}_{z^N}^N\circ G_N^{-1}$ the distribution  of $M^N$ which is the pushforward of $\mathbb{P}_{z^N}^N$ under the mapping $G_N$. 
 
Recalling $(\ref{local-emp-meas-def})$, consider the $(\mathcal{M}_1(\mathcal{Z}))^{2r}$-valued empirical vector process defined as
\begin{equation*}
\begin{split}
\mu^N: t\in [0,T]\longrightarrow \mu^N(t)=\left(\mu_1^{c,N}(t),\mu_1^{p,N}(t),\cdots,\mu_r^{c,N}(t),\mu_r^{p,N}(t)\right).
\end{split}
\end{equation*}

Notice that $\mu^N(0)=\nu_N$ and $\mu^N(t)$ is the projection $\pi_t(M^N)$, at time $t$, of $M^N$, that is,
\begin{align*}
\mu^N=\pi(M^N)=\pi(G_N(X^N)).
\end{align*} 

Denote by $p_{\nu_N}^N$ the distribution of $\mu^N$. The flow $\mu^N$ takes values in the product space $\big(\mathcal{D}([0,T],\mathcal{M}_1(\mathcal{Z}))\big)^{2r}$. Also, denote by $P_{z_n}$ the $n$-th particle's law  with initial condition $z_n$ in the case of non-interaction. Thus, the law of the entire non-interacting system is given by $\mathbb{P}_{z^n}^{0,N}=\otimes_{n=1}^NP_{z_n}$. Moreover, denote by $P_{\nu^N}^{0,N}$ the distribution of the corresponding empirical vector where  $\nu_N$ is the initial empirical vector. Therefore the Radon-Nikodym derivative $dP_{\nu^N}^{N}/dP_{\nu^N}^{0,N}$ at any $\mathbf{Q}=(Q_{1}^c,Q_{1}^p,\cdots, Q_{r}^c,Q_{r}^p)\in \big( \mathcal{M}_1(\mathcal{D}([0,T],\mathcal{Z}))\big)^{2r}$ is given by (see \cite[eqn. $(4.10)$]{Daw+Sid+zha2020})
\begin{equation}
\begin{split}
\frac{dP_{\nu^N}^{N}}{dP_{\nu^N}^{0,N}}(\mathbf{Q})&=\exp\bigg\{\sum_{j=1}^r\bigg[N_j^c\int_{D([0,T],\mathcal{Z})}h_1(x, \pi(Q_{j}^c),\pi(Q_{j}^p))Q_{j}^c(dx)\\
&\qquad\qquad+N_j^p\int_{D([0,T],\mathcal{Z})}h_2(x, \pi(Q_{j}^c), \pi(Q_{1}^p),\ldots, \pi(Q_{r}^p))Q_{j}^p(dx)\bigg]\bigg\}\\
                                                   &=\exp\big\{N h(\mathbf{Q})\big\},
\label{rad-nik}
\end{split}
\end{equation}
with
\begin{equation}
\begin{split}
h(\mathbf{Q})&=\sum_{j=1}^r\bigg[\frac{N_j^c}{N}\int_{D([0,T],\mathcal{Z})}h_1(x,\pi(Q_{j}^c),\pi(Q_{j}^p))Q_{j}^c(dx)\\
&\qquad\qquad+\frac{N_j^p}{N}\int_{D([0,T],\mathcal{Z})}h_2(x,\pi(Q_{j}^c), \pi(Q_{1}^p),\ldots, \pi(Q_{r}^p))Q_{j}^p(dx)\bigg],
\label{h-func}
\end{split}
\end{equation}
where, for any $\eta,\rho_1,\ldots, \rho_r$ in $\mathcal{D}([0,T], \mathcal{M}_1 (\mathcal{Z}))$,
\begin{align}
h_1(x,\eta,\rho_j)&= \sum_{0\leq t\leq T}\mathds{1}_{\{x_t\neq x_{t-}\}}\log \bigg(\lambda^c_{x_{t-},x_t}(\eta (t-),\rho_j(t-))\bigg)\nonumber \\
       &\qquad-\int_0^T\bigg(\sum_{z:(x_t,z)\in\mathcal{E}}\lambda^c_{x_t,z}(\eta(t),\rho_j(t))-1\bigg)dt,
       \label{h1-func}
\end{align}
and
\begin{align}
h_2(x,\eta,\rho_1,\ldots,\rho_r)&=\sum_{0\leq t\leq T}\mathds{1}_{\{x_t\neq x_{t-}\}}\log \bigg(\lambda^p_{x_t,x_{t-}}(\eta (t-),\rho_1(t-),\ldots,\rho_r(t-)\bigg)\nonumber \\
             &\qquad-\int_0^T\bigg(\sum_{z:(x_t,z)\in\mathcal{E}}\lambda^p_{x_t,z}(\eta (t),\rho_1(t),\ldots,\rho_r(t))-1\bigg)dt.
             \label{h2-func}
\end{align}

Equip the Skorokhod space $\mathcal{D}([0,T],\mathcal{M}_1(\mathcal{Z}))$ with the metric
\begin{align}
\rho_T (\mu,\nu)=\sup_{0\leq t\leq T}\rho_0 (\mu_t,\nu_t),\quad\mu,\nu\in\mathcal{D}([0,T],\mathcal{M}_1(\mathcal{Z})),
\label{rho-metric}
\end{align}
where $\rho_0(\alpha,\beta)$, for $\alpha,\beta\in\mathcal{M}_1(\mathcal{Z})$, is a metric on $\mathcal{M}_1(\mathcal{Z})$ which generates the weak topology on $\mathcal{M}_1(\mathcal{Z})$. Moreover, define by $\rho^{2r}_0(\cdot,\cdot)$ the product metric that generates the weak topology on the product space $(\mathcal{M}_1(\mathcal{Z}))^{2r}$. Furthermore, let the product space $\big(\mathcal{D}([0,T],\mathcal{M}_1(\mathcal{Z}))\big)^{2r}$ be equipped with the product topology induced by the product metric $\rho^{2r}_T=\max\{\rho_T,\cdots,\rho_T\}$. 

 For any $\mathbf{\xi}=(\xi_1^c,\xi_1^p,\ldots,\xi_r^c,\xi_r^p)\in\big(\mathcal{M}_1(\mathcal{Z})\big)^{2r}$, define the rate matrices 
\begin{align*}
A^{j,c}_{\xi}=\left(\lambda^c_{z,z'}(\xi_j^c,\xi_j^p)\right)_{(z,z')\in\mathcal{Z}\times\mathcal{Z}}\qquad\text{and}\quad A^{j,p}_{\xi}=\left(\lambda^p_{z,z'}(\xi_j^c,\xi_1^p\ldots,\xi^p_r)\right)_{(z,z')\in\mathcal{Z}\times\mathcal{Z}},
\end{align*}
where $\lambda^c_{z,z}(\xi_j^c,\xi_j^p)=-\sum_{z'\neq z}\lambda^c_{z,z'}(\xi_j^c,\xi_j^p)$ and $\lambda^p_{z,z}(\xi_j^c,\xi_1^p,\ldots,\xi_r^p)=-\sum_{z'\neq z}\lambda^p_{z,z'}(\xi_j^c,\xi_1^p,\ldots,\xi_r^p)$. Assuming that the initial condition $\mu^N(0)=\nu^N$ converges weakly to some $\nu\in (\mathcal{M}_1(\mathcal{Z}))^{2r}$ then, from the law of large numbers (cf. \cite[Coro. 3.1]{Daw+Sid+zha2020}), one can deduce that, as $N\rightarrow\infty$, the sequence $(\mu^N,N\geq 1)$ converges weakly towards the solution $\mu$ of the following McKean-Vlasov system of equations
\begin{equation}
\begin{split}
\left\{ \begin{array}{lcl}\dot{\mu}_j^c(t)=A^{j,c^*}_{\mu (t)}\mu_j^c(t), & &  \\
 \dot{\mu}_j^p(t)=A^{j,p^*}_{\mu(t)}\mu_j^p(t),  & & \\
 \mu_j^c(0)=\nu_j^c,\mu_j^p(0)=\nu_j^p, & & \\
 1\leq j\leq r, & &
 \end{array}\right.
\end{split}
\label{McKean-Vlas-syst}
\end{equation}
 where $A^*$ is the adjunct/transpose of the matrix $A$ and $\dot{\mu}(t)= \frac{\partial}{\partial t} \mu(t)$. Note that the Lipschitz property of the functions $\lambda^c_{z,z'}$ and $\lambda^p_{z,z'}$ guarantees that $(\ref{McKean-Vlas-syst})$ is well-posed. 
 \\
\indent Denote by $\tau$  the log-Laplace transform of the centered Poisson distribution with parameter $1$ given by $\tau(u)=e^u-u-1$, and let $\tau^*$ be its Legendre transform defined by 
\begin{equation}
\begin{split}
\tau^*(u)=\left\{\begin{tabular}{lll}
$(u+1) \log (u+1)-u$& \text{if} & $u>-1$, \\
$1$               & \text{if} & $u=-1$,\\
$+\infty$         & \text{if}  & $u<-1$. \\
\end{tabular}
\right.
\end{split}
\label{Legen-tran}
\end{equation}

Define, for any $\theta\in\mathcal{M}(\mathcal{Z})$,
\begin{align*}
|||\theta|||^{j,c}_{\mu(t)}&=\sup_{\Phi:\mathcal{Z}\rightarrow\mathbb{R}} \bigg\{ \sum_{z\in\mathcal{Z}}\theta(z)\cdot\Phi(z)-\sum_{z':(z,z')\in\mathcal{E}} \tau \big(\Phi(z')- \Phi(z)\big)\cdot\mu_j^c(t)(z)\cdot\lambda^c_{zz'}\big(\mu_j^c(t),\mu_j^p(t)\big) \bigg\},\\
|||\theta|||^{j,p}_{\mu(t)}&=\sup_{\Phi:\mathcal{Z}\rightarrow\mathbb{R}} \bigg\{ \sum_{z\in\mathcal{Z}}\theta(z)\cdot\Phi(z)-\sum_{z':(z,z')\in\mathcal{E}} \tau \big(\Phi(z')- \Phi(z)\big)\cdot\mu_j^p(t)(z)\cdot\lambda^p_{zz'}\big(\mu_j^c(t),\mu_1^p(t),\ldots,\mu_r^p(t)\big) \bigg\},
\end{align*}
and let us introduce, for each $\nu\in(\mathcal{M}_1(\mathcal{Z}))^{2r}$ and according to \cite[eqn.~(4.9)]{Daw+Gart87}, the functional $S(\mu|\nu)$ defined from  $(\mathcal{D}([0, T ], \mathcal{M}_1 (Z)))^{2r}$ to $[0,\infty]$ by setting
\begin{equation}
\begin{split}
S_{[0,T]}(\mu|\nu)= \sum_{j=1}^r\bigg[\alpha_jp_j^c \int_0^T ||| \dot{\mu}_j^c(t)-A_{\mu_j^c(t)}^{*}\mu_j^c(t)|||_{\mu(t)}dt  +\alpha_jp_j^p\int_0^T ||| \dot{\mu}_j^p(t)-A_{\mu_j^p(t)}^{*}\mu_j^p(t)|||_{\mu(t)}dt \bigg]
\label{rate-emp-proc}
\end{split}
\end{equation}
 if $\mu(0)=\nu$ and $\mu_j^c,\mu_j^p$ are absolutely continuous in the sense of \cite[Def.~4.1]{Daw+Gart87} for all $1\leq j\leq r$, and $S_{\nu}(\mu)=+\infty$ otherwise. 

\subsection{Large deviations over finite time intervals}
We recall here the large deviations principle for the sequence of probability measures $(p^N_{\nu_N},N\geq 1)$ over finite time intervals.  

\begin{theo}
Suppose that $\nu_N\rightarrow\nu$ weakly. The sequence of probability measures $(p^N_{\nu_N},N\geq 1)$ obeys a large deviations principle in the space $\big(\mathcal{D}([0,T],\mathcal{M}_1(\mathcal{Z}))\big)^{2r}$, with speed $N$, and rate function $S_{[0,T]}(\mu|\nu)$ given by  $(\ref{rate-emp-proc})$. \\
\indent Moreover, if a path $\mu\in\big(\mathcal{D}([0,T],\mathcal{M}_1(\mathcal{Z}))\big)^{2r}$ satisfies $S_{[0,T]}(\mu|\nu)<\infty$, then, for $1\leq j\leq r$, $\mu_j^c$ and $\mu_j^p$  are absolutely continuous and there exist rate matrices $L_{j,c}(t)=(l_{z,z'}^{j,c}(t), t \in [0, T ], (z, z')\in\mathcal{E})$ and $L_{j,p}(t)=(l_{z,z'}^{j,p}(t), t \in [0, T ], (z, z')\in\mathcal{E})$ such that
\begin{equation}
\begin{split}
\left\{ \begin{array}{l}
\dot{\mu}_j^c(t)={L_{j,c}(t)}^{*}\mu_j^c(t),\\
\dot{\mu}_j^p(t)={L_{j,p}(t)}^{*}\mu_j^p(t),\\
1\leq j \leq r, t\in [0,T],
\end{array}\right.
\label{no-Mcke-sys}
\end{split}
\end{equation}
 and the good rate function $S_{[0,T]}(\mu|\nu)$ is given by
\begin{equation}
\begin{split}
\sum_{j=1}^r\bigg[&\alpha_jp_j^c \int_0^T \bigg( \sum_{(z,z')\in\mathcal{E}}(\mu_j^c(t)(z))\lambda^c_{z,z'}\left(\mu_j^c(t),\mu_j^p(t)\right) \tau^*\bigg(\frac{l_{z,z'}^{j,c}(t)}{\lambda^c_{z,z'}\left(\mu_j^c(t),\mu_j^p(t)\right)}-1\bigg) \bigg)dt \\
 &+\alpha_jp_j^p\int_0^T  \bigg( \sum_{(z,z')\in\mathcal{E}}(\mu_j^p(t)(z))\lambda^p_{z,z'}\left(\mu_j^c(t),\mu_1^p(t),\ldots,\mu_r^p(t)\right) \tau^*\bigg(\frac{l_{z,z'}^{j,p}(t)}{\lambda^p_{z,z'}\left(\mu_j^c(t),\mu_1^p(t),\ldots,\mu_r^p(t)\right)}-1\bigg) \bigg) dt \bigg].
\label{rate-emp-proc-2}
\end{split}
\end{equation}
Conversely, if a path $\mu\in\big(\mathcal{D}([0,T],\mathcal{M}_1(\mathcal{Z}))\big)^{2r}$ satisfies the following: $\mu$ is absolutely continuous, $\mu(0)= \nu$, and there exist time-varying rate matrices $L_{j,c}(t)$ and $L_{j,p}(t)$ for $1\leq j\leq r$ such that $\mu=(\mu_j^c,\mu_j^p, 1\leq j\leq r)$ satisfies $(\ref{no-Mcke-sys})$, then the good rate function $S_{[0,T ]}(\mu|\nu)$ evaluated at $\mu$ is given by $(\ref{rate-emp-proc-2})$. 
\label{large-dev-emp-proc}
\end{theo}

\proof See \cite[Th. 4.2]{Daw+Sid+zha2020}. \carre

\begin{rem}
\begin{itemize}
\item The action functional $S$ characterizes the difficulty of the passage of $\mu^N$ near $\mu$ in the time interval $[0, T]$. Indeed, according to Theorem \ref{large-dev-emp-proc}, the probability of such a passage behaves like $\exp(- NS_{[0,T]}(\mu|\nu))$ as $N\rightarrow+\infty$.
\item Observe from $(\ref{rate-emp-proc-2})$ that if the rate function $S_{[0,T ]} (\mu|\nu)= 0$ then $\mu$ must be the solution to the McKean-Vlasov system $(\ref{McKean-Vlas-syst})$ with initial condition $\mu(0)= \nu$. 
\end{itemize}
\label{rem-ldp-theo}
\end{rem}

The following result states the uniform large deviations principle for $(p^N_{\nu_N},N\geq 1)$ with respect to the initial condition $\nu$ over compact sets. 
\begin{coro}
For any compact set $K\subset\big(\mathcal{M}_1(\mathcal{Z}))\big)^{2r}$, any closed set $F\subset\big(\mathcal{D}([0,T],\mathcal{M}_1(\mathcal{Z}))\big)^{2r}$, and any open set $G\subset\big(\mathcal{D}([0,T],\mathcal{M}_1(\mathcal{Z}))\big)^{2r}$, we have
\begin{align}
\limsup_{N\rightarrow\infty}\frac{1}{N} \log \sup_{\nu\in K} p_{\nu}^{N} \big({\mu^N\in F}\big) \leq - \inf_{\nu\in K} \inf_{\substack{ \mu\in F}} S_{[0,T]} (\mu|\nu),
\end{align}
and,
\begin{align}
\liminf_{N\rightarrow\infty}\frac{1}{N} \log \inf_{\nu\in K} p_{\nu}^{N} \big({\mu^N\in G}\big) \geq - \sup_{\nu\in K}\inf_{\substack{ \mu\in G}} S_{[0,T]} (\mu|\nu).
\end{align}
\label{coro-ldp-emp-proc}
\end{coro}

\proof See \cite[Coro. 5.6.15]{Dem+Zeit2010}. \carre

\subsection{Large deviations at initial and terminal times}

Fix $T> 0$. Recall that $p_{\nu_N}^N$ is the distribution of the empirical process $\mu^N\in(\mathcal{D}([0,T],\mathcal{M}_1(\mathcal{Z})))^{2r}$ with initial conditions given by the empirical vector $\mu^N(0)=\nu_N$. Let $\mu^N(T)$ be the empirical measure at time $T$ and let $p_{\nu_N,T}^N$ be its distribution. The next result states the large deviations principle for the sequence $(p_{\nu_N,T}^N,N\geq 1)$.

\begin{theo}
Suppose that $\nu_N\rightarrow\nu$ weakly. The sequence of probability measures $(p^N_{\nu_N,T},N\geq 1)$ obeys a large deviations principle in the space $\big(\mathcal{M}_1(\mathcal{Z})\big)^{2r}$ with speed $N$ and good rate function 
\begin{align}
S_T(\xi|\nu)=\inf\{S_{[0,T]}(\mu|\nu):\mu(0)=\nu,\mu(T)=\xi\}.
\label{T-rate}
\end{align}
Moreover, $S_T(\xi|\nu)$ is bounded for all $\xi,\nu\in(\mathcal{M}_1(\mathcal{Z}))^{2r}$ and its infimum is achieved.
\label{LDP-term-time}
 \end{theo}

\proof Recall that the space $\big(\mathcal{D}([0,T],\mathcal{M}_1(\mathcal{Z}))\big)^{2r}$ is equipped with the product topology induced by the product metric $\rho^{2r}_T=\max\{\rho_T,\cdots,\rho_T\}$ where  
\begin{align*}
\rho_T (\mu,\nu)=\sup_{0\leq t\leq T}\rho_0 (\mu_t,\nu_t),\quad\mu,\nu\in\mathcal{D}([0,T],\mathcal{M}_1(\mathcal{Z})).
\end{align*}
 Thus, it is easy to see that the application $\mu\in\big(\mathcal{D}([0,T],\mathcal{M}_1(\mathcal{Z}))\big)^{2r}\rightarrow\mu(T)\in\big(\mathcal{M}_1(\mathcal{Z})\big)^{2r}$ is continuous. Therefore, a simple application of the contraction principle (cf. \cite[Th. 4.2.1]{Dem+Zeit2010}) proves the validity of the large deviations principle. In addition, since $S_T(\xi|\nu)$ is a good rate function, its infimum is achieved over closed sets. Finally, the boundedness follows from Lemma \ref{T-rate-prop} below.   \carre
\\
\\
We next state some useful technical results.  

\begin{lem}
The following statements hold.
\begin{enumerate}
\item There exists a constant $C_1(T)<+\infty$ such that for any $\xi,\nu\in\mathcal{M}_1(\mathcal{Z})^{2r}$, there is a piecewise linear and continuous path $\mu$, with $\mu(0)=\nu$ and $\mu(T)=\xi$, having constant velocity in each linear segment and satisfies $S_{[0,T ]} (\mu|\nu) \leq C_1 (T )$.
\item For any $\xi,\nu\in(\mathcal{M}_1(\mathcal{Z}))^{2r}$, we have $S_T(\xi|\nu)\leq C_1(T)$. 
\item There exists a constant $C_3$ such that, for each $\varepsilon>0$, there exists a $\delta\in(0,\varepsilon)$ such that: $\rho^{2r}_0(\nu,\xi)<\delta$ implies that $S_{\varepsilon}(\xi|\nu)\leq C_3\varepsilon$. 
\end{enumerate}
\label{T-rate-prop}
\end{lem}

\proof See Appendix \ref{T-rate-prop-proof}. \carre

\begin{lem}
Let, for $1\leq j\leq r$, $L_{j,c}(t)$ and $L_{j,p}(t)$ be rate matrices such that the solution $\mu: [0, T ]\rightarrow (\mathcal{M}_1(\mathcal{Z}))^{2r}$ to the system
\begin{align*}
\left\{ \begin{array}{l}
\dot{\mu}_j^c(t)={L_{j,c}(t)}^{*}\mu_j^c(t),\\
\dot{\mu}_j^p(t)={L_{j,p}(t)}^{*}\mu_j^p(t),\\
1\leq j \leq r, t\in [0,T],
\end{array}\right.
\end{align*}
with $\mu(0)= \nu$ satisfies $S_{[0,T]} (\mu|\nu)<+\infty$. Then there exists a constant $\mathfrak{K} < +\infty$ such that
\begin{align*}
S_{[0,T ]} (\mu|\nu)\geq \sum_{j=1}^r\bigg[\alpha_jp_j^c \int_0^T& \bigg( \sum_{(z,z')\in\mathcal{E}}(\mu_j^c(t)(z)) l_{z,z'}^{j,c}(t) \bigg)dt+\alpha_jp_j^p\int_0^T  \bigg( \sum_{(z,z')\in\mathcal{E}}(\mu_j^p(t)(z))l_{z,z'}^{j,p}(t)  \bigg) dt \bigg]-\mathfrak{K}T
\end{align*}
\label{unif-cont-lem-1}
\end{lem}

\proof See Appendix \ref{unif-cont-lem-1-proof}. \carre

\begin{lem}
Let, for $1\leq j\leq r$, $L_{j,c}(t)$ and $L_{j,p}(t)$ be rate matrices such that the solution  $\mu: [0, T ]\rightarrow (\mathcal{M}_1(\mathcal{Z}))^{2r}$ to the system
\begin{align*}
\left\{ \begin{array}{l}
\dot{\mu}_j^c(t)={L_{j,c}(t)}^{*}\mu_j^c(t),\\
\dot{\mu}_j^p(t)={L_{j,p}(t)}^{*}\mu_j^p(t),\\
1\leq j \leq r, t\in [0,T],
\end{array}\right.
\end{align*}
with $\mu(0)= \nu$ has $S_{[0,T]} (\mu|\nu)<+\infty$ and $\mu(T)= \xi$. Let $0 < \beta < +\infty$ be a time scaling and $T'= T /\beta$. Consider the path $\{\tilde{\mu}(t)= \mu(\beta t)| t \in [0, T']\}$ having $\tilde{\mu}(t)=\nu$ and $\tilde{\mu}(T)=\xi$.
Then $\tilde{\mu}$ satisfies 
\begin{align*}
\left\{ \begin{array}{l}
\dot{\tilde{\mu}}_j^c(t)={\tilde{L}_{j,c}(t)}^{*}\tilde{\mu}_j^c(t),\\
\dot{\tilde{\mu}}_j^p(t)={\tilde{L}_{j,p}(t)}^{*}\tilde{\mu}_j^p(t),\\
1\leq j \leq r, t\in [0,T],
\end{array}\right.
\end{align*}
where $\tilde{L}_{j,c}(t) = \beta L_{j,c}(\beta t)$ and $\tilde{L}_{j,p}(t) = \beta L_{j,p}(\beta t)$ for $1\leq j\leq r$. Furthermore, the scaled path $\tilde{\mu}:[0, T']\rightarrow(\mathcal{M}_1(\mathcal{Z}))^{2r}$ satisfies
\begin{equation}
\begin{split}
S_{[0,T']} (\tilde{\mu}|\nu) \leq S_{[0,T]} (\mu|\nu)+\frac{|1-\beta|}{\beta} CT|\mathcal{E}|+ |\log\beta|\sum_{j=1}^r\bigg[ &\alpha_jp_j^c \int_0^{T} \bigg( \sum_{(z,z')\in\mathcal{E}}(\mu_j^c( t)(z))l_{z,z'}^{j,c}(t)\bigg)dt \\
 &+\alpha_jp_j^p \int_0^{T} \bigg( \sum_{(z,z')\in\mathcal{E}}(\mu_j^p( t)(z))l_{z,z'}^{j,p}(t)\bigg)dt \bigg].
\end{split}
\label{bound-time-scal}
\end{equation}
\label{unif-cont-lem-2}
\end{lem}

\proof See Appendix \ref{unif-cont-lem-2-proof}. \carre

\begin{lem}
The mapping $(\nu, \xi)\rightarrow S_T(\xi|\nu)$ is uniformly continuous.
\label{unif-cont-map}
\end{lem}

\proof This is a mild generalization of \cite[Lem. 3.3]{Bork+Sund2012}. Fix $T>0$, $0<\varepsilon<T/4$, and let $\delta\in(0,\varepsilon)$. From Lemma $\ref{T-rate-prop}$ there exists a constant $C_3$ such that, for any $\nu,\xi$ with $\rho_0^{2r}(\nu,\xi)<\delta$, we have $S_{\varepsilon}(\xi|\nu)\leq \varepsilon C_3$. Let $(\nu_1,\xi_1),(\nu_2,\xi_2)$ be two points in the space $(\mathcal{M}_1(\mathcal{Z}))^{2r}\times(\mathcal{M}_1(\mathcal{Z}))^{2r}$ equipped with the metric $\bar{\rho}$ defined  by
\begin{align*}
\bar{\rho}((\nu_1,\xi_1),(\nu_2,\xi_2))=\max\{\rho^{2r}_0(\nu_1,\nu_2),\rho^{2r}_0(\xi_1,\xi_2)\},
\end{align*}
and suppose that $\bar{\rho}((\nu_1,\xi_1),(\nu_2,\xi_2))<\delta$. Thus it suffices to show that $S_T(\xi_1|\nu_1)$ and $S_T(\xi_2|\nu_2)$ are close to each other. First note that $\rho^{2r}_0(\nu_1,\nu_2)<\delta$ and $\rho^{2r}_0(\xi_1,\xi_2)<\delta$. Therefore, from Lemma $\ref{T-rate-prop}$ there exists a path going from $\nu_1$ to $\nu_2$ with cost $S_{\varepsilon}(\nu_2|\nu_1)\leq C_3\varepsilon$ in time $\varepsilon$, and from $\xi_2$ to $\xi_1$ with cost $S_{\varepsilon}(\xi_1|\xi_2)\leq C_3\varepsilon$ in time $\varepsilon$. Moreover, denote by $\mu$ the minimum cost path from $\nu_2$ to $\xi_2$ in time $T$ with cost $S_T(\xi_2|\nu_2)$ and then consider the path from $\nu_1$ to $\xi_1$ as follows:
\begin{itemize}
\item Traverse the path from $\nu_1$ to $\nu_2$ in time $[0, \varepsilon]$ with cost at most $C_3\varepsilon$.
\item Consider the  path $\tilde{\mu}:[0,T-2\varepsilon]\rightarrow (\mathcal{M}_1(\mathcal{Z}))^{2r}$ given by $\tilde{\mu}(t)= \mu(\alpha t)$ with $\alpha= T /(T-2\varepsilon)$ where $\mu$ is the optimal $[0, T ]$-path $\mu$ from $\nu_2$ to $\xi_2$. Therefore, travel from $\nu_2$ to $\xi_2$ in the duration $[\varepsilon, T-\varepsilon]$ along the path $\tilde{\mu}$.
\item  Traverse the path from $\xi_2$ to $\xi_1$ in time $[0, \varepsilon]$ with cost at most $C_3\varepsilon$. 
\end{itemize}

Thus the minimum cost for traversal from $\nu_1$ to $\xi_1$ is at most the sum of the previous paths. Hence by Lemmas $\ref{unif-cont-lem-1}$ and $\ref{unif-cont-lem-2}$ one obtains  
\begin{align*}
S_T (\xi_1|\nu_1) \leq 2C_3 \varepsilon+S_{[0,T]} (\xi_2|\nu_2)+2 C\varepsilon|\mathcal{E}|+ \bigg(\log\frac{T}{T-2\varepsilon}\bigg) (S_{[0,T ]} (\xi_2|\nu_2)+\mathfrak{K}T).
\end{align*}

Note that $\log u \leq u- 1$ for $u > 0$. Thus using again Lemma $\ref{T-rate-prop}$ leads to 
\begin{align*}
S_T (\xi_1|\nu_1) \leq 2C_3 \varepsilon+S_{[0,T]} (\xi_2|\nu_2)+2 C\varepsilon|\mathcal{E}|+ \bigg(\frac{T}{T-2\varepsilon}-1\bigg) (C_1(T)+\mathfrak{K}T).
\end{align*}
Now, obbserve that
\begin{align*}
\varepsilon <T/4\Rightarrow \frac{T}{T-2\varepsilon}-1=\frac{T\varepsilon}{T-2\varepsilon}\leq\frac{4\varepsilon}{T},
\end{align*}
from which we deduce that 
\begin{align*}
S_T (\xi_1|\nu_1) \leq S_{[0,T]} (\xi_2|\nu_2)+\varepsilon C_4(T),
\end{align*}
with $C_4(T)=2C_3+2 C|\mathcal{E}|+\frac{4}{T}(C_1(T)+\mathfrak{K}T)$. Finally, reversing the roles of $(\nu_1, \xi_1 )$ and $(\nu_2, \xi_2 )$ gives to us
\begin{align*}
|S_T (\xi_1|\nu_1)- S_T(\xi_2|\nu_2)| \leq C_4 (T)\varepsilon,
\end{align*}
which concludes the proof. \carre
\\
\\
Denote by $\wp^N_0$ the law of the initial empirical measure vector $\mu^N(0)=\nu_N$, and let $\wp^N_{0,T}$ denote the joint law of $(\mu^N(0), \mu^N (T))$. We next establish the large deviations principle for the sequence $(\wp^N_{0,T}, N \geq 1)$. 

\begin{theo}
Suppose that the sequence $(\wp^N_{0}, N \geq 1)$ satisfies the large deviations principle with speed $N$ and good rate function $s: (\mathcal{M}_1 (\mathcal{Z}))^{2r}\rightarrow [0, +\infty]$. Then the sequence of joint laws $(\wp^N_{0,T}, N \geq 1)$ satisfies the large deviation principle with speed $N$ and good rate function
\begin{align}
S_{0,T} (\xi,\nu)= s(\nu)+ S_T (\xi|\nu).
\label{ldp-rate-init-term}
\end{align}
\label{ldp-init-term}
\end{theo}

\proof Here $\wp^N_{0,T}$ is the joint distribution of $(\mu^N(0), \mu^N (T))$, $\wp^N_0$ is the distribution of $\mu^N(0)=\nu_N$ and $p_{\nu^N_N,T}$  is the conditional distribution of $\mu^N(T)$ given $\mu^N(0)$. Thus one can write
\begin{align*}
d\wp^N_{0,T}(\nu,\xi)=d\wp^N_0(\nu)\times dp^N_{\nu_N,T}(\xi).
\end{align*}

Since both the sequences $(\wp^N_{0,T},N\geq 1)$ and  $(p_{\nu^N_N,T}, N \geq 1)$ obey large deviations principle, one can apply \cite[Prop. 3.25]{Feng+Kurt2006} in order to derive the rate function corresponding to the large deviations principle of the sequence $((\mu^N(0), \mu^N (T)),N\geq 1)$ in the product space. To this end, we first verify that the conditions of application of \cite[Prop. 3.25]{Feng+Kurt2006} are satisfied. Let $\nu_N\rightarrow\nu$ weakly. By Theorem $\ref{LDP-term-time}$, the sequence of laws of the terminal measure $(p^N_{\nu_N,T}, N \geq 1)$ satisfies the large deviations principle with speed $N$ and good rate function $S_T(\xi|\nu)$. Moreover, recall that $p_{\nu_N,T}=P^N_{\nu_N}\circ\pi_T^{-1}$. From \cite[Lem. 4.6]{Daw+Sid+zha2020} we have that for any $\alpha>0$,
\begin{align}
\limsup_{N\rightarrow\infty}\frac{1}{N}\log\int_{\mathcal{M}_{1,\varphi (\mathcal{X})}\times\cdots\times\mathcal{M}_{1,\varphi (\mathcal{X})}}\exp\left\{ N\alpha |h|   \right\} dP_{\nu_N}^{0,N}<\infty.
\label{vara-cond-1}
\end{align}

 Then, using $(\ref{vara-cond-1})$ and the Radon-Nikodym derivative $(\ref{rad-nik})$, one can easily verify that, for any $f\in C_b ((\mathcal{M}_1(\mathcal{Z}))^{2r}$ and $\alpha>0$,
\begin{align*}
\limsup_{N\rightarrow\infty}\frac{1}{N}\log\int_{(\mathcal{M}_1(\mathcal{Z}))^{2r}}e^{N\alpha|f|}dp^N_{\nu_N,T}<\infty.
\end{align*}
Moreover, using Varadhan's lemma (see \cite[Prop. 2.5]{Leonard95}) we obtain, for every $f\in C_b ((\mathcal{M}_1(\mathcal{Z}))^{2r}$, 
\begin{align}
\lim_{N\rightarrow\infty}\frac{1}{N}\log\int_{(\mathcal{M}_1(\mathcal{Z}))^{2r}}e^{N f}dp^N_{\nu_N,T}=\sup_{\xi\in(\mathcal{M}_1(\mathcal{Z}))^{2r}}[f(\xi)-S_T(\xi|\nu)].
\label{conv-vara-loc}
\end{align}
Furthermore, by defining 
\begin{align}
\Lambda(f|\nu)=\sup_{\xi\in(\mathcal{M}_1(\mathcal{Z}))^{2r}}[f(\xi)-S_T(\xi|\nu)],
\label{Lambd-funct}
\end{align}
the Bryc's Inverse Varadhan Lemma (see \cite[Th. 4.4.2]{Dem+Zeit2010}) gives
\begin{align*}
S_T (\xi|\nu)=\sup_{f \in C_b ((\mathcal{M}_1(\mathcal{Z}))^{2r})} [f (\xi)-\Lambda(f |\nu)].
\end{align*}
One can observe that, for any $f \in C_b ((\mathcal{M}_1(\mathcal{Z}))^{2r}$, the mapping $\nu\in(\mathcal{M}_1(\mathcal{Z})^{2r})\rightarrow\Lambda(f|\nu)\in\mathbb{R}$ is continuous. Indeed, for any $f \in C_b ((\mathcal{M}_1 (\mathcal{Z}))^{2r})$, define the mapping
\begin{align*}
\eta: (\nu, \xi)\in(\mathcal{M}_1 (\mathcal{Z}))^{2r}\times(\mathcal{M}_1 (\mathcal{Z}))^{2r}\rightarrow [f (\xi)- S_T (\xi|\nu)]\in\mathbb{R}.
\end{align*} 
Since $f$ is continuous, and $S_T (\xi|\nu)$ is also continuous by Lemma \ref{unif-cont-map}, the mapping $\eta(\cdot,\cdot)$ is jointly continuous. Let $\nu_N\rightarrow\nu$ weakly and for each $N$, let $\xi_N$ denote a point where the supremum in $(\ref{Lambd-funct})$, corresponding to $\nu_N$, is attained. This is consistent since the space $(\mathcal{M}_1(\mathcal{Z}))^{2r}$ is compact. Thus  $\Lambda(f|\nu_N)=\eta(\nu_N, \xi_N)$ for each $N$. Moreover, by the same compactness argument, the sequence $((\nu_N,\xi_N),N\geq 1)$ has a convergent subsequence that converge to $(\nu,\xi)$, for some $\xi\in\mathcal{M}_1 (\mathcal{Z}))^{2r}$. Let $((\nu_{N'},\xi_{N'}),N'\geq 1)$ denote this subsequence. Therefore $(\nu_{N'},\xi_{N'})\rightarrow (\nu, \xi)$ as $N'\rightarrow +\infty$. By the continuity of $\eta$ we have, as $N'\rightarrow +\infty$,
\begin{align*}
\Lambda(f |\nu_{N'})= \eta(\nu_{N'}, \xi_{N'})\rightarrow \eta(\nu, \xi).
\end{align*}

From the definition of $\xi_N$, one can observe that, for any $\xi'$, $\eta(\nu_{N'}, \xi') \leq \eta(\nu_{N'}, \xi_{N'})$. Hence by the continuity of $\eta$
\begin{align*}
\eta(\nu, \xi')=\lim_{N'\rightarrow\infty}  \eta(\nu_{N'}, \xi') \leq \limsup_{N'\rightarrow\infty} \eta(\nu_{N'}, \xi_{N'})= \eta(\nu, \xi).
\end{align*}

Thus  $\Lambda(f |\nu)=\eta(\nu, \xi)$, from which we deduce that $\Lambda(f |\nu)$ is continuous since $\eta(\nu, \xi)$ is jointly continuous. Thus, the mapping $\nu\in(\mathcal{M}_1(\mathcal{Z})^{2r})\rightarrow\Lambda(f|\nu)\in\mathbb{R}$ is indeed continuous. 

Now, notice that when there are $N$ particles in the system, the initial empirical measure $\mu^N(0)=\nu_N$ takes values in the product space $(\mathcal{M}^N_1(\mathcal{Z}))^{2r}$ where 
\begin{align*}
\mathcal{M}^N_1(\mathcal{Z})=\bigg\{\frac{1}{N}\sum_{n=1}^N\delta_{a_n}\bigg|a^N=(a_n,1\leq n\leq N)\in\mathcal{Z}^N\bigg\},
\end{align*}
which is a compact subset of $\mathcal{M}_1(\mathcal{Z})$. Using Theorem \ref{large-dev-emp-proc} and the continuity of $\nu\in(\mathcal{M}_1(\mathcal{Z})^{2r})\rightarrow\Lambda(f|\nu)\in\mathbb{R}$, one can prove that the convergence in $(\ref{conv-vara-loc})$ is uniform for $\nu$ over compact subsets of $(\mathcal{M}_1(\mathcal{Z}))^{2r}$, namely,
\begin{align*}
\lim_{N\rightarrow+\infty}\sup_{\nu\in(\mathcal{M}^N_1(\mathcal{Z}))^{2r}}\bigg|\frac{1}{N} \log\int_{(\mathcal{M}_1(\mathcal{Z}))^{2r}}e^{N\alpha f}dp^N_{\nu_N,T}-\Lambda(f|\nu)\bigg|= 0,
\end{align*}
for any $f \in C_b ((\mathcal{M}_1(\mathcal{Z}))^{2r}$ (see \cite[Lem. 8.2]{Bork+Sund2012} for a detailed proof). Finally one has to verify that the sequence $(\wp^N_{0,T}, N \geq 1)$ is exponentially tight. Since the sequence $(\mu^N(0), \mu^N (T))$ takes values in a product space, it is enough to verify that each marginal distribution is exponentially
tight. This is in fact straightforward since, by the compactness of the space $(\mathcal{M}_1(\mathcal{Z}))^{2r}$, the exponential tightness of the sequences $(\wp^N_0,N\geq 1)$ and $(\wp^N_{0,T},N\geq 1)$ is implied by the goodness of the rate functions $S_T(\mu|\nu)$ and $s(\nu)$ (see \cite[Exercise 1.2.19]{Dem+Zeit2010}).

We are now in position to apply \cite[Prop. 3.25]{Feng+Kurt2006}. Indeed, we have shown that the sequence $(\wp^N_{0,T}, N \geq 1)$ is exponentially tight. Moreover, the convergence in convergence in $(\ref{conv-vara-loc})$ is uniform for $\nu$ in compact subsets of $(\mathcal{M}_1(\mathcal{Z}))^{2r}$. Furthermore, the function $\Lambda(f|\nu)$ is continuous in $\nu$. Hence, since $(\wp^N_0,N\geq 1)$ satisfies the large deviations principle with good rate function $s(\nu)$, then $(\wp^N_{0,T}, N \geq 1)$ satisfies the large deviations principle with good rate function given by $(\ref{ldp-rate-init-term})$. This concludes the proof. \carre

\section{Large deviations of the invariant measure} 
 \label{LDP-IM-sec}
By Assumption \ref{ass-prin}, the graph $(\mathcal{Z},\mathcal{E})$ of the allowed transitions is irreducible. Moreover, the state space $\mathcal{Z}$ is finite, therefore, for each fixed total number of particles $N$, there exists a unique invariant measure for the Markov process $X^N=(X_n, X_m,n\in C_j^c, m\in C_j^p,1\leq j\leq r)$. Hence there is a unique invariant measure, denoted by $\wp^N$, for the $(\mathcal{M}_1 (\mathcal{Z}))^{2r}$-valued Markov process $\mu^N$. The goal of this section is to investigate the large deviations properties of the sequence $(\wp^N, N\geq 1)$ under two separate scenarios. First, we consider the case where the limiting Mckean-Vlasov system $(\ref{McKean-Vlas-syst})$ has a unique globally asymptotically stable equilibrium. Then we treat the general case with multiple $\omega$-limit sets.

\subsection{Unique globally asymptotically stable equilibrium} 

We first establish the large deviations principle for the invariant measure $(\wp^N, N\geq 1)$ in the case where the limiting McKean-Vlasov system has a unique globally asymptotically stable equilibrium $\xi_0$. 
 
\begin{theo}
Suppose that Assumption $\ref{ass-prin}$ holds true. Moreover, suppose that the McKean-Vlasov system $(\ref{McKean-Vlas-syst})$  has a unique globally asymptotically stable equilibrium $\xi_0$. Then the sequence $(\wp^{N},N\geq 1)$ satisfies a large deviations principle with speed $N$ and a good rate function $s$ given by   
  \begin{equation}
\begin{split}
s(\xi)=\inf_{\hat{\mu}}\sum_{j=1}^r&\bigg[\alpha_jp_j^c \int_0^{+\infty} \bigg( \sum_{(z,z')\in\mathcal{E}}(\hat{\mu}_j^c(t)(z))\lambda^c_{z,z'}\left(\hat{\mu}_j^c(t),\hat{\mu}_j^p(t)\right) \tau^*\bigg(\frac{\hat{l}_{z,z'}^{j,c}(t)}{\lambda^c_{z,z'}\left(\hat{\mu}_j^c(t),\hat{\mu}_j^p(t)\right)}-1\bigg) \bigg)dt \\
 &+\alpha_jp_j^p\int_0^{+\infty}  \bigg( \sum_{(z,z')\in\mathcal{E}}(\hat{\mu}_j^p(t)(z))\lambda^p_{z,z'}\left(\hat{\mu}_j^c(t),\hat{\mu}_1^p(t),\ldots,\hat{\mu}_r^p(t)\right) \\
 &\qquad\qquad\qquad\qquad\times\tau^*\bigg(\frac{\hat{l}_{z,z'}^{j,p}(t)}{\lambda^p_{z,z'}\left(\hat{\mu}_j^c(t),\hat{\mu}_1^p(t),\ldots,\hat{\mu}_r^p(t)\right)}-1\bigg) \bigg) dt \bigg],
 \label{infi-rate}
\end{split}
\end{equation}
where the infimum is over all the paths $\hat{\mu}$ that are solutions to the reversed-time system 
\begin{equation}
\begin{split}
\left\{ \begin{array}{l}
\dot{\hat{\mu}}_j^{c}(t)=-{\hat{L}_{j,c}(t)}^{*}\hat{\mu}_j^{c}(t),\\
\dot{\hat{\mu}}_j^{p}(t)=-{\hat{L}_{j,p}(t)}^{*}\hat{\mu}_j^{p}(t),\\
1\leq j\leq r,
\end{array}\right.
\end{split}
\end{equation}
for some family of rate matrices $\hat{L}_{j,c}$ and $\hat{L}_{j,p}$, with initial condition $\mu(0)=\xi$, terminal condition $\lim_{t\rightarrow\infty}\mu(t)=\xi_0$, and $\mu(t)\in(\mathcal{M}_1(\mathcal{Z}))^{2r}$ for all $t\geq 0$.
\label{LDP-invariant-GAS}
\end{theo} 
 
The rest of this section is dedicated to the proof of Theorem $\ref{LDP-invariant-GAS}$. The proof is based on the control-theoretic approach introduced in \cite{Bis+Bork2011} for small noise diffusions and in \cite{Bork+Sund2012} for finite-state mean-field systems on complete graphs. We proceed through several lemmas. We start by establishing a subsequential large deviations principle. 
 \begin{lem}
 Let $(N \geq 1)$ be a sequence of natural numbers going to $+\infty$. Then there exists a subsequence $(N_k, k \geq 1)$ such that $(\wp^{N_k},N_k\geq 1)$ satisfies a large deviations principle with speed $N_k$ and a good rate function $s$ verifying
 \begin{align}
 s(\xi)=\inf_{\nu\in(\mathcal{M}_1 (\mathcal{Z}))^{2r}}[s(\nu)+ S_T (\xi|\nu)]\quad\text{for every $T>0$}.
 \label{rate-inv-meas-1}
 \end{align}
 Moreover, $s \geq 0$ and there exists some $\nu^*\in(\mathcal{M}_1(\mathcal{Z}))^{2r}$ such that $s(\nu^*)=0$. 
 \label{lem-rate-inv-meas-1}
 \end{lem}
 
 \proof The space $\mathcal{Z}$ being compact since finite then, the space $\mathcal{M}_1(\mathcal{Z})$ is also compact. Hence, the product space $(\mathcal{M}_1 (\mathcal{Z}))^{2r}$ is also compact. Moreover, the space $(\mathcal{M}_1 (\mathcal{Z}))^{2r}$ endowed with the product metric $\rho^{2r}_0(\cdot,\cdot)$ is a metric space and thus it has a countable basis. Therefore, by \cite[Lem.~4.1.23]{Dem+Zeit2010}, there exists a sequence $N_k\rightarrow 0$ such that $(\wp^{N_k},N_k\geq 1)$ satisfies the weak large deviations principle in $(\mathcal{M}_1 (\mathcal{Z}))^{2r}$. Moreover, since $(\mathcal{M}_1 (\mathcal{Z}))^{2r}$ is compact, $(\wp^{N_k},N_k\geq 1)$ satisfies the strong large deviations principle with a good rate function $s:(\mathcal{M}_1 (\mathcal{Z}))^{2r}\rightarrow[0,\infty]$ and speed $N_k$ (see, e.g. \cite[Lem. 1.2.18]{Dem+Zeit2010}).

 Fix an arbitrary $T>0$. Recall that $\wp^N_{0}$ is the probability measure of the initial empirical measure vector $\mu^N(0)$. Set $\wp^N_{0}=\wp^N$. Then, by Theorem \ref{ldp-init-term}, the sequence of joint laws $(\wp^{N_k}_{0,T}, N_k \geq 1)$ satisfies the large deviations principle along the subsequence $(N_k, k \geq 1)$ with speed $N_k$ and good rate function $S_{0,T} (\xi,\nu)= s(\nu)+ S_T (\xi|\nu)$. Using the continuity of the projection and the contraction principle \cite[Th.~4.2.1]{Dem+Zeit2010}, the sequence of the terminal probability distributions $(\wp^{N_k}_{T}, N_k \geq 1)$ satisfies a large deviations principle with the good rate function 
 \begin{align*}
 \inf_{\nu\in(\mathcal{M}_1 (\mathcal{Z}))^{2r}}S_{0,T} (\xi,\nu)= \inf_{\nu\in(\mathcal{M}_1 (\mathcal{Z}))^{2r}}[s(\nu)+ S_T (\xi|\nu)].
 \end{align*}
 
 Since $\wp^N$ is the invariant measure and that $\wp^N_{0}=\wp^N$, one has $\wp^N_{T}=\wp^N$. Thus, by the uniqueness of the rate function, we deduce that
 \begin{align*}
   \inf_{\nu\in(\mathcal{M}_1 (\mathcal{Z}))^{2r}}S_{0,T} (\xi,\nu)=s(\xi).
 \end{align*}
  
 Finally, since every rate function is nonnegative, we have $s\geq 0$. In addition, $s$ being a good rate function, it attains its minimum and thus, there exists some $\nu^*$ such that $s(\nu^*)=0$. This concludes the proof. \carre
\\
\\
Notice that equation $(\ref{rate-inv-meas-1})$ has multiple solutions. The trivial $s\equiv 0$ is one of them. To see this, take the McKean-Vlasov path $\mu$ of duration $T$ given by $(\ref{McKean-Vlas-syst})$ starting at some initial condition $\nu$ and ending at $\xi$, then $S_T(\xi|\nu)=0$ and the right hand side of $(\ref{rate-inv-meas-1})$ vanishes when $s\equiv 0$. Therefore, one shall identify more conditions satisfied by the rate function $s$. 

Observe from $(\ref{T-rate})$ that $(\ref{rate-inv-meas-1})$ is equivalent to 
 \begin{equation}
 \begin{split}
 s(\xi)&=\inf_{\nu\in(\mathcal{M}_1 (\mathcal{Z}))^{2r}}\{s(\nu)+ \inf_{\mu} [S_{[0,T]}(\mu|\mu(0)):\mu (0)=\nu, \mu (T)=\xi ]\}\quad\text{for every $T> 0$}.\\
 &=\inf_{\mu|\mu(T)=\xi}\{s(\mu(0))+  S_{[0,T]}(\mu|\mu(0))\},\quad\text{for every $T> 0$}.
\end{split} 
 \end{equation}
Hence, one can see $(\ref{rate-inv-meas-1})$ as a Bellman equation associated with an optimal control problem, with $s$ being the corresponding \textit{value function} of a minimization problem over path space with paths defined on $[0,mT]$, for $m\geq 1$, and ending at $\xi$. Therefore, one must determine the optimal control problem for which equation $(\ref{rate-inv-meas-1})$ is the dynamic programming equation and $\xi$ is the terminal condition. Since the terminal condition is fixed, we shall define translated and reversed-time paths that start from $\xi$ at time $0$. In particular, for $m\geq 1$ and a path $\mu\in\big(\mathcal{D}([0,mT],\mathcal{M}_1(\mathcal{Z}))\big)^{2r}$ satisfying 
\begin{equation}
\begin{split}
\left\{ \begin{array}{l}
\dot{\mu}_j^c(t)={L_{j,c}(t)}^{*}\mu_j^c(t),\\
\dot{\mu}_j^p(t)={L_{j,p}(t)}^{*}\mu_j^p(t),\\
1\leq j\leq r, t\in [0,mT],
\end{array}\right.
\end{split}
\end{equation}
  with $L_{j,c}(t)$ (resp. $L_{j,p}(t)$) being the rate matrix associated with the time-varying rates $(l_{z,z'}^{j,c}(t), (z, z')\in\mathcal{E})$ (resp. $(l_{z,z'}^{j,p}(t), (z, z')\in\mathcal{E})$), the corresponding translated and reversed path $\hat{\mu}$ satisfies the following reversed-time dynamical system 
\begin{equation}
\begin{split}
\left\{ \begin{array}{l}
\dot{\hat{\mu}}_j^c(t)=-{\hat{L}_{j,c}(t)}^{*}\hat{\mu}_j^c(t),\\
\dot{\hat{\mu}}_j^p(t)=-{\hat{L}_{j,p}(t)}^{*}\hat{\mu}_j^p(t),\\
1\leq j\leq r, t\in [0,mT],
\end{array}\right.
\end{split}
\end{equation}
with initial condition $\hat{\mu}(0)=\xi$, where
\begin{align*}
\hat{\mu}(t)&=\mu(mT-t),\\
\hat{L}_{j,c}(t)&=L_{j,c}(mT-t),\\
\hat{L}_{j,p}(t)&=L_{j,p}(mT-t),
\end{align*}
for $t\in [0,mT]$. The latter quantities are defined where time flows in the opposite direction regarding the direction under the McKean-Vlasov dynamics $(\ref{McKean-Vlas-syst})$. Thus, one can consider the set of rate matrices $\hat{L} (t)=(\hat{L}_{j,c}(t),\hat{L}_{j,p}(t), 1\leq j\leq r)$ as the control at time $t$ when the state is $\hat{\mu}(t)$. Define
\begin{equation}
\begin{split}
r(\hat{\mu}(t),\hat{L}(t),t)=&\sum_{j=1}^r\bigg[\alpha_jp_j^c  \bigg( \sum_{(z,z')\in\mathcal{E}}(\hat{\mu}_j^c(t)(z))\lambda^c_{z,z'}\left(\hat{\mu}_j^c(t),\hat{\mu}_j^p(t)\right) \tau^*\bigg(\frac{\hat{l}_{z,z'}^{j,c}(t)}{\lambda^c_{z,z'}\left(\hat{\mu}_j^c(t),\hat{\mu}_j^p(t)\right)}-1\bigg) \bigg)\\
 &+\alpha_jp_j^p  \bigg( \sum_{(z,z')\in\mathcal{E}}(\hat{\mu}_j^p(t)(z))\lambda^p_{z,z'}\left(\hat{\mu}_j^c(t),\hat{\mu}_1^p(t),\ldots,\hat{\mu}_r^p(t)\right)\tau^*\bigg(\frac{\hat{l}_{z,z'}^{j,p}(t)}{\lambda^p_{z,z'}\left(\hat{\mu}_j^c(t),\hat{\mu}_1^p(t),\ldots,\hat{\mu}_r^p(t)\right)}-1\bigg) \bigg) \bigg],
 \label{cost-funct}
\end{split}
\end{equation}
thence the total cost over $[0,mT]$ is 
\begin{align*}
\int_{0}^{mT}r(\hat{\mu}(t),\hat{L}(t),t)dt,
\end{align*}
which is simply $S_{[0,mT]}(\mu|\mu(0))$ given by $(\ref{rate-emp-proc-2})$, as can be verified using a simple change of variable $t\rightarrow mT-t$. Thence, the equation $(\ref{rate-inv-meas-1})$ is equivalent to 
 \begin{align}
 s(\xi)=\inf_{\hat{L}}\{s(\hat{\mu}(mT))+ \int_{0}^{mT}r(\hat{\mu}(t),\hat{L}(t),t)dt\}\quad\text{for every $m\geq 1$},
 \label{rate-inv-meas-1-equiv}
 \end{align}
which is the Bellman equation associated with the optimal control problem 
\begin{align*}
\min_{\hat{L}}\int_{0}^{mT}r(\hat{\mu}(t),\hat{L}(t),t)dt, \quad m\geq 1,
\end{align*}
subject to 
\begin{equation}
\begin{split}
\left\{ \begin{array}{l}
\dot{\hat{\mu}}_j^c(t)=-{\hat{L}_{j,c}(t)}^{*}\hat{\mu}_j^c(t),\\
\dot{\hat{\mu}}_j^p(t)=-{\hat{L}_{j,p}(t)}^{*}\hat{\mu}_j^p(t),\\
1\leq j \leq r, t\in [0,mT], m\geq 1,
\end{array}\right.
\end{split}
\end{equation}
and $\hat{\mu}(0)=\xi$. The next lemma establishes the existence of one optimal path $\hat{\mu}$ of infinite duration starting at $\xi$.

\begin{lem}
For each $\xi\in(\mathcal{M}_1 (\mathcal{Z}))^{2r}$, there exists a path $\hat{\mu}=(\hat{\mu}_1^c,\hat{\mu}_1^p,\ldots,\hat{\mu}_r^c,\hat{\mu}_r^p): [0, +\infty)\rightarrow (\mathcal{M}_1(\mathcal{Z}))^{2r}$ and families of rate matrices  $\hat{L}_{j,c}(t)$ and $\hat{L}_{j,p}(t)$ defined on $[0,\infty)$ such that, for $t\in[0,\infty)$,
\begin{equation}
\begin{split}
\left\{ \begin{array}{l}
\dot{\hat{\mu}}_j^c(t)=-{\hat{L}_{j,c}(t)}^{*}\hat{\mu}_j^c(t),\\
\dot{\hat{\mu}}_j^p(t)=-{\hat{L}_{j,p}(t)}^{*}\hat{\mu}_j^p(t),\\
1\leq j \leq r,
\label{no-Mck-hat}
\end{array}\right.
\end{split}
\end{equation}
with initial condition $\hat{\mu}(0)=\xi$, and 
\begin{equation}
\begin{split}
s(\xi)=s(\hat{\mu}(mT))+ \int_{0}^{mT}r(\hat{\mu}(t),\hat{L}(t),t)dt,\quad\mbox{for all $m\geq 1$}.
 \label{rate-inv-meas-2}
\end{split}
\end{equation}
\label{opt-path-lem}
\end{lem}

\proof Let $\hat{\mu}=(\hat{\mu}_1^c,\hat{\mu}_1^p,\ldots,\hat{\mu}_r^c,\hat{\mu}_r^p):[0,\infty)\rightarrow(\mathcal{M}_1(\mathcal{Z}))^{2r}$ be an infinite duration path, and let $\psi_m$ be its restriction to the finite interval $[0,mT]$ given by
\begin{align*}
\hat{\mu}\rightarrow\psi_m\hat{\mu}=\hat{\mu}^{(m)}:[0,mT]\rightarrow(\mathcal{M}_1(\mathcal{Z}))^{2r}.
\end{align*} 
 Equip the space of infinite duration paths with the following metric
 \begin{align*}
 \rho_{\infty}(\hat{\mu},\hat{\nu})=\sum_{m=1}^{\infty}2^{-m}(\rho^{2r}_{mT}(\psi_m\hat{\mu},\psi_m\hat{\nu})\wedge 1).
 \end{align*} 
 One can easily observe that the restriction $\psi_m$ is continuous for all $m$. Moreover, define the reversed restriction by
\begin{align*}
\mu^{(m)} (t)=\hat{\mu}^{(m)} (mT-t), \quad t \in [0,mT].
\end{align*}
 For any $B\in [0,\infty)$, consider the sets
\begin{align*}
\Gamma_{\infty}(B)&=\left\{\hat{\mu}:[0,\infty)\rightarrow(\mathcal{M}_1(\mathcal{Z}))^{2r}\big|\sup_{m\geq 1}S_{[0,mT]}(\mu^{(m)}|\mu^{(m)}(0))\leq B\right\},\\
\Gamma_{m}(B)&=\left\{\hat{\eta}:[0,mT]\rightarrow(\mathcal{M}_1(\mathcal{Z}))^{2r}\big|S_{[0,mT]}(\eta|\eta(0))\leq B\right\}.
\end{align*}  
Thus $\Gamma_{\infty}(B)$ is the set of infinite duration paths $\hat{\mu}$ with the corresponding restrictions $\mu^{(m)}$ on time intervals $[0,mT]$ having the costs bounded by $B$, for any $m\geq 1$. On the other hand,  the set $\Gamma_{m}(B)$ comports the paths $\hat{\eta}$ of duration $[0,mT]$ with corresponding reversed path $\eta(t)=\hat{\eta}(mT-t)$ having cost bounded by $B$. We next prove that the set $\Gamma_{\infty}(B)$ is compact for any $B\in [0,\infty)$. 

First, observe that $\Gamma_{\infty}(B)$ is a subset of a metric space, thus it is enough to prove that it is sequentially compact. Take an infinite sequence $(\hat{\mu}_n,n\geq 1)\subset\Gamma_{\infty}(B)$ and the corresponding restriction sequence  $(\psi_m\hat{\mu}_n,n\geq 1)\subset\Gamma_{m}(B)$. Note that the sets $\Gamma_{m}(B)$, for $m\geq 1$, are compact. Indeed, $S_{[0,mT]}$ being a good rate function, the corresponding level set $\left\{\eta:[0,mT]\rightarrow(\mathcal{M}_1(\mathcal{Z}))^{2r}\big|S_{[0,mT]}(\eta|\eta(0))\leq B\right\}$ is compact and then, since $\hat{\eta}(t)=\eta(mT-t)$ for each $t\in [0,mT]$, we have that $\Gamma_{m}(B)$ is compact since it is a continuous image of a compact. Therefore, one can find an infinite subset $\mathbb{V}_1\subset\mathbb{N}$ such that $(\psi_1\hat{\mu}_n,n\in\mathbb{V}_1)\subset\Gamma_{1}(B)$ converges. Moreover, one can take a further subsequence represented by the infinite subset $\mathbb{V}_2 \subset \mathbb{V}_1$ such that $(\psi_2\hat{\mu}_n,n\in\mathbb{V}_2)\subset\Gamma_{2}(B)$ converges. We can continue this procedure for all $m\geq 1$ and eventually take the subsequence along the diagonal. This subsequence converges for every interval $[0,mT]$. Take $\hat{\mu}(t)$ its pointwise limit for every $t$. The restriction $\psi_m\hat{\mu}\in\Gamma_m(B)$, and the corresponding $\mu^{(m)}$ satisfies $S_{[0,mT]}(\mu^{(m)}|\mu^{(m)}(0))\leq B$ and thus, $\hat{\mu}\in\Gamma_{\infty}(B)$. This guarantees that $\Gamma_{\infty}(B)$ is sequentially compact and thus compact.  
 
Fix $\xi\in\left(\mathcal{M}_1(\mathcal{Z})\right)^{2r}$. From Lemma $\ref{lem-rate-inv-meas-1}$, there exists $\nu^*$ such that $s(\nu^*)=0$ and then, by Lemma $\ref{T-rate-prop}$,  
 \begin{align}
 s(\xi)\leq S_T(\xi|\nu^*)\leq C_1(T).
\label{bound-s-C}
\end{align}  
Take in the sequel $B=C_1(T)$. Starting from any location $\nu\in(\mathcal{M}_1(\mathcal{Z}))^{2r}$, the minimum cost $S_{mT}(\xi|\nu)$ of transporting the system from $\nu$ to $\xi$ in $mT$ units of time is also bounded by $C_1(T)$. Indeed, consider the path consisting of traversing the McKean-Vlasov path with initial condition $\nu$ in $(m-1)T$ units of time, the corresponding cost is zero. Then, go to $\xi$ in $T$ units of time. The cost of the last traversal is bounded by $B=C_1(T)$ as stated by Lemma $\ref{T-rate-prop}$. Now, if we consider the translated and reversed time paths $\hat{\mu}$ that start from $\xi$ at time $0$, they have a cost at most $B$ and stay within $\left(\mathcal{M}_1(\mathcal{Z})\right)^{2r}$ for the duration $[0,mT]$. Let
\begin{align*}
\Gamma_m^*=\bigcup_{\nu\in(\mathcal{M}_1(\mathcal{Z}))^{2r}}\big\{\hat{\mu}:[0,mT]\rightarrow(\mathcal{M}_1(\mathcal{Z}))^{2r}\big|\mu(t)=\hat{\mu}(mT-t),\mu(0)=\nu,\mu(mT)=\xi,S_{[0,mT]}(\mu|\nu)=S_{mT}(\xi|\nu)\big\}
\end{align*}  
be the subset of $\Gamma_m$ consisting of the collection of all minimum cost reversed and translated paths $\hat{\mu}$ on $[0,mT]$, starting from $\xi$ to every location in $(\mathcal{M}_1(\mathcal{Z}))^{2r}$. Thanks to the arguments above, the set $\Gamma_m^*$ is not empty for all $m\geq 1$ and the minimum cost $S_{mT}(\xi|\nu)$ is at most $B$ for every $\nu$. 

Next, we prove that $\Gamma_m^*$ is compact. To this end, and since $\Gamma_m^*$ is a subset of a compact set, it suffices to show that it is closed. Let $\hat{\mu}$ be a point of the closure of $\Gamma_m^*$, then one can find  a sequence $(\hat{\mu}^{(k)}, k \geq 1)\subset\Gamma_m^*$ such that $\lim_{k\rightarrow\infty}\hat{\mu}^{(k)}= \hat{\mu}$. By definition of $\Gamma_m^*$ we have $\hat{\mu}(0)= \xi$. In addition, set $\hat{\mu}(mT)=\nu$ for some $\nu\in(\mathcal{M}_1(\mathcal{Z}))^{2r}$. Consider the corresponding translated and reversed paths $(\mu^{(k)}, k \geq 1)$ and $\mu$. Therefore, using the lower semi-continuity property of the good rate function $S_{[0,mT]}(\cdot|\nu)$ we obtain
\begin{align*}
S_{[0,mT]}(\mu|\nu)&\leq\liminf_{k\rightarrow\infty}S_{[0,mT]}(\mu^{(k)}|\mu^{(k)}(0))\\
                   &=\liminf_{k\rightarrow\infty}S_{mT}(\xi|\mu^{(k)}(0))\\
                   &=S_{mT}(\xi|\nu),
\end{align*}
where the first equality follows from the definition of $\Gamma^*_m$, and the last equality is found using the continuity of the rate function $S_{mT}$ in its arguments (cf. Lemma $\ref{unif-cont-map}$). On the other hand, by definition, $S_{mT}(\xi|\nu)$ is the minimum cost of transporting the system from $\nu$ to $\xi$ in $mT$ units of time, we then deduce that the path $\mu$  must satisfy  $S_{[0,mT]}(\mu|\nu)=S_{mT}(\xi|\nu)$ and thus $\hat{\mu}\in\Gamma_m^*$. Hence $\Gamma_m^*$ is closed and thus compact.  
  
We are now in the position to conclude the proof. Given that $\Gamma_m^*$ is nonempty and closed, by the continuity of the restriction $\psi_m$, the image set $\psi_m^{-1}\Gamma_m^*$ is nonempty and closed. Moreover, since $\Gamma_{\infty}(B)$ is compact, the intersection $\psi_m^{-1}\Gamma_m^*\cap \Gamma_{\infty}(B)$ is also compact. Furthermore, note that $\psi_{m+1}^{-1}\Gamma_{m+1}^*\subset\psi_m^{-1}\Gamma_m^*$. Indeed, take an optimal path $\hat{\mu}\in\psi_{m+1}^{-1}\Gamma_{m+1}^*$ that starts from $\xi$ and passes by  $\nu$ at time $(m+1)T$ and by $\nu'$ at time $mT$ then, its restriction to the interval $[0,mT]$ is necessarily optimal (if not $\hat{\mu}$ would not be optimal on $[0,(m+1)T]$). Therefore $\hat{\mu}\in\psi_{m}^{-1}\Gamma_{m}^*$ and thus $(\psi_m^{-1}\Gamma_m^*,m\geq 1)$ is a nested decreasing sequence of subsets. Hence, $(\psi_m^{-1}\Gamma_m^*\cap \Gamma_{\infty},m\geq 1)$ is in turn a nested decreasing sequence of nonempty, compact and closed subsets and thus, by Cantor's intersection theorem, its intersection is not empty. Take $\hat{\mu}\in\bigcap_{m\geq 1}(\psi_m^{-1}\Gamma_m^*\cap\Gamma_{\infty})$ and let $\eta^{(m)}(t)=\hat{\mu}(mT+T-t)$, for $m\geq 1$  and $t\in [0,T]$, be its reversal restriction on the interval $ [mT,mT+T]$, and let $\eta^{(0)}(t)=\hat{\mu}(T-t)$ $t\in [0,T]$ be the reversal restriction of $\hat{\mu}$ on $[0,T]$. Thus $S_{[0,T]}(\eta^{(m)}|\eta^{(m)}(0))\leq B<\infty$ for all $m\geq 0$. Thence, from Theorem \ref{large-dev-emp-proc}, there exist families of rate matrices  $(L^{(m)}_{j,c}(t),t\in [0,T])$ and $(L^{(m)}_{j,p}(t),t\in [0,T])$ such that $\eta^{(m)}=(\eta^{(m),c}_1,\eta^{(m),p}_1,\ldots,\eta^{(m),c}_r,\eta^{(m),p}_r)$ satisfies $(\ref{no-Mcke-sys})$ on $[0,T]$ with initial condition $\eta^{(m)} (0)$. Define for all $m\geq 0$, $\hat{L}_{j,c}(mT+t)=L^{(m)}_{j,c}(T-t)$ and $\hat{L}_{j,p}(mT+t)=L_{j,p}^{(m)}(T-t)$, for $t\in[0,T]$. Thus, the rate matrices $\hat{L}_{j,c}$ and $\hat{L}_{j,p}$ are defined on $[0,\infty)$ such that $\hat{\mu}$ satisfies $(\ref{no-Mck-hat})$ on $[0,\infty)$ with initial condition $\hat{\mu}(0)=\xi$. The equality in $(\ref{rate-inv-meas-2})$ follows since the path $\hat{\mu}$ satisfies $(\ref{rate-inv-meas-1-equiv})$ on any interval $[0,mT]$ by the principle of optimality. This concludes the proof of the lemma. \carre 
\\
\\ 
 Next, we show that the optimal path given in Lemma $\ref{opt-path-lem}$ must end up in an invariant set of the dynamics given by the time-reversed McKean-Vlasov system.

 \begin{lem}
 Let 
 \begin{align*}
 \Omega=\bigcap_{t>0}\overline{\{\hat{\mu}(t'),t'>t\}}
 \end{align*}
 be the $\omega$-limit set of the path $\hat{\mu}$ given in Lemma $\ref{opt-path-lem}$. Then, $\Omega$ is contained in an $\omega$-limit set of the reversed-time McKean-Vlasov system
\begin{equation}
\begin{split}
\left\{ \begin{array}{lcl}\dot{\hat{\mu}}_j^c(t)=-A^{j,c^*}_{\hat{\mu (t)}}\hat{\mu}_j^c(t), & &  \\
 \dot{\hat{\mu}}_j^p(t)=-A^{j,p^*}_{\hat{\mu}(t)}\hat{\mu}_j^p(t),  & & \\
 1\leq j\leq r, t\geq 0. & &
 \end{array}\right.
\end{split}
\label{McKean-Vlas-syst2}
\end{equation}
\label{omega-set-lem}
 \end{lem}
 
 \proof The path $\hat{\mu}$ given in Lemma $\ref{opt-path-lem}$ remains in $(\mathcal{M}_1(\mathcal{Z}))^{2r}$ and satisfies the dynamics $(\ref{no-Mck-hat})$ with initial condition $\hat{\mu}(0)=\xi$, for some $\xi\in (\mathcal{M}_1(\mathcal{Z}))^{2r} $. Moreover, notice that the integral term in $(\ref{rate-inv-meas-2})$ is increasing with $m$ since the integrand is nonnegative. Hence, since the equality is valid for all $m\geq 1$, the term $s(\hat{\mu}(mT))$ must decreases as $m$ increases. But we know from Lemma $\ref{T-rate-prop}$ and $(\ref{bound-s-C})$ that $s\in[0,C_1(t)]$ and then, there exists $s^*$ such that $s\rightarrow s^*$ as $m\rightarrow\infty$. 
 
Take a subsequence $(\hat{\mu}(m_kT),m_k\geq 1)$ of $(\hat{\mu}(mT),m\geq 1)$ and let $\xi'$ be its limit. Moreover, denote by $\nu$ the limit of the subsequence  $(\hat{\mu}(m_kT+T),m_k\geq 1)$, and consider the path of duration $T$ given by
\begin{align*}
\mu^{(m_k)} (t) = \hat{\mu}(m_k T+T-t), \quad t \in [0, T]. 
\end{align*}
Thus the following convergence 
 \begin{align*}
 (\mu^{(m_k)} (0),\mu^{(m_k)} (T))\rightarrow(\nu,\xi')
\end{align*}  
 holds as $k\rightarrow\infty$. In addition, both $s(\mu^{(m_k)} (0))$ and $s(\mu^{(m_k)} (T)))$ go to $s^*$ as $k\rightarrow\infty$ since $(\ref{rate-inv-meas-2})$ is valid for all $m\geq 1$. Therefore, taking limits in $(\ref{rate-inv-meas-2})$ as $k\rightarrow\infty$ with both $\hat{\mu}(m_k T+T)$ and $\hat{\mu}(m_k T)$ we deduce that
 \begin{align*}
\limsup_{k\rightarrow\infty}
\int_{m_kT}^{m_kT+T} r(\hat{\mu}(t),\hat{L}(t),t)dt=0,
\end{align*} 
  which, by a simple change of variable $t\rightarrow m_kT+T-t$, gives in turn 
  \begin{align*}
\limsup_{k\rightarrow\infty} S_{[0,T]}(\mu^{(m_k)}|\mu^{(m_k)}(0))=0.
\end{align*}  
 By the nonegativity of the rate function $S_T$, one further finds  
  \begin{align*}
\lim_{k\rightarrow\infty} S_{T}(\mu^{(m_k)}(T)|\mu^{(m_k)}(0))=0.
\end{align*}

 In addition, Lemma $\ref{unif-cont-map}$ tells us that the rate function $S_T$ is uniformly continuous in both its arguments. Thence, we deduce that $S_T(\xi'|\nu)=0$. This means that the path, say $\mu$, that goes from $\nu$ to $\xi'$ in $T$ units of time with no cost $S_{[0,T]}(\mu|\nu)= S_T(\xi'|\nu)=0$, is necessarily the McKean-Vlasov path that satisfies $(\ref{McKean-Vlas-syst})$ on $[0,T]$, with initial condition $\mu(0)=\nu$ and terminal condition $\mu(T)=\xi'$ (see Remark \ref{rem-ldp-theo}). Hence, the reversed-time path $\bar{\mu}(t)=\mu(T-t)$ satisfies $(\ref{McKean-Vlas-syst2})$ for $t\in [0,T]$ with initial condition $\bar{\mu}(0)=\xi'$. This proves that $\Omega$, the $\omega$-limit set of $\hat{\mu}$, is contained in an $\omega$-limit set of the reversed-time McKean-Vlasov dynamics $(\ref{McKean-Vlas-syst2})$. The lemma is proved. \carre
 \\
 \\
 The following result proves that the rate function $s$ vanishes at the equilibrium point $\xi_0$. 
 
 \begin{lem}
 If the McKean-Vlasov system $(\ref{McKean-Vlas-syst})$ has a unique globally asymptotically stable equilibrium $\xi_0$, then $s(\xi_0)=0$.
  \label{s-zero-lem}
 \end{lem}

 \proof Consider the Mckean-Vlasov dynamics $(\ref{McKean-Vlas-syst})$ with initial condition $\mu(0)=\nu^*$ satisfying $s(\nu^*)=0$ (see Lemma $\ref{lem-rate-inv-meas-1}$). Since $\xi_0$ is the unique asymptotically stable equilibrium, $\lim_{t\rightarrow\infty}\mu(t)=\xi_0$. Moreover, the McKean-Vlasov path $\mu$ has no cost, i.e., $S_T(\mu(T)|\nu^*)=0$ for each $T>0$. Therefore, from equation $(\ref{rate-inv-meas-1})$, one obtains 
 \begin{align*}
 s(\mu(T))\leq s(\nu^*)+S_T(\mu(T)|\nu^*)=s(\nu^*).
 \end{align*}
 
 Taking the limit as $T\rightarrow\infty$, and using the lower semicontinuity of $s$, one gets 
 \begin{align*}
 0 \leq s(\xi_0) \leq \liminf_{T\rightarrow\infty} s(\mu(T)) \leq s(\nu^*) = 0,
 \end{align*}
hence $s(\xi_0)=0$. \carre 
\\
\\
Finally, the next result gives the unique form of the rate function.   

\begin{lem}
 If the McKean-Vlasov system $(\ref{McKean-Vlas-syst})$ has a unique globally asymptotically stable equilibrium $\xi_0$, then the solution to $(\ref{rate-inv-meas-1})$ and $(\ref{rate-inv-meas-2})$ is unique and is given by $(\ref{infi-rate})$. 
 \label{rate-uniq-GAS}
\end{lem} 
 \proof  From Lemma $\ref{omega-set-lem}$, the $\omega$-limit set $\Omega$ of the path $\hat{\mu}$ given in Lemma $\ref{opt-path-lem}$ is contained in an $\omega$-limit set of the reversed-time McKean-Vlasov dynamics $(\ref{McKean-Vlas-syst2})$, which is both positively and negatively invariant. Moreover, the path $\hat{\mu}$ stays within $(\mathcal{M}_1(\mathcal{Z}))^{2r}$ (see the proof of Lemma $\ref{opt-path-lem}$) and thus, $\Omega\subset(\mathcal{M}_1(\mathcal{Z}))^{2r}$. Since $(\mathcal{M}_1(\mathcal{Z}))^{2r}$ is compact and that $\Omega$ is closed being an $\omega$-limit set then, it is non-empty, compact and connected. But by the assumption, the forward McKean-Vlasov system $(\ref{McKean-Vlas-syst})$ possesses a unique globally asymptotically stable equilibrium $\xi_0$, thus $\xi_0$ is the unique non-empty, compact, connected set in $(\mathcal{M}_1(\mathcal{Z}))^{2r}$ which is both positively and negatively invariant for the dynamics in $(\ref{McKean-Vlas-syst2})$. Thus we necessarily have $\Omega=\{\xi_0\}$. Therefore, letting $m\rightarrow\infty$ in $(\ref{rate-inv-meas-2})$ we obtain, by lower semicontinuity of the rate function $s$, 
\begin{equation}
\begin{split}
s(\xi)\geq s(\xi_0)+ \int_0^{\infty}r(\hat{\mu}(t),\hat{L}(t),t)dt.
\end{split}
\end{equation} 
But, from Lemma $\ref{s-zero-lem}$, we have $s(\xi_0)=0$, then 
 \begin{equation}
\begin{split}
s(\xi)\geq \int_0^{\infty}r(\hat{\mu}(t),\hat{L}(t),t)dt.
\end{split}
\end{equation}
Finally, applying $(\ref{rate-inv-meas-1})$ shows that $s(\xi)$ is upper bounded by the right hand side of $(\ref{infi-rate})$ and thus the equality holds. The uniqueness follows since the rate function is unique (\cite[Sect. 4.1.1]{Dem+Zeit2010}). \carre
 \\
 \\
  \textbf{Proof of Theorem \ref{LDP-invariant-GAS}.} We are now ready to conclude the proof of Theorem $\ref{LDP-invariant-GAS}$. Take an arbitrary sequence of positive numbers going to $+\infty$. From Lemma $\ref{lem-rate-inv-meas-1}$ there exists a subsequence $(N_k, k \geq 1)$ such that $(\wp^{N_k},N_k\geq 1)$ satisfies the large deviations principle with speed $N_k$ and a good rate function $s$ that satisfies equation $(\ref{rate-inv-meas-1})$. Moreover, from Lemma \ref{rate-uniq-GAS}, the rate function $s$ is uniquely determined by equation $(\ref{infi-rate})$. Therefore, for every sequence, there exists a further subsequence $(N_k , k \geq 1)$ such that   $(\wp^{N_k},N_k\geq 1)$ satisfies the large deviations principle with speed $N_k$ and the same good rate function $s$ given by $(\ref{infi-rate})$. Hence, since the rate function  $s$ corresponding to the subsequence $(\wp^{N_k},N_k\geq 1)$ satisfying a large deviations principle does not depend on this subsequence, one can conclude that the sequence of probability measures $(\wp^{N},N\geq 1)$ satisfies the large deviations principle with speed $N$ and good rate function $s$ defined by $(\ref{infi-rate})$. \carre

\subsection{Multiple $\omega$-limit sets}  
 
 We consider now the general case where the McKean-Vlasov system $(\ref{McKean-Vlas-syst})$ admits multiple $\omega$-limit sets. To this end, we rely on the classical Freidlin-Wentzell program (cf. \cite{Freid+Wentz2012}). The same strategy was adapted in \cite{Bork+Sund2012} for finite-state mean-field systems on complete graphs. We begin by introducing the central concepts. 
 
 First, recall the important notion of quasipotential $V: (\mathcal{M}_1(\mathcal{Z}))^{2r}\times(\mathcal{M}_1(\mathcal{Z}))^{2r}\rightarrow [0,\infty)$  defined, for any $\nu,\xi\in(\mathcal{M}_1(\mathcal{Z}))^{2r}$, by 
\begin{align}
V (\xi|\nu)= \inf\{S_{[0,T ]} (\mu|\nu): \mu(0)=\nu, \mu(T)=\xi, T > 0\}.
\label{Freid-Went-Pot}
\end{align}

Roughly speaking, $V (\xi|\nu )$ measures the difficulty for the empirical vector process to move from $\nu$ to $\xi$ in a finite time interval.   

\begin{lem}
The mapping $(\nu,\xi)\rightarrow V(\xi|\nu)$ is uniformly continuous.
\label{V-unif-cont}
\end{lem}

\proof Replacing the metric $\rho_0(\cdot,\cdot)$ by the product metric $\rho^r_0(\cdot,\cdot)$, the proof of \cite[Lem. 3.4]{Bork+Sund2012} holds verbatim. \carre 
\\
\\
Using the Freidlin-Wentzell quasipotential, we define the following equivalence relation on $(\mathcal{M}_1 (\mathcal{Z}))^{2r}$:
\begin{align}
\nu\sim\xi\quad\mbox{if}\quad V (\xi|\nu)= V (\nu|\xi)=0.
\label{equiv-rel-pot}
\end{align}

We make throughout the section the following Freidlin-Wentzell assumptions \cite[Chap.6, Sect.2]{Freid+Wentz2012}.  
\begin{ass}
In $(\mathcal{M}_1 (\mathcal{Z}))^{2r}$ there exist a finite number of compact sets $K_1,K_2,\ldots,K_l$ such that
\begin{enumerate}
\item For any two points $\nu_1$ and $\nu_2$ belonging to the same compact we have  $\nu_1\sim\nu_2$.
\item For each $i\neq j$, $\nu_1\in K_i$ and $\nu_2 \in K_j$ implies $\nu_1\nsim\nu_2$.
\item Every $\omega$-limit set of the McKean-Vlasov system $(\ref{McKean-Vlas-syst})$ lies completely in one of the compact sets $K_i$.
\end{enumerate}
\label{ass-mult-omeg}
\end{ass}
 
For any compacts $K_i$ and $K_j,i\neq j$, let us introduce
\begin{align}
V (K_i, K_j )= \inf\left\{S_{[0,T ]} (\mu|\mu(0)): \mu(0)\in K_i, \mu(T)\in K_j,T>0 \right\},
\label{V-comp}
\end{align}
which represents the minimum cost of going from the compact $K_i$ to $K_j$ in a finite time interval. Next, define the minimum cost of going from $K_i$ to $K_j$ without touching the other compact sets $K_k, k\neq i, j$ by
\begin{align}
\tilde{V} (K_i, K_j )= \inf\left\{S_{[0,T ]} (\mu|\mu(0)): \mu(0)\in K_i, \mu(T ) \in K_j, \mu(t) \notin\cup_{k \neq i,j}K_k\quad 
\mbox{for all $0\leq t \leq T$}, T > 0\right\},
\label{V-tilde-comp}
\end{align}
and use the convention that if there are no such paths, we set $\tilde{V} (K_i, K_j)=+\infty$. 
 
 \begin{defi}
 Let $L$ be a finite set and let $W\subset L$ be a subset of $L$. An oriented graph consisting of edges $(m,n)$ $(m\in L/W , n\in L, n\neq m)$ is called a $W$-graph if it satisfies the following conditions:
\begin{enumerate}
\item  Every point $m \in L/ W$ is the initial point of exactly one edge;
\item  There are no closed cycles in the graph.
\end{enumerate}
\end{defi}
The second condition above can be replaced by the following condition: for any point $m \in L/ W$, there exists a sequence of edges leading from it to some point $n\in W$. Denote by  $G(W )$ the set of $W$-graphs. Take $L=\{1, 2,\ldots, l\}$ the indices corresponding to the compact sets $K_1,K_2,\ldots, K_l$ given in Assumption \ref{ass-mult-omeg}, we define the following quantity
\begin{align}
W(K_i):=\min_{g\in G\{i\} }\sum_{(i,j)\in g}\tilde{V}(K_i,K_j)=\min_{g\in G\{i\} }\sum_{(i,j)\in g}V(K_i,K_j),
\label{W-func}
\end{align}
where $G\{i\}$ is the $W$-graph corresponding to $W=i$ with $i\in\{1,\ldots,l\}$. The second equality holds by \cite[Ch. 6, Lem 4.1]{Freid+Wentz2012}. We now state the main result of this section. 
 
\begin{theo}
Suppose that both Assumptions $\ref{ass-prin}$ and $\ref{ass-mult-omeg}$ hold true. Then the sequence $(\wp^{N},N\geq 1)$ satisfies the large deviations principle with speed $N$ and a good rate function $s$ given by 
  \begin{equation}
\begin{split}
s(\xi)=\inf_{1\leq l'\leq l}\inf_{\hat{\mu}}\bigg[s_{l'}+\sum_{j=1}^r&\bigg[\alpha_jp_j^c \int_0^{+\infty} \bigg( \sum_{(z,z')\in\mathcal{E}}(\hat{\mu}_j^c(t)(z))\lambda^c_{z,z'}\left(\hat{\mu}_j^c(t),\hat{\mu}_j^p(t)\right) \tau^*\bigg(\frac{\hat{l}_{z,z'}^{j,c}(t)}{\lambda^c_{z,z'}\left(\hat{\mu}_j^c(t),\hat{\mu}_j^p(t)\right)}-1\bigg) \bigg)dt \\
 &+\alpha_jp_j^p\int_0^{+\infty}  \bigg( \sum_{(z,z')\in\mathcal{E}}(\hat{\mu}_j^p(t)(z))\lambda^p_{z,z'}\left(\hat{\mu}_j^c(t),\hat{\mu}_1^p(t),\ldots,\hat{\mu}_r^p(t)\right) \\
 &\qquad\qquad\qquad\qquad\times\tau^*\bigg(\frac{\hat{l}_{z,z'}^{j,p}(t)}{\lambda^p_{z,z'}\left(\hat{\mu}_j^c(t),\hat{\mu}_1^p(t),\ldots,\hat{\mu}_r^p(t)\right)}-1\bigg) \bigg) dt \bigg]\bigg],
 \label{infi-rate-2}
\end{split}
\end{equation}
where $s_{l'}=W(K_{l'})-\min_{l'} W(K_{l'})$ with $W(K_{l'})$ given by $(\ref{W-func})$, and the second infimum is over all $\hat{\mu}$ that are solutions to the dynamical system 
\begin{equation}
\begin{split}
\left\{ \begin{array}{l}
\dot{\hat{\mu}}_j^{c}(t)=-{\hat{L}_{j,c}(t)}^{*}\hat{\mu}_j^{c}(t),\\
\dot{\hat{\mu}}_j^{p}(t)=-{\hat{L}_{j,p}(t)}^{*}\hat{\mu}_j^{p}(t),\\
1\leq j\leq r,
\end{array}\right.
\end{split}
\end{equation}
for some family of rate matrices $\hat{L}_{j,c}$ and $\hat{L}_{j,p}$, with initial condition $\mu(0)=\xi$, terminal condition $\lim_{t\rightarrow\infty}\mu(t)\in K_{l'}$, and $\mu(t)\in(\mathcal{M}_1(\mathcal{Z}))^{2r}$ for all $t\geq 0$.
\label{LDP-invariant-MEP}
\end{theo}

{\bf Proof of Thereom \ref{LDP-invariant-MEP}.} Again, the proof is split into several lemmas.

\begin{lem}
Assume Assumption $\ref{ass-mult-omeg}$ holds true. The rate function $s$ satisfies the following assertions:
\begin{itemize}
\item There exists $\xi_0$ in some compact $K_{i_0}$, with $1\leq i_0 \leq l$, that satisfies $s(\xi_0 )=0$;
\item The rate function $s$ is constant on each of the compacts $K_1,\ldots,K_l$, i.e., for all  $\xi\in K_i$, $s(\xi)=s_i$, for some nonnegative real numbers $s_1, s_2,\ldots, s_l$.
\end{itemize}
\label{rate-func-comp}
\end{lem} 

\proof Let $\mu$ be the McKean-Vlasov path given by $(\ref{McKean-Vlas-syst})$ and starting at some $\nu^*$, with $s(\nu^*)=0$. Note that the existence of such a point is guaranteed by Lemma $\ref{lem-rate-inv-meas-1}$. Take $\xi_0\in(\mathcal{M}_1 (\mathcal{Z}))^{2r}$ such that $\lim_{t\rightarrow}\mu(t)=\xi_0$. By Assumption $\ref{ass-mult-omeg}$, $\xi_0\in K_{i_0}$ for some $1\leq i_0\leq l$. Using the same arguments as in the proof of Lemma $\ref{s-zero-lem}$, the first statement  follows. 

Fix $1\leq i\leq l$ and let $\nu,\xi\in K_i$. Thus $\nu\sim\xi$. Hence, by $(\ref{Freid-Went-Pot})$ and $(\ref{equiv-rel-pot})$, for any $\varepsilon>0$, there exists $T>0$ and a path of length $T$ starting at $\nu$ and ending at $\xi$ such that $S_T(\xi|\nu)\leq\varepsilon$. Consequently, using $(\ref{rate-inv-meas-1})$, one obtains 
 \begin{align*}
 s(\xi)\leq s(\nu)+S_T(\xi|\nu)\leq s(\nu)+\varepsilon.
 \end{align*}
Reversing the roles of $\nu$ and $\xi$, one gets  
  \begin{align*}
 s(\nu)\leq s(\xi)+\varepsilon.
 \end{align*} 
Combining the two last inequalities leads to $|s(\nu)-s(\xi)|\leq\varepsilon$. Letting $\varepsilon\downarrow 0$ we obtain $s(\xi)=s(\nu)$, which proves the second statement. \carre

 \begin{lem}
 Assume Assumption $\ref{ass-mult-omeg}$ holds true. Then, the rate function $s$, the solution to the equations $(\ref{rate-inv-meas-1})$ and $(\ref{rate-inv-meas-2})$  is uniquely given by $(\ref{infi-rate-2})$. 
 \label{rate-uniq-MES}
 \end{lem}
 
 \proof Using the same arguments as in the proof of Lemma $\ref{rate-uniq-GAS}$ together with Assumption $\ref{ass-mult-omeg}$, one can conclude that the $\omega$-limit set $\Omega$ of the path $\hat{\mu}$ given in Lemma $\ref{opt-path-lem}$ is contained in some  compact $K_{l'}$, with $l'\in 1,\ldots,l$. Therefore, letting $m\uparrow\infty$ in equation $(\ref{rate-inv-meas-2})$, $s(\hat{\mu}(mT))\rightarrow s_{l'}$ for some $l'\in \{1,\ldots,l\}$ which gives that $s(\xi)$ is lower bounded by the right hand side of $(\ref{infi-rate-2})$. The upper bound is obtained by equation $(\ref{rate-inv-meas-1})$, and thus the proof follows. \carre

\begin{rem}
Using Lemma $\ref{rate-func-comp}$, Theorem $\ref{inv-meas-ineq}$, and the definition of the large deviations principle, it is easy to see that the values $s_1,\ldots,s_l$ of the rate function $s$ at the compact sets $K_i$ are given by $W(K_i)-\min_iW(K_i)$, with $W(K_i)$ defined by $(\ref{W-func})$. 
\end{rem} 
 
 \paragraph{Proof of Theorem $\ref{LDP-invariant-MEP}$} Take an arbitrary sequence of positive numbers going to $+\infty$. From Lemma $\ref{lem-rate-inv-meas-1}$ there exists a subsequence $(N_k, k \geq 1)$ such that $(\wp^{N_k},N_k\geq 1)$ satisfies the large deviations principle with speed $N_k$ and a good rate function $s$ that satisfies equation $(\ref{rate-inv-meas-1})$. Moreover, from Lemma $\ref{rate-uniq-MES}$, the rate function $s$ is uniquely determined by equation $(\ref{infi-rate-2})$ where $s_{l'}=W(K_{l'})-\min_{l'}W(K_{l'})$ for $1\leq l'\leq l$. Therefore, for every sequence, there exists a further subsequence $(N_k , k \geq 1)$ such that   $(\wp^{N_k},N_k\geq 1)$ satisfies the large deviations principle with speed $N_k$ and the same good rate function $s$ given by $(\ref{infi-rate-2})$. Hence, since the rate function  $s$ of the subsequence $(\wp^{N_k},N_k\geq 1)$ satisfying a large deviations principle does not depend on this subsequence, we conclude that the sequence of probability measures $(\wp^{N},N\geq 1)$ satisfies the large deviations principle with speed $N$ and good rate function $s$ defined by $(\ref{infi-rate-2})$. \carre

\section{Metastability and convergence to the invariant measure}
\label{Meta-sect}
We study in this section the {\textit metastable phenomena} which occur when the total number $ N $ of particles in the system as well as the time $ t $ are large. First, let us briefly summarize the main results of the previous sections and their consequences.

From the law of large numbers (cf. \cite[Coro. 3.1]{Daw+Sid+zha2020}), as $N\rightarrow\infty$ and for converging initial conditions $\mu^N(0)=\nu^N\rightarrow\nu$, the sequence $(\mu^N,N\geq 1)$ converges weakly and uniformly over any finite time interval $[0,T]$, towards the deterministic solution $\mu$ of the McKean-Vlasov system in $(\ref{McKean-Vlas-syst})$ with initial condition $\nu$. This suggests that, when $N\rightarrow\infty$ and over any finite time interval $[0,T]$, one can approximate the trajectories of the empirical vector process $\mu^N$ by the solution to the McKean-Vlasov system in $(\ref{McKean-Vlas-syst})$ with initial condition $\nu$. 

Moreover, given that the graph of allowed transitions $(\mathcal{Z},\mathcal{E})$ is irreducible, there is a unique invariant measure $\wp^N$ for $\mu^N$. Therefore, for each $N$, the distribution of $\mu^N$ converges towards $\wp^N$ as $t\rightarrow\infty$. Thence, the large $N$ behavior is described by the large deviations properties of the invariant measure $\wp^N$ established in Theorems $\ref{LDP-invariant-GAS}$ and $\ref{LDP-invariant-MEP}$. Thence, the natural question one might ask is whether or not we can interchange the $N\rightarrow\infty$ and $t\rightarrow\infty$ limits, namely 
\begin{align*}
\lim_{t\rightarrow\infty}\lim_{N\rightarrow\infty}\mu^N(t)\overset{?}{=} \lim_{N\rightarrow\infty}\lim_{t\rightarrow\infty}\mu^N(t).
\end{align*}

This classical question is in fact related to the limiting behavior of the McKean-Vlasov system $(\ref{McKean-Vlas-syst})$. In particular, two distinct cases must be considered. The most simple situation is when the McKean-Vlasov system $(\ref{McKean-Vlas-syst})$ has a unique globally asymptotically stable equilibrium $\xi_0\in (\mathcal{M}_1(\mathcal{Z}))^{2r}$. In this case, one can prove that the unique invariant measure $\wp^N$ of the empirical process vector $\mu^N$ converges towards the point mass $\delta_{\xi_0}$ as $N\rightarrow\infty$. Moreover, since $\wp^N$ is unique, it is independent of the initial condition. This leads to a justification of interchange of the limits $N\rightarrow\infty$ and $t\rightarrow\infty$. A detailed discussion of this scenario is given in \cite{Ben+Leboud2008}. 

The second and more complicated case, which we are interested in here, is when the McKean-Vlasov system $(\ref{McKean-Vlas-syst})$ has multiple $\omega$-limit sets depending on the initial condition. In this case, starting at a given $\nu\in (\mathcal{M}_1(\mathcal{Z}))^{2r}$, and letting $t\rightarrow\infty$, the solution to the McKean-Vlasov system in $(\ref{McKean-Vlas-syst})$ goes to an $\omega$-limit set corresponding to this initial condition. Moreover, recalling that the Birkhoff center of the solution $\mu$ to the McKean-Vlasov system $(\ref{McKean-Vlas-syst})$ is the closure of the set of recurrent points, i.e., the set of points $\nu\in(\mathcal{M}_1(\mathcal{Z}))^{2r}$ such that $\nu\in\omega(\nu)$, it is well known that the support of any limit point of $\wp^N$, as $N\rightarrow\infty$, is a compact subset of the Birkhoff center of $\mu$. See, e.g., \cite[Th. 3]{Ben+Leboud2008}. Therefore, most of the time, for large but finite $N$, the empirical vector process $\mu^N (t)$ remains close to the Birkhoff center of $\mu$. However, the difficulty here is that the Birkhoff center contains multiple $\omega$-limit sets, stable equilibrium, and/or limit cycles depending on the initial condition. Thence, \textit{metastable phenomena} are likely to arise. Here is an example.  

Let the initial condition $\mu^N(0)=\nu^N$ converge weakly to a given $\nu\in(\mathcal{M}_1(\mathcal{Z}))^{2r}$ as $N\rightarrow\infty$. Then on any finite time horizon, for a large but finite $N$, the empirical vector process $\mu^N$ would track the solution of the McKean-Vlasov equation $(\ref{McKean-Vlas-syst})$ starting at $\nu$ with a high probability. Therefore, as $t$ becomes large, $\mu^N$ would enter a neighborhood of the $\omega$-limit set of  $(\ref{McKean-Vlas-syst})$ corresponding to the initial condition $\nu$. However, since $N$ is finite, the process can exit the basin of attraction of this $\omega$-limit set, and likely remains in a neighborhood of another $\omega$-limit set for a large amount of time before transiting again to the next one, and so on. This is an example of metastability. The goal of this section is thus to study such phenomena. The main ingredient is the large deviations properties developed in Sections \ref{LDP-EP-sec} and \ref{LDP-IM-sec}.  

The metastable phenomena have first been studied for diffusion processes with a small noise parameter. The two main references are \cite[Ch. 6]{Freid+Wentz2012}, in which the authors studied these phenomena under the hypothesis summarized in Assumption \ref{ass-mult-omeg}, and \cite{Hwang+Sheu90} where slightly more general assumptions have been made. More recently, an extension to finite-state mean-field models on complete graphs was established in \cite{Yaso+Sund2019}. We propose in this section an extension of the aforementioned results to mean-field models with jumps on block-structured graphs detailed in Section \ref{model}. The main idea is to consider the empirical vector process $\mu^N(t)$ as a small noise perturbation of the deterministic solution $\mu$ to the McKean-Vlasov system $(\ref{McKean-Vlas-syst})$. Here, $N^{-1}$ plays the role of the small noise parameter $\varepsilon$ considered in \cite{Freid+Wentz2012,Hwang+Sheu90} in the sense that, as $N\rightarrow\infty$, $\varepsilon\rightarrow 0$, and we recover the "non-perturbed" McKean-Vlasov equation $(\ref{McKean-Vlas-syst})$. Thence, under Assumption \ref{ass-mult-omeg}, one considers an embedded Markov chain $Z_n$, for which the state space is the union of small neighborhoods of the compact sets $K_i, 1 \leq i \leq l$, and whose transitions probabilities allow us to estimate the exit and entering times of the empirical vector process $\mu^N$ in the neighborhood of the compact sets $K_i$, thus describing the metastability of the finite $N$-particles system.

\subsection{Metastable phenomena estimates} 
   
Let us introduce some additional notations. For a set $A\subset\left(\mathcal{M}_1(\mathcal{Z})\right)^{2r}$, let $[A]_{\delta}$ denote the open $\delta$-neighborhood of $A$, and $\overline{[A]}_{\delta}$ denote its closure. Moreover, let the stopping time $\tau_A=\inf\{t > 0 | \mu_N (t)\notin A\}$ denote the first exit time from $A$. Recall Assumption \ref{ass-mult-omeg} and let $r_0$ and $r_1$ be two positive numbers such that $0 < r_0< \frac{1}{2} \min_{i,j} \rho^{2r}_0 (K_i, K_j)$ and $0 < r_1 < r_0$, where we recall that $\rho^{2r}_0(\cdot,\cdot)$ is the product metric that equips the product space $(\mathcal{M}_1(\mathcal{Z}))^{2r}$. Denote by $C= (\mathcal{M}_1(\mathcal{Z}))^{2r}\setminus \bigcup_{i=1}^l [K_i]_{r_0}$ the set $(\mathcal{M}_1(\mathcal{Z}))^{2r}$ from which we delete the $r_0$-neighborhoods of $K_i, i= 1,\ldots, l$, and let $\Gamma_i= \overline{[K_i]_{r_0}}$ be the closure of $[K_i]_{r_0}$. Furthermore, denote by  $\gamma_i= [K_i]_{r_1}$ the $r_1$-neighborhood of $K_i$, and $\gamma=\bigcup_{i=1}^l \gamma_i$. Consider the following stopping times:  
\begin{align*}
\tau_0= 0,\quad \sigma_n= \inf\{t\geq \tau_n | \mu^N (t) \in C\},\quad \tau_n= \inf\{t \geq \sigma_{n-1} | \mu^N (t)\in \partial\gamma\},
\end{align*}
 and consider the embedded Markov chain $Z_n=\mu^N(\tau_n)$ defined at hitting times of small neighborhood of the stable limit sets $K_i,i=1,\ldots,l$. The one-step transition probabilities of $Z_n$ are given by
\begin{align*}
P (\nu,\partial\gamma_j)=P_{\nu}\big(Z_n\in \partial\gamma_j \big)=P_{\nu}(\mu^N(\tau_n)\in\partial\gamma_j), 
\end{align*}
when $Z_{n-1}=\nu\in\partial\gamma_i$. The upper and lower bound estimates of these transition probabilities are given in Lemma \ref{Z-chain-prob-lem} and play a central role in the remainder of this section.

We first give an estimate of the stopping time $\tau_1$ of the first reentrance into the $r_1$-neighborhood of one of the compact $K_i$. Note that similar results were established in \cite[Lem. 1.3]{Hwang+Sheu90} for small noise diffusion processes, and in \cite[Lem. 3.6]{Yaso+Sund2019} for complete interaction mean-field systems with jumps. 

\begin{lem} Given $\varepsilon>0$ and $r_0$ small enough, there exists $N_0\geq 1$ such that for any $N\geq N_0$ we have
\label{tau1-estim-lem}
\begin{align}
\mathbb{E}_{\nu}\tau_ 1\leq\exp\{N\varepsilon\},\quad\text{for any $\nu\in\bar{\gamma}$}.
\end{align}
\end{lem}
 \proof  First notice that $\mathbb{E}_{\nu}\tau_1=\mathbb{E}_{\nu}\sigma_0+\mathbb{E}_{\nu}[\tau_1-\sigma_0]$. By Lemma $\ref{exit-mom1-lem}$, there exists $\delta>0$ and $N_0\geq 1$ such that for all $r_0\leq\delta$, $\nu\in\gamma$ and $N\geq N_0$ we have  
\begin{align*}
\mathbb{E}_{\nu}\sigma_0<\exp\{N\varepsilon\}.
\end{align*} 
   
 Moreover $\mathbb{E}_{\nu}[\tau_1-\sigma_0]=\mathbb{E}_{\nu}[\mathbb{E}_{\mu^N(\sigma_0)}\tau_F]$, where $F=(\mathcal{M}_1^N (\mathcal{Z}))^{2r}\setminus\bar{\gamma}$. Since the compact set $F$ does not contain any $\omega$-limit set, Corollary $\ref{coro-appen}$ allows to deduce that there exists a constant $\kappa>0$ such that $\mathbb{E}_{\mu^N(\sigma_0)}\tau_F\leq\kappa$ and thus $\mathbb{E}_{\nu}[\tau_1-\sigma_0]\leq\kappa$. Taking $N$ large enough gives $(\ref{tau1-estim-lem})$. \carre
\\
\\
Let $L=\{1, 2,\ldots, l\}$ be the indices corresponding to the compact sets $K_1,K_2,\ldots, K_l$ given in Assumption \ref{ass-mult-omeg}. For any $W\subset L$, introduce the following stopping times:
\begin{align*}
\hat{\tau}_W&=\inf\left\{t>0: \mu^N(t)\in\cup_{i\in W}\gamma_i\right\},\\
\bar{\tau}_W&=\inf\left\{t> 0: \mu^N(t)\in\cup_{i\in L/W}\gamma_i\right\}.
\end{align*}

For a $W$-graph $g$, set $\tilde{V}(g)=\sum_{(m,n)\in g}\tilde{V}(K_m,K_n)$. For $i\in L/W$ and $j \in W$, let $G_{i,j} (W)$ denote the set of $W$-graphs in which there is a sequence of arrows leading from $i$ to $j$. For any subset $W\subset L$, define
 \begin{align*}
I_i(W ):= \min\{\tilde{V}(g): g \in G(W)\}- \min\{\tilde{V}(g):& g \in G(W\cup \{i\})\quad\mbox{or}\\
                          & g \in G_{i,j} (W\cup \{j\}), i\neq j, j\in L\setminus W \}. 
 \end{align*}

The next lemma gives upper and lower bound estimates for the mean entrance time into small neighborhoods of a set of compacts indexed by $W\subset L$, starting from the neighborhood of a given compact $K_i$, with $i\notin W$. Similar estimates have been obtained in \cite[Part I, Lem. 1.6]{Hwang+Sheu90} and \cite[Lem. 3.10]{Yaso+Sund2019} for small noise
diffusion processes and finite-state mean-field systems on complete graphs respectively.   
 
 \begin{lem} 
 Let $W\subset L$, and let $i \in L \setminus W$. Given $\varepsilon > 0$, there exist $\delta > 0$ and $N_0 \geq 1$ such that for any $r_0 \leq \delta$, $\nu\in\gamma_i$ and $N\geq  N_0$, we have
 \begin{align*}
\exp\{N (I_i (W)-\varepsilon)\}\leq \mathbb{E}_{\nu} [\hat{\tau}_W] \leq \exp\{N (I_i (W ) + \varepsilon)\}. 
 \end{align*}
 \label{mean-entr-esti}
 \end{lem}
 \proof By the strong Markov property, one obtains 
\begin{align*}
\mathbb{E}_{\nu}[\hat{\tau}_W]=\mathbb{E}_{\nu}[\tau_v]=\sum_{m=1}^{\infty}\mathbb{E}[\mathds{1}_{v=m}\tau_m]=\sum_{m=1}^{\infty}m\mathbb{E}[\mathds{1}_{v=m}\mathbb{E}_{\mu^N(\tau_{m-1})}[\tau_1]],
\end{align*} 
 where $v$ is the hitting time of the Markov chain $Z_n$ into the set $\gamma_W=\cup_{i\in W}\gamma_i$. Moreover, using Lemma $\ref{tau1-estim-lem}$, for sufficiently small $r_0$ and sufficiently large $N$, we obtain 
 \begin{align*}
 \mathbb{E}_{\nu}[\hat{\tau}_W]\leq \exp\{N\varepsilon\}\sum_{m=1}^{\infty}m\mathbb{E}[\mathds{1}_{v=m}].
 \end{align*}
 
 Notice that the sum term in the last inequality corresponds to the expectation of the number of steps from $\nu$ until the first entrance in $\cup_{i\in W}\gamma_i$. Thus using the upper bound estimate given in \cite[Ch.6, Lem. 3.4]{Freid+Wentz2012} we find
 \begin{align*}
 \mathbb{E}_{\nu}[\hat{\tau}_W]\leq \exp\{N(I_i(W)+\varepsilon)\}.
 \end{align*}
 
Furthermore, Lemma $\ref{exit-mom2-lem}$ allows us to deduce that, for all sufficiently small $r_0$ and sufficiently large $N$, 
 \begin{align*}
 \mathbb{E}_{\nu}[\tau_1]\geq \exp\{-N\epsilon\}
 \end{align*}
 for all $\nu\in\gamma$. Hence, using the lower bound estimate in \cite[Ch.6, Lem. 3.4]{Freid+Wentz2012}, we deduce that 
 \begin{align*}
 \mathbb{E}_{\nu}[\hat{\tau}_W]\geq \exp\{N(I_i(W)-\varepsilon)\}
 \end{align*}
 for all $\nu \in \gamma_i$ and sufficiency large $N$. The lemma is proved. \carre 
 \\
 \\
Define 
\begin{align*}
I_{i,j} (W ) := \min\{ \tilde{V} (g): g \in G_{i,j} (W)\} - \min\{ \tilde{V}(g): g \in G(W)\}.  
\end{align*}
The following result gives the estimate of the probability that the first entry of $\mu^N$ into a neighborhood of a set $W \subset L$ takes place via a given compact set $K_j$, with $j\in W$, starting from a neighborhood of $K_i$, with $i\in L\setminus W$. 

\begin{lem}
Let $W\subset L$, $i\in L\setminus W$ and $j\in W$. For any $\varepsilon>0$, there exists $\delta>0$ such that for any $0<r_0<\delta$, sufficiently large $N$, and all $\nu\in\gamma_i$ we have 
\begin{align*}
\exp\{-N(I_{ij}(W)+\varepsilon)\}\leq P_{\nu}(\mu^N(\hat{\tau}_W)\in\gamma_j )\leq \exp\{-N(I_{ij}(W)-\varepsilon)\}.
\end{align*}

\label{entr-gamma-j}
\end{lem}
\proof This follows by applying \cite[Ch.6 Lem 3.3]{Freid+Wentz2012} and making use of the estimate in $(\ref{Z-chain-prob})$. \carre
\\
\\
We now introduce the important notion of \textit{cycle} which, roughly speaking, describes how the process runs through its lifetime. Indeed, as explained above, for large but finite $N$, and over large time intervals, there are passages of the process $\mu^N$ between the neighborhoods of the compact sets $K_i$. The cycles then describe the most probable order in which the trajectories of $\mu^N$ traverse these neighborhoods, and the time required to go from one compact to another. Notice nevertheless that the notion of cycle used here was introduced in \cite{Hwang+Sheu90}, which is slightly different from the classical Freidlin and Wentzell notion of cycle \cite{Freid+Wentz2012}. 
 
Recall the definition of $\tilde{V}(K_i,K_j)$ in $(\ref{V-tilde-comp})$, and set $\tilde{V}(K_i)=\min_{j\neq i}\tilde{V}(K_i,K_j)$. From the estimate in $(\ref{Z-chain-prob})$ of the one step transition probabilities of the Markov chain $Z_n$, one can notice that, starting from a neighborhood of a given compact $K_i$, the most  likely set that will be visited by the process $\mu^N$, for large enough $N$, is the one that reaches the minimum $\tilde{V}(K_i)$. Denote by $i\rightarrow j$ if $\tilde{V}(K_i)=\tilde{V}(K_i,K_j)$. This gives us an oriented graph structure where the nodes are given by the set $L=\{1,\ldots,l\}$, and the edges are given by the relation $i\rightarrow j$, for $i,j\in L$. We refer to this graph by $L$. Moreover, for $i,j\in L$, we say that $i\Rightarrow j$ if there exists a sequence of arrows leading from $i$ to $j$, i.e., if there exist $i_1, i_2,\ldots, i_n$ in $L$ such that $i\rightarrow i_1\rightarrow i_2\rightarrow\cdots\rightarrow i_n\rightarrow j$. 

\begin{defi}
A cycle $\pi$ in $L$ is a subgraph of $L$ satisfying:
\begin{enumerate}
\item $i\in\pi$ and $i\rightarrow j$ imply $j\in\pi$,
\item For any $i\neq j$ in $\pi$, $i\Rightarrow j$ and $j\Rightarrow i$.
\end{enumerate}
\end{defi} 
 
For a proof of existence of cycles in $L$, one can consult \cite[Lem. A1]{Hwang+Sheu90}. Next, we describe the decomposition of the set $L$ into a hierarchy of cycles. Setting $L_0=L$, the cycles of rank $1$, or the $1$-cycles, are defined as
 \begin{align*}
 L_1 := \{\pi| \mbox{$\pi$ is a cycle in $L_0$}\} \cup \{i \in L_0|\mbox{$i$ is not in any cycle}\}.
\end{align*}  
 We use the superscript $1$ to refer to the $1$-cycles, namely $\pi^1$. For $\pi^1_1,\pi^1_2 \in L_1$ two 1-cycles, with $\pi^1_1 \neq \pi^1_2$, let $\hat{V}(\pi^1_1 ):= \max\{\tilde{V}(K)| K \in \pi^1_1 \}$ and
\begin{align*}
\tilde{V}(\pi^1_1, \pi^1_2) := \hat{V}(\pi^1_1 )+ \min\{\tilde{V}(K, K')- \tilde{V}(K)| K\in \pi^1_1, K' \in\pi^1_2 \},
\end{align*}
 \begin{align*}
\tilde{V}(\pi^1_1):= \min\{\tilde{V}(\pi^1_1, \pi^1_2): \pi^1_2 \in L_1, \pi^1_2\neq \pi^1_1 \}.
 \end{align*}
 We say that $\pi^1_1\rightarrow\pi^1_2$ if $\tilde{V}(\pi^1_1)= \tilde{V}(\pi_1^1, \pi^1_2)$, and $\pi^1_1\Rightarrow\pi^1_2$ if there is a sequence of arrows leading from $\pi^1_1$ to $\pi^1_2$. This gives a cycle of second order, or a $2$-cycle, and we use the notation $\pi^2$ to denote them. In this way, $\tilde{V}(\pi_1)$ represents the exit rate from a cycle $\pi^1_1$ as specified by Lemma \ref{exit-time-cycle} below.  
 
 By recurrence, assuming that we have defined up to $m$-cycles, we set
 \begin{align*}
 L_{m}= \{\pi^{m}|\mbox{$\pi^{m}$ is a cycle in $L_{m-1}$}\}\cup \{\pi^{m-1} \in L_{m-1}| \mbox{$\pi^{m-1}$ is not in any $m$-cycle}\}.
 \end{align*}
  For $\pi_1^{m},\pi_2^{m} \in L_m$ , $\pi_1^{m} \neq \pi_2^{m}$, let $\hat{V}(\pi_1^{m} ):= \max\{\tilde{V}(\pi^{m-1})| \pi^{m-1} \in \pi_1^{m} \}$ and
\begin{align*}
\tilde{V}(\pi_1^{m}, \pi_2^{m}) := \hat{V}(\pi_1^{m})+ \min\{\tilde{V}(\pi_1^{m-1}, \pi_2^{m-1})- \tilde{V}(\pi_1^{m-1})| \pi_1^{m-1}\in \pi_1^{m}, \pi_2^{m-1}\in \pi_2^{m} \},
\end{align*}
 \begin{align*}
\tilde{V}(\pi_1^{m}):= \min\{\tilde{V}(\pi_1^{m}, \pi_2^{m}): \pi_2^{m} \in L_{m}, \pi_2^{m}\neq \pi_1^{m} \}.
 \end{align*}
 
 We say that $\pi_1^{m} \rightarrow \pi_2^{m}$ if $\tilde{V} (\pi_1^{m})= \tilde{V} (\pi_1^{m}, \pi_2^{m})$, which gives us the $(m+1)$-cycle. We keep going until we reach a given $m$ for which the set of $m$-cycle is a singleton. From now on, we will use, with slight abuse of notation, the notation $\pi^k$ to refer both to the $k$-cycle $\pi^k$ and the set of elements  of $L$ constituting it. Recall that for $W\subset L$, $\gamma_W=\cup_{i\in W} \gamma_i$. Before stating some important results about cycles, we first introduce an example to help the reader have a clear picture about this notion.

\begin{ex} Let the set of indices be $L=\{1,2,\ldots,8\}$, and consider the matrix corresponding to the values of $\tilde{V}(K_i,K_j)$, for $i, j\in L$

 \begin{center}
 $
\left(\begin{array}{c c c c c c c c}
	0 & 2 & 4 & 6 & 7 & 6& 12 & 8 \\
	6 & 0 & 1 & 8 & 9& 11 & 13 & 15\\
	5 & 7 & 0 & 10 & 11 & 8 & 9 & 11 \\
	5 & 10 & 20 & 0 &3 & 4& 8& 9\\ 
	10 & 11 & 12 & 7 & 0& 18 & 16& 21 \\
	7 & 11 & 13 & 9 & 11 & 0 & 8 & 6 \\
	8 & 9& 14& 8& 13& 4& 0& 10\\
	15& 12& 9& 7& 10& 11& 5& 0  
	\end{array}\right)
$
 \end{center}
 
 Set $L_0=L$. Using the definition above we find three 1-ycles: $\pi^1_1=\{1,2,3\}$ characterized by the edges $1\rightarrow 2$, $2\rightarrow 3$, and  $3\rightarrow 1$, $\pi^1_2=\{4,5\}$ characterized by the edges $4\rightarrow 5$,  and $5\rightarrow 4$, and  finally $\pi^1_3=\{6,7,8\}$ characterized by the edges $6\rightarrow 8$, $8\rightarrow 7$ and $7\rightarrow 6$. Thus $L_1=\{\pi_1^1,\pi_2^1,\pi_3^1\}$. 
 
  Now in order to find the 2-cycles, we use again the construction above. Straightforward computations give: $\hat{V}(\pi_1^1)=5$, $\tilde{V}(\pi_1^1,\pi_2^1)=9$, $\tilde{V}(\pi_1^1,\pi_3^1)=8$. Thus we deduce that $\pi_1^1\rightarrow\pi_1^3$. Moreover we obtain $\hat{V}(\pi_3^1)=6$, $\tilde{V}(\pi_3^1,\pi_1^1)=7$, $\tilde{V}(\pi_3^1,\pi_2^1)=8$, and thus $\pi^1_3\rightarrow\pi_1^1$. Finally we find $\hat{V}(\pi_2^1)=7$, $\tilde{V}(\pi_2^1,\pi_1^1)=9$, $\tilde{V}(\pi_2^1,\pi_3^1)=8$ from which we deduce that $\pi_2^1\rightarrow\pi_3^1$. Hence we have $L_2=\{\{\pi_1^1,\pi_3^1\},\{\pi_2^1\}\}$ represents the collection of $2$-cycles. Denote by $\pi^2_1=\{\pi_1^1,\pi_3^1\}$, $\pi_2^2=\{\pi_2^1\}$ the two elements of $L_2$. Now, in order to identify the set of $3$-cycles, we find again by simple calculations of the following quantities: $\hat{V}(\pi^2_1)=8$, and $\tilde{V}(\pi_1^2,\pi_2^2)=10$. Moreover, $\hat{V}(\pi^2_2)=7$ and  $\tilde{V}(\pi_2^2,\pi_1^2)=8$. Therefore, $\pi_1^2\rightarrow\pi_2^2$ and $\pi_2^2\rightarrow\pi_1^2$, which give us the unique $3$-cycle $\pi^3=\{\pi^2_1,\pi_2^2\}$, the only element of $L_3$. Now we stop since $L_3$ is a singleton. See Figure \ref{cycle-fig} for an illustration.   
\label{cycle-ex}
\end{ex}
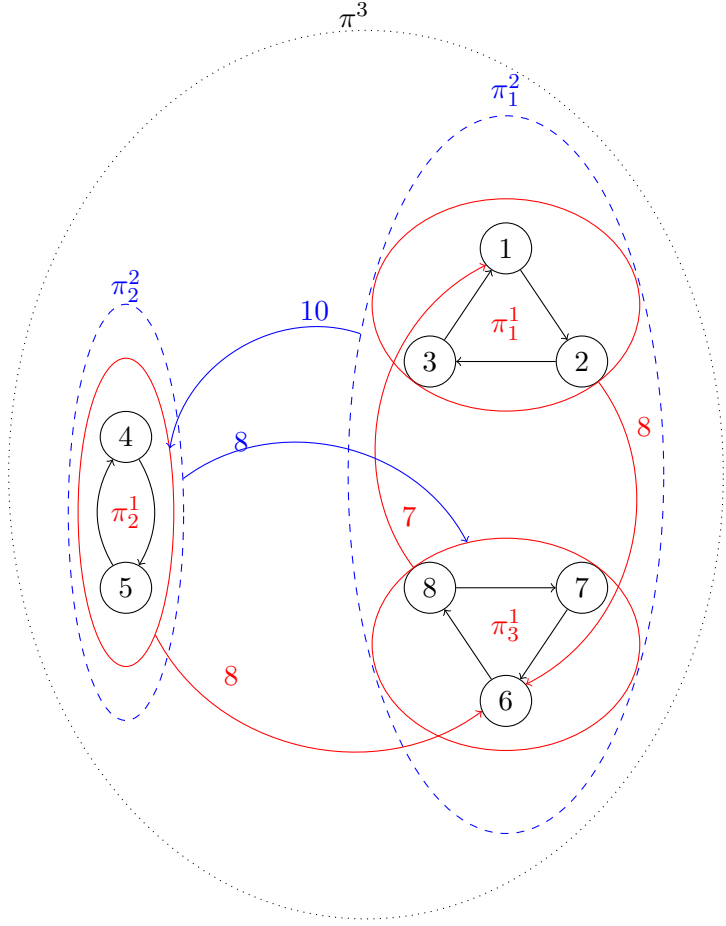
\begin{figure}[h]
\centering
\begin{tikzpicture}
\draw (0,3.5)node[circle,draw](1){$1$}; 
\draw (1,2)node[circle,draw](2){$2$}; 
\draw (-1,2)node[circle,draw](3){$3$}; 
\draw (0,2.5)node (){\textcolor{red}{$\pi_1^1$}}; 
\draw (0,5.6)node (){\textcolor{blue}{$\pi_1^2$}}; 
\draw (0,-1.5)node (){\textcolor{red}{$\pi_3^1$}}; 
\draw[->] (1) -- (2);
\draw[->] (2) -- (3);
\draw[->] (3) -- (1);
\draw (0,-2.5)node[circle,draw](6){$6$}; 
\draw (1,-1)node[circle,draw](7){$7$}; 
\draw (-1,-1)node[circle,draw](8){$8$}; 
\draw[->] (6) -- (8);
\draw[->] (8) -- (7);
\draw[->] (7) -- (6);
\draw (-5,1)node[circle,draw](4){$4$}; 
\draw (-5,-1)node[circle,draw](5){$5$}; 
\draw (-5,0)node (){\textcolor{red}{$\pi_2^1$}}; 
\draw (-5,3)node (){\textcolor{blue}{$\pi_2^2$}}; 
\draw (-2,6.6)node (){\textcolor{black}{$\pi^3$}}; 
\path[->]
        (4) edge[bend left=30] node[pos=0.5,above] {} (5);
\path[->]
        (5) edge[bend left=30] node[pos=0.5,above] {} (4);
\node[ellipse, draw=red, fit=(1) (2) (3),inner sep=-1mm] (pi11) {};
\node[ellipse, draw=red, fit=(6) (7) (8),inner sep=-1mm] (pi13) {};
\node[ellipse, dashed, draw=blue, fit=(pi11) (pi13), inner sep=-3mm] (pi21)  {};
\node[ellipse, draw=red, fit=(4) (5),inner sep=1mm] (pi12) {};
\path[->]
        (pi13) edge[bend left=50, draw=red] node[pos=0.1,above] {$\quad\textcolor{red}{7}$} (1);
\path[->]
        (pi11) edge[bend left=50, draw=red] node[pos=0.2,above] {$\quad\textcolor{red}{8}$} (6);        
\path[->]
        (pi12) edge[bend right=50,draw=red] node[pos=0.2,above] {$\qquad\textcolor{red}{8}$} (6);  
\node[ellipse, dashed, draw=blue, fit=(pi12), inner sep=-1mm] (pi22)  {};
\path[->]
        (pi21) edge[bend right=50, draw=blue] node[pos=0.3,above] {$\qquad\textcolor{blue}{10}$} (pi12); 
\path[->]
        (pi22) edge[bend left=50, draw=blue] node[pos=0.1,above] {$\qquad\textcolor{blue}{8}$} (pi13);   
\node[ellipse, dotted, draw=black, fit=(pi22) (pi21), inner sep=-6mm] (pi3)  {};        
\end{tikzpicture}
\caption{The hierarchy of cycles in Example \ref{cycle-ex}}
\label{cycle-fig}
\end{figure}

The next result gives an estimate of the mean exit time from a cycle.
\begin{lem}
 Let $\pi^k$ be a $k$-cycle, and let  $i \in \pi^k$. Moreover, denote by $W=L\setminus\pi^k$ the set of compacts not contained in $\pi^k$. Then, given $\varepsilon > 0$, there exist $\delta > 0$ and $N_0 \geq 1$ such that for any $r_1 \leq \delta$, $\nu\in\gamma_i$ and $N\geq  N_0$, we have
 \begin{align*}
\exp\{N (\tilde{V}(\pi^k)-\varepsilon)\}\leq \mathbb{E}_{\nu} [\hat{\tau}_W] \leq \exp\{N (\tilde{V}(\pi^k) + \varepsilon)\}. 
 \end{align*}
 \label{exit-time-cycle}
\end{lem} 
 
 \proof We have, from \cite[Lem. A3]{Hwang+Sheu90} and \cite[Cor. A4]{Hwang+Sheu90}, that $I_i(W)=\tilde{V}(\pi^k)$. Using this together with Lemma $\ref{mean-entr-esti}$ leads to the result. \carre
 \\
 \\
The next result gives an estimate of the probability of going from one $k$-cycle to another without passing by any other element of $L^k$.  
\begin{lem}
Let $\pi_1^k$ , $\pi_2^k$ be two distinct $k$-cycles, and let $i\in\pi_1^k$. Moreover, denote $W=L\setminus\pi_1^k$. Then, given $\varepsilon > 0$, there exist $\rho > 0$ and $N_0 \geq 1$ such that, for all $r_1 \leq\rho$, $\nu \in \gamma_i$, and $N\geq N_0$, we have
\begin{align*}
\exp\{-N(\tilde{V} (\pi_1^k, \pi_2^k)- \tilde{V} (\pi_1^k)+ \varepsilon)\} \leq P_{\nu}(\mu^N(\hat{\tau}_W)\in\gamma_{\pi_2^k})
\leq \exp\{-N ( \tilde{V}(\pi_1^k, \pi_2^k)-\tilde{V}(\pi_1^k)-\varepsilon)\}.
\end{align*}
\end{lem} 
 
 \proof From Lemma $\ref{entr-gamma-j}$ we have, for each $j\in\pi_2^k$ and large enough $N$,
 \begin{align*}
\exp\{-N (I_{i,j} (W)+\varepsilon)\}\leq P_{\nu}(\mu^N(\hat{\tau}_W)\in\gamma_j )\leq \exp\{-N (I_{i,j} (W)-\varepsilon)\}.
\end{align*}
 
 Therefore, summing over the disjoint compacts $K_i$ we obtain,  
  \begin{align*}
\sum_{j\in\pi_2^k} \exp\{-N (I_{i,j} (W)+\varepsilon)\}\leq P_{\nu}(\mu^N(\hat{\tau}_W)\in\gamma_{\pi_2^k} )\leq\sum_{j\in\pi_2^k} \exp\{-N (I_{i,j} (W)-\varepsilon)\}.
\end{align*}
From the sums above we select the term decreasing more slowly than the remaining ones, which gives us
  \begin{align*}
\exp\{-N (\min\{I_{i,j} (W):j\in\pi_2^k\}+\varepsilon)\}\leq P_{\nu}(\mu^N(\hat{\tau}_W)\in\gamma_{\pi_2^k} )\leq \exp\{-N (\min\{I_{i,j} (W):j\in\pi_2^k\}-\varepsilon)\}.
\end{align*}

Finally, from \cite[Lem. A5]{Hwang+Sheu90} we have that $\min\{I_{ij}(W):j\in\pi_2^k\}=\tilde{V}(\pi_1^k, \pi_2^k)-\tilde{V}(\pi_1^k)$, which leads to the stated result.  \carre

\subsection{Convergence to the invariant measure}
As aforementioned, there exists, under Assumption \ref{ass-prin}, a unique invariant measure $\wp^N$ for the empirical vector process $\mu^N$. Thus, the distribution of $\mu^N$ converges, as time $t\rightarrow\infty$, towards $\wp^N$. However, one might investigate the corresponding rate of convergence. Therefore, using the results from the previous section we show that, when the time is of order $\exp\{(\Lambda+\delta) N\}$, with $\delta>0$ and $\Lambda$ a suitable constant detailed above, the empirical measure is very close to its invariant measure $\wp^N$. Interestingly, despite the heterogeneity introduced by the block structure, a similar constant appears in the case of small noise diffusion processes \cite[Lem. 1.3]{Hwang+Sheu90}, and in the case of homogeneous mean-field systems with jumps in \cite[Lem. 3.6]{Yaso+Sund2019}. Before stating our main result, let us introduce further notations and intermediate results. 
 
 Let $i_0\in L$ be such that $\min\{ \tilde{V} (g) : g \in G(i_0 )\} = \min\{\tilde{V}(g) : g \in G(i), i \in L\}$. Moreover, define
 \begin{align*}
\Lambda := \min\{\tilde{V} (g): g\in G(i), i \in L\}- \min\{ \tilde{V}(g): g \in G(i, j), i, j \in L, i \neq j\}. 
 \end{align*}
Let $P_T (\nu, \cdot) = P_{\nu} (\mu^N (T )\in\cdot)$ denote the transition probability kernel associated with the empirical process $\mu^N$. The next result gives a lower bound for the transition probability $P_T (\nu, K_{i_0})$ of reaching a small neighborhood of $K_{i_0}$ when $T$ is of order $\exp\{N (\Lambda- \delta_0 )\}$, for some $\delta_0 > 0$.
 
\begin{theo} Given $\varepsilon>0$, there exist $\delta_0 > 0$, $r> 0$ and $N_0\geq 1$ such that, for all $r_1 \leq r$,
$N \geq N_0$, and $\nu\in (\mathcal{M}_1(\mathcal{Z}))^{2r}$, we have
\begin{align*}
P_{T_0}(\nu, \gamma_{i_0}) \geq \exp\{-N \varepsilon\},
\end{align*}
where $T_0=\exp\{N (\Lambda-\delta_0)\}$. Furthermore, there exist $\nu_0\in (\mathcal{M}_1(\mathcal{Z}))^{2r}$ and $\beta> 0$ such that, for all $N \geq N_0$ and $\nu\in [\nu_0]_{r_1}$,
\begin{align*}
P_{T_0} (\nu,\gamma_{i_0}) \leq \exp\{-N \beta\}.
\end{align*}
 \label{entr-i-zero}
 \end{theo}

 \proof Replacing the space $\mathcal{M}_1(\mathcal{Z})$ by the product space $(\mathcal{M}_1(\mathcal{Z}))^{2r}$ together with the corresponding metrics, the proof given in \cite[Th. 3.21]{Yaso+Sund2019} adapted from \cite[Th. 2.3, part I]{Hwang+Sheu90} holds verbatim. \carre
 
 \begin{coro}
 \label{cor-entr-i-zero}
 Under the conditions of Theorem $\ref{entr-i-zero}$, for all $\nu\in (\mathcal{M}_1(\mathcal{Z}))^{2r}$, $\xi \in \gamma_{i_0}$, and $N$ sufficiently large, we have
 \begin{align*}
P_{T_0} (\nu, \xi) \geq \exp\{-2N \varepsilon\}. 
 \end{align*}
 \end{coro}
 
 \proof The proof follows verbatim the proof given in \cite[Coro. 3.22]{Yaso+Sund2019}. \carre 
 \\
 \\
 We state now the main result of this section, which gives the time scale at which the empirical vector converges towards its invariant measure.    
 
\begin{theo} 
 There exists a constant $\Lambda\geq 0$ such that, for any $\delta>0$, there exist $\varepsilon > 0$ and $N_0 \geq 1$
such that, for all $\nu\in (\mathcal{M}_1 (\mathcal{Z}))^{2r}$ and $N \geq N_0$,
\begin{align*}
\left|\mathbb{E}_{\nu} \left[f (\mu^N(T ))\right]-\int f d\wp^N \right| \leq \|f\|_{\infty}  \exp\{- \exp(N \varepsilon)\},
\end{align*}
where $T = \exp\{N (\Lambda + \delta)\}$ and $f$ is any bounded and measurable function on $(\mathcal{M}_1(\mathcal{Z}))^{2r}$. 
\label{rate-conv-invar-meas} 
 \end{theo}
 \proof Let $\varepsilon>0$, and let $T_0$, $\delta_0$, $r$, $r_1$ and $N_0 \geq 1$ be as in the statement of Theorem $\ref{entr-i-zero}$. By Corollary $\ref{cor-entr-i-zero}$ we have that, for any $\nu\in (\mathcal{M}_1 (\mathcal{Z}))^{2r}$ and $\xi\in [ K_{i_0}]_{r_1}$, 
 \begin{align}
P_{T_0} (\nu, \xi)\geq \exp\{-2N \varepsilon\}.
\label{ineq-1-conv-emp}
\end{align}
 Moreover, by the Markov property and using Corollaries $\ref{cor-entr-i-zero}$ and $\ref{coro-ldp-emp-proc}$ we have that, for any $\nu\in (\mathcal{M}_1 (\mathcal{Z}))^{2r}$, $\xi\notin [ K_{i_0}]_{r_1}$, and some fixed $t$, 
\begin{equation}
\begin{split}
P_{T_0} (\nu, \xi)&\geq\sum_{\nu'\in[K_{i_0}]} P_{T_0-t} (\nu,\nu')P_t(\nu',\xi)\\
&\geq \exp\{-2N\varepsilon\}\inf_{\nu'\in[K_{i_0}]}P_t (\nu' , \xi) \\
&\geq \exp\{-2N \varepsilon\} \exp\{-N \sup_{\nu'\in[K_{i_0}]} S_t (\xi|\nu')\},
\label{ineq-1-conv-emp}
\end{split}
\end{equation}
for sufficiently large $N$. Notice that, thanks to $(\ref{ineq-1-conv-emp})$ and since $S_T(\xi,\nu)$ is bounded for all $\nu,\xi\in (\mathcal{M}_1(\mathcal{Z}))^{2r}$, one can find a positive function $U(\xi)$ such that $U(\xi)=0$ for $\xi\in [ K_{i_0}]_{r_1}$, and $U(\xi)\geq \sup_{\nu'\in[K_{i_0}]} S_t (\xi|\nu')$ for $\xi\notin [ K_{i_0}]_{r_1}$. 

Define  $\pi_N (\xi) = c_N \exp\{-N U (\xi)\}$ such that $U(\xi)=0$ for $\xi\in [ K_{i_0}]_{r_1}$, and 
 \begin{align*}
P_{T_0} (\nu, \xi)\geq c_N \exp\{-N U(\xi)\}
\end{align*} 
 for all $\nu\in (\mathcal{M}_1(\mathcal{Z}))^{2r}$, $\xi\notin [ K_{i_0}]$, and sufficiently large $N$, where $c_N$ is a constant chosen such that $\pi_N$ is a probability measure on the product space $(\mathcal{M}_1(\mathcal{Z}))^{2r}$. Moreover, set $Q_{T_0}(\nu,\cdot)=P_{T_0} (\nu, \cdot)/ \pi_N (\cdot)$. Hence, by $(\ref{ineq-1-conv-emp})$, $Q_{T_0}(\nu,\xi)\geq 1$ for any $\nu,\xi\in (\mathcal{M}_1(\mathcal{Z}))^{2r}$. Let $f$ be a bounded measurable function on $(\mathcal{M}_1(\mathcal{Z}))^{2r}$. Therefore, for any $\nu_1,\nu_2\in (\mathcal{M}_1(\mathcal{Z}))^{2r}$ and sufficiently large $N$ we have 
 \begin{equation}
 \begin{split}
 \mathbb{E}_{\nu_1} (f (\mu^N (T_0)))-  \mathbb{E}_{\nu_2} (f (\mu^N (T_0)))&=\int_{(\mathcal{M}_1(\mathcal{Z}))^{2r}}f(\xi)P_{T_0}(\nu_1,d\xi)-\int_{(\mathcal{M}_1(\mathcal{Z}))^{2r}}f(\xi)P_{T_0}(\nu_2,d\xi)\\
 &=\int_{(\mathcal{M}_1(\mathcal{Z}))^{2r}}f(\xi)Q_{T_0}(\nu_1,d\xi)\pi_N(d\xi)-\int_{(\mathcal{M}_1(\mathcal{Z}))^{2r}}f(\xi)Q_{T_0}(\nu_2,d\xi)\pi_N(d\xi)\\
  &=\int_{(\mathcal{M}_1(\mathcal{Z}))^{2r}}f(\xi)(Q_{T_0}(\nu_1,d\xi)-\exp\{-2N\varepsilon\})\pi_N(d\xi)\\
  &\qquad\qquad-\int_{(\mathcal{M}_1(\mathcal{Z}))^{2r}}f(\xi)(Q_{T_0}(\nu_2,d\xi)-\exp\{-2N\varepsilon\})\pi_N(d\xi)\\
  &\leq \sup_{\xi\in(\mathcal{M}_1(\mathcal{Z}))^{2r}}f(\xi)\bigg(1-\exp\{-2N\varepsilon\}\bigg)\\
  &\qquad\qquad-\inf_{\xi\in(\mathcal{M}_1(\mathcal{Z}))^{2r}}f(\xi)\bigg(1-\exp\{-2N\varepsilon\}\bigg)\\
  &= \big(1-\exp\{-2N\varepsilon\}\big)|\sup_{\xi}f(\xi)-\inf_{\xi}f(\xi)|.
 \end{split}
\end{equation}  
 Thence we find, 
  \begin{equation}
 \begin{split}
 \sup_{\nu_1,\nu_2}\left|\mathbb{E}_{\nu_1} (f (\mu^N (T_0)))-  \mathbb{E}_{\nu_2} (f (\mu^N (T_0)))\right|\leq \big(1-\exp\{-2N\varepsilon\}\big)\|f\|_{\infty}.
  \end{split}
\end{equation}  
 Hence, by repeating the previous steps $k$ times using the Chapman-Kolmogorov property given that $\mu^N$ is Markov, we find 
  \begin{equation}
 \begin{split}
 \sup_{\nu_1,\nu_2}\left|\mathbb{E}_{\nu_1} (f (\mu^N (kT_0)))-  \mathbb{E}_{\nu_2} (f (\mu^N (kT_0)))\right|\leq \bigg(1-\exp\{-2N\varepsilon\}\bigg)^k\|f\|_{\infty},
  \end{split}
\end{equation}  
 from which we deduce that 
  \begin{equation}
 \begin{split}
 \sup_{\nu}\left|\mathbb{E}_{\nu} (f (\mu^N (kT_0)))-  \langle \wp^N,f\rangle\right|\leq \bigg(1-\exp\{-2N\varepsilon\}\bigg)^k\|f\|_{\infty}.
  \end{split}
\end{equation}
 Choose $k = \exp\{N (\delta_0 + \delta)\}$. Thus $kT_0=\exp\{N(\delta+\Lambda)\}$. Moreover, using the property $(1+\frac{x}{n})^n\rightarrow\exp\{x\}$ as $n\rightarrow\infty$, one obtains for large $k$ 
 \begin{equation}
 \begin{split}
 \sup_{\nu}\left|\mathbb{E}_{\nu} (f (\mu^N (kT_0)))-  \langle \wp^N,f\rangle\right|\leq \exp\bigg\{-\exp\bigg\{N(-2\varepsilon+\delta_0+\delta)\bigg\}\bigg\}\|f\|_{\infty}.
  \end{split}
\end{equation} 
 Choosing $\varepsilon$ small enough such that $\varepsilon'=-2\varepsilon+\delta_0+\delta>0$, one finally obtains  
  \begin{equation}
 \begin{split}
 \sup_{\nu}\left|\mathbb{E}_{\nu} (f (\mu^N (T)))-  \langle \wp^N,f\rangle\right|\leq \exp\{-\exp\{N\varepsilon'\}\}\|f\|_{\infty},
  \end{split}
\end{equation}
 for large enough $N$, with  $T=\exp\{N(\delta+\Lambda)\}$. The theorem is proved. \carre


\begin{appendices}

\section{Freidlin-Wentzell program}
\label{Freid-Went-appe}
We give here generalizations to our setting of series of Lemmas introduced in \cite{Freid+Wentz2012} in the case of diffusion processes and generalizations in \cite{Bork+Sund2012} to the case of jump processes in one homogeneous population. These results play an important role in the study of the large-time behavior of the system. Since this is a slight generalization, we will give full proofs only when needed and refer to the previous references otherwise.

\begin{lem}[\cite{Freid+Wentz2012}, Chap.6 Lem. 1.2] For any $ \varepsilon>0$ and any compact set $K\subset(\mathcal{M}_1(\mathcal{Z}))^{2r}$, there exists a $T_0$ such that for any $\nu,\xi\in K$ there exists a function $\mu(t)$, with $t\in[0,T]$, satisfying $\mu(0)=\nu$, $\mu(T)= \xi$, $T\leq T_0$ with $S_{[0,T]} (\mu|\nu)\leq V (\xi|\nu)+\varepsilon$. 
\end{lem} 
\proof The proof of \cite[Lem. A.1.]{Bork+Sund2012} holds verbatim by replacing the metric $\rho_0(\cdot,\cdot)$ by the product metric $\rho^{2r}_0(\cdot,\cdot)$. \carre  

\begin{lem}[\cite{Freid+Wentz2012}, Chap.6 Lem. 1.6]
Let all points of a compact set $K\subset\left(\mathcal{M}_1(\mathcal{Z})\right)^{2r}$ be equivalent to each other, but not to any other point in $\left(\mathcal{M}_1(\mathcal{Z})\right)^{2r}$. Then, for any $\varepsilon > 0$, $\delta > 0$, $\nu$, $\xi\in K$, there exist a $T > 0$ and a function $\mu(t)$ defined on $[0,T]$, with $\mu(0)=\nu$, $\mu(T)=\xi$, $\mu(t) \in [K]_{\delta}$ for all $t\in [0,T]$, and $S_{[0,T ]}(\mu|\nu)<\varepsilon$. 
\end{lem}
\proof Using Lemma \ref{T-rate-prop}, the proof follows verbatim the proof of \cite[Lem. A.2.]{Bork+Sund2012}. \carre
\\
\\
Let $A\in(\mathcal{M}_1(\mathcal{Z}))^{2r}$, and define the stopping time
\begin{align*}
\tau_A=\inf\{t > 0 | \mu_N (t)\notin A\},
\end{align*}
which gives the first exit time from $A$. The law of this exit time depends on $N$ and $\nu^N$ through the law $p_{\nu^N}^N$ of $ \mu^N$.

\begin{lem}[\cite{Freid+Wentz2012}, Chap.6 Lem. 1.7]
Let all points of a compact set $K\subset(\mathcal{M}_1(\mathcal{Z}))^{2r}$ be equivalent to each other, and let $K\neq (\mathcal{M}_1(\mathcal{Z}))^{2r}$. For any $\varepsilon>0$, there exists a $\delta>0$ such that, for all sufficiently large $N$ and all $\nu\in [K]_{\delta}$, we have
\begin{align}
\mathbb{E}\tau_{[K]_{\delta}}  < e^{N\varepsilon},
\label{exit-mom1}
\end{align}
where the expectation is with respect to the measure $p^N_{\nu}$. 
\label{exit-mom1-lem}
\end{lem}

\proof Using Corollary \ref{coro-ldp-emp-proc}, the proof of \cite[Lem. A.3.]{Bork+Sund2012} holds verbatim by replacing the metric $\rho_0(\cdot,\cdot)$ by the product metric $\rho^{2r}_0(\cdot,\cdot)$. \carre

\begin{lem}[\cite{Freid+Wentz2012}, Chap.6 Lem. 1.8]
Let $K$ be an arbitrary compact subset of $(\mathcal{M}_1(\mathcal{Z}))^{2r}$, and let $G$ be a neighborhood of $K$. For any $\varepsilon > 0$, there exists a $\delta> 0$ such that, for all sufficiently large $N$ and all $\nu$ belonging to  $\overline{[K]}_{\delta}$, we have
\begin{align}
\mathbb{E}_{\nu}\left[\int_{[0,\tau_G]}\mathds{1}_{\overline{[K]}_{\delta}}(\mu^N(t))\right]\geq e^{-\varepsilon N},
\label{exit-mom2}
\end{align}
where the expectation is with respect to the measure $p_{\nu}^N$.
\label{exit-mom2-lem}
\end{lem}

\proof The proof of \cite[Lem. A.4.]{Bork+Sund2012} holds verbatim by replacing the metric $\rho_0(\cdot,\cdot)$ by the product metric $\rho^{2r}_0(\cdot,\cdot)$. \carre

\begin{lem}[\cite{Freid+Wentz2012}, Chap.6 Lem. 1.9]
Let $K$ be a compact subset of $(\mathcal{M}_1(\mathcal{Z}))^{2r}$ not containing any $\omega$-limit set entirely. Then, there exist positive constants $c$ and $T_0$ such that, for all sufficiently large $N$, any $T > T_0$ , and any $\nu\in K$, we have
\begin{align*}
p_{\nu}^N \{\tau_K>T\} \leq e^{-Nc(T-T_0)}.
\end{align*}
\end{lem}

\proof The proof of \cite[Lem. A.5.]{Bork+Sund2012} holds verbatim by replacing the metric $\rho_0(\cdot,\cdot)$ by the product metric $\rho^{2r}_0(\cdot,\cdot)$. \carre 

\begin{coro} Let $K$ be a compact set not containing any $\omega$-limit set entirely. There exists a positive integer $N_0$ and a positive constant $c$ such that for $N\geq N_0$ and any $\nu\in K$, we have
\begin{align*}
E[\tau_K]\leq T_0 + 1/(cN_0 )
\end{align*}

where the expectation is with respect to the measure $p^N_{\nu}$. 
\label{coro-appen}
\end{coro}

\proof See \cite{Freid+Wentz2012}, Chap.6, p. 149. \carre
\\
\\
 Recall that $r_0$ and $r_1$ are positive numbers such that $0 < r_0< \frac{1}{2} \min_{i,j} \rho^r_0 (K_i, K_j)$ and $0 < r_1 < r_0$. We denote by $C= (\mathcal{M}_1(\mathcal{Z}))^{2r}\setminus \bigcup_{i=1}^l [K_i]_{r_0}$, and let $\Gamma_i= \overline{[K_i]_{r_0}}$ be the closure of $[K_i]_{r_0}$, for $i=1,\ldots,l$. Moreover, we denote by  $\gamma_i= [K_i]_{r_1}$ the $r_1$-neighborhood of $K_i$ and $\gamma=\bigcup_{i=1}^l \gamma_i$. Moreover, recall the following stopping times $\tau_0= 0$, $\sigma_n= \inf\{t\geq \tau_n | \mu^N (t) \in C\}$ and $\tau_n= \inf\{t \geq \sigma_{n-1} | \mu^N (t)\in \gamma\}$, and consider the embedded Markov chain of states at hitting times of neighborhood of the stable limit sets $Z_n=\mu^N(\tau_n)$. Let
\begin{align*}
p^{N} (\nu,\overline{\gamma_j})=p^N_{\nu}\big(Z_n\in \overline{\gamma_i} \big), 
\end{align*}
when $Z_{n-1}=\nu$ is the one-step transition probability of $Z_n$.

\begin{lem}[\cite{Freid+Wentz2012}, Chap.6 Lem. 2.1]
For any $\varepsilon>0$, there is a small enough $r_0>0$ such that for any $r_2$ satisfying $0<r_2<r_0$, there is an $r_1$ satisfying $0<r_1<r_2$ such that for all sufficiently large $N$, for all $\nu\in [K_i]_{r_2}$, the one-step transition probabilities of $Z_n$ satisfy
\begin{align}
\exp\big\{-N( \tilde{V}(K_i, K_j)+\varepsilon)\big\}\leq p^{N}(\nu, \bar{\gamma_j})\leq \exp\big\{-N(\tilde{V} (K_i,K_j)- \varepsilon)\big\}.
\label{Z-chain-prob}
\end{align}
\label{Z-chain-prob-lem}
\end{lem}

\proof Replacing the space $\mathcal{M}_1(\mathcal{Z})$ by the product space $(\mathcal{M}_1(\mathcal{Z}))^{2r}$ endowed  with the product norm $\rho_0^{2r}$, the proof of \cite[Lem. A.6]{Bork+Sund2012} holds verbatim. \carre  
\\
\\
Finally, the next result gives the limiting behavior of the empirical measure vector $\mu^N$ when we first let $t$ go to infinity, then let $N$ go to infinity. 

\begin{theo}[\cite{Freid+Wentz2012}, Ch. 6, Th. 4.1.]
Assume Assumption $\ref{ass-mult-omeg}$ holds true. Then for any $\varepsilon> 0$, there exists $r_1 > 0$, which can be chosen arbitrarily small, such that the $\wp^N$-measure of the $r_1$-neighborhood $\gamma_i$ of the compact $K_i$ satisfies
\begin{align*}
\exp\left\{-N\left(W (K_i)- \min_i W (K_i)+\varepsilon\right)\right\}\leq \wp^N\left(\gamma_i\right)\leq \exp\left\{-N\left(W (K_i)- \min_i W (K_i)-\varepsilon\right)\right\}
\end{align*}
for all sufficiently large $N$. 
\label{inv-meas-ineq}
\end{theo}

\proof Choose small positive $0<r_1<r_2<r_0$ such that the inequalities in $(\ref{Z-chain-prob})$, $(\ref{exit-mom1})$ and  $(\ref{exit-mom2})$ are satisfied for sufficiently large $N$ with $\varepsilon/4l$ replacing $\varepsilon$. One can notice that the Markov chain $Z_n$ is irreducible and thus has an invariant measure. By \cite[Ch.6, Lem 3.2]{Freid+Wentz2012}, and by selecting the exponential terms that decrease more slowly, when $N$ is large enough, the values of the normalized invariant measure $\vartheta^N$ of the chain $Z_n$ lie in the interval 
\begin{align*}
\exp\left\{-N\left( W (K_i)- \min_i W (K_i) \pm\frac{l-1}{2l}\varepsilon\right)\right\}.
\end{align*}
Recall the following formula which expresses, up to a factor, the invariant measure $\wp^N$ of the process $\mu^N$ in terms of the invariant measure $\vartheta^N$ of the chain $Z_n$ on $\partial \gamma$ (see \cite[Ch.6, eqn. (4.1)]{Freid+Wentz2012}):
\begin{align}
\wp^N(B) =\int_{\partial \gamma} \vartheta^N(d\nu)\mathbb{E}_{\nu}\left[\int_0^{\tau_1}\mathds{1}_{B}(\mu^N(t))dt\right],
\label{khas-form}
\end{align}
where the expectation is  with respect to the measure $p_{\nu}^N$ with starting conditions $\nu$. Using this formula we find, for any $\gamma_i$, with $i=1,\ldots,l$,
\begin{align*}
\wp^N(\gamma_i) =\int_{\partial \gamma} \vartheta(d\nu)\mathbb{E}\left[\int_0^{\tau_1}\mathds{1}_{\gamma_i}(\mu^N(t))dt\right]=\int_{\partial \gamma_i} \vartheta(d\nu)\mathbb{E}\left[\int_0^{\tau_1}\mathds{1}_{\gamma_i}(\mu^N(t))dt\right].
\end{align*}
Using the estimates in $(\ref{exit-mom2})$ and $(\ref{exit-mom1})$ we find that
\begin{align}
 \exp\left\{-N\left( W (K_i)- \min_i W (K_i) +\frac{2l-1}{4l}\varepsilon\right)\right\}\leq \wp^N(\gamma_i)\leq  \exp\left\{-N\left( W (K_i)- \min_i W (K_i) -\frac{2l-1}{4l}\varepsilon\right)\right\}.
 \label{bound-inv1}
\end{align}
By summing over $i=1,\ldots,l$, we get 
 \begin{align}
\wp^N(\mathcal{M}_1((\mathcal{Z}))^{2r}) \geq \exp\left\{-N\left( \frac{2l-1}{4l}\varepsilon\right)\right\}.
\label{low-inv}
\end{align}
 
 Using again the formula $(\ref{khas-form})$ together with the definition of the stopping times $\sigma_n$ and $\tau_n$ we find 
 \begin{align*}
 \wp^N(\mathcal{M}_1((\mathcal{Z}))^{2r})=\int_{\partial \gamma} \vartheta^N(d\nu)\mathbb{E}_{\nu}\left[\tau_1\right]&=\int_{\partial \gamma} \vartheta^N(d\nu)\left(\mathbb{E}_{\nu}\left[\sigma_0\right]+\mathbb{E}_{\nu}\left[\mathbb{E}_{\mu^N(\sigma_0)}[\tau_1]\right]\right)\\
 &\leq \sup_{\nu\in\partial \gamma}\mathbb{E}_{\nu}\left[\sigma_0\right]+\sup_{x\in C}\mathbb{E}_{x}[\tau_1].
 \end{align*}
By inequality $(\ref{exit-mom1})$ we find that $\sup_{\nu\in\partial \gamma}\mathbb{E}_{\nu}\left[\sigma_0\right]<\exp\{N\frac{\varepsilon}{4l}\}$ and by Corollary $\ref{coro-appen}$ we have that $\sup_{x\in C}\mathbb{E}_{x}[\tau_1]$ is bounded above by some constant. Using this together with the lower bound in $(\ref{low-inv})$ we obtain, by dividing $\wp^N(\gamma_i)$ in $(\ref{bound-inv1})$ by $(\wp^N(\mathcal{M}_1((\mathcal{Z}))^{2r}))$, the lower and upper bounds for the normalized measure $\wp^N$ of the $r_1$-neighborhood $\gamma_i$. The theorem is proved. \carre

\section{Technical proofs}
\label{proof-appen}
We give here the proof of some technical results. 

\subsection{Proof of Lemma \ref{T-rate-prop}}
\label{T-rate-prop-proof}
This is a mild generalization of \cite[Lem. 3.2]{Bork+Sund2012} to the multi-populations case. We prove each of the three assertions respectively. 
\begin{enumerate}
\item Let us ignore for the moment the $\mathcal{E}$ constraint by supposing that all possible transitions are allowed. Let $\xi,\nu\in(\mathcal{M}_1(\mathcal{Z}))^{2r}$ with $\mu(0)=\nu$ and $\mu(T)=\xi$ and consider the constant velocity path given by
\begin{align}
\mu(t)=(1-\frac{t}{T})\nu+\frac{t}{T}\xi,\quad\forall t\in[0,T],
\label{const-vel-path}
\end{align}
from which we can see that 
\begin{align}
\dot{\mu}(t)=\frac{\xi-\nu}{T}, \forall t\in[0,T],
\label{deriv-path}
\end{align}
a constant velocity. Note that $(\ref{const-vel-path})$ can be rewritten as follows: for all $1\leq j\leq r$ and $t$ in $[0,T]$ 
\begin{align}
\mu_j^c(t)=(1-\frac{t}{T})\nu_j^c+\frac{t}{T}\xi_j^c,\\
\mu_j^p(t)=(1-\frac{t}{T})\nu_j^p+\frac{t}{T}\xi_j^p.
\label{const-vel-path-2}
\end{align}

 Now we construct rate matrices $L_{j,c}(t)=(l^{j,c}_{z,z'}(t))_{z,z'\in\mathcal{Z}}$ and $L_{j,p}(t)=(l^{j,p}_{z,z'}(t))_{z,z'\in\mathcal{Z}}$ that ensure the traversal of this constant velocity path. Note that there is conservation of mass between the initial and terminal states, thus for all $1\leq j\leq r$, $\sum_z\nu_j^c(z)=\sum_z\xi_j^c(z)$ and $\sum_z\nu_j^p(z)=\sum_z\xi_j^p(z)$. Therefore, on can construct mass transport parameters ${g^{j,c}_{z,z'}}$ and ${g^{j,p}_{z,z'}}$ such that for any two states $z,z'\in\mathcal{Z}$ with $\nu_j^c(z) > \xi_j^c(z)$ (resp. $\nu_j^p(z) > \xi_j^p(z)$) and $z'$ with $\nu_j^c(z')<\xi_j^c(z')$ (resp. $\nu_j^p(z')<\xi_j^p(z')$), the quantity $g^{j,c}_{z,z'}$ (resp. $g^{j,p}_{z,z'}$) is the fraction of the excess of mass $\nu_j^{c}(z)-\xi_j^{c}(z)$ (resp. $\nu_j^{p}(z)-\xi_j^{p}(z)$) that goes from $z$ to $z'$. In particular, the coefficients $g^{j,\iota}_{z,z'}$ for $\iota\in\{c,p\}$ satisfy, for all $1\leq j\leq r$, the following conditions:
\begin{align}
 g^{j,\iota}_{z,z'} \in [0,1]\quad\text{for all $z,z'\in\mathcal{Z}$},
\end{align}
\begin{align}
 g^{j,\iota}_{z,z'}=0\quad\text{if $\nu_j^{\iota}(z)\leq\xi_j^{\iota}(z)$ or $\nu_j^{\iota}(z')\geq\xi_j^{\iota}(z')$},
\end{align}
\begin{align}
\sum_{z':\nu_j^{\iota}(z')<\xi_j^{\iota}(z')} g^{j,\iota}_{z,z'}=1\quad\text{if $\nu_j^{\iota}(z)> \xi_j^{\iota}(z)$},
\end{align}
and
\begin{align}
 \sum_{z:\nu_j^{\iota}(z)>\xi_j^{\iota}(z)} [\nu_j^{\iota}(z)-\xi_j^{\iota}(z)]g^{j,\iota}_{z,z'}=\xi_j^{\iota}(z')-\nu_j^{\iota}(z')\text{ if $\nu_j^{\iota}(z') < \xi_j^{\iota}(z')$}.
\end{align}
 
The first condition follows from the definition of the fractions $g^{j,\iota}_{z,z'}$. The second attests that no mass is transferred from a state with no mass excess, and no mass is received by a state with mass excess. The third point tells us that there is no mass destruction, and the last point stipulates that no mass is created. By using these mass transfer parameters, we construct the rate matrices $L_{j,c}(t)$ and $L_{j,p}(t)$ as follows. The diagonal elements are given, for $\iota\in\{c,p\}$, by 
\begin{align}
l^{j,\iota}_{z,z}(t)=\frac{-(\nu_j^{\iota}(z)-\xi_j^{\iota}(z))}{T(\mu_j^{\iota}(t)(z))},
\label{diag-rate}
\end{align} 

for all ${z}\in\mathcal{Z}$ that satisfy $\nu_j^{\iota}(z)> \xi_j^{\iota}(z)$ and $l^{j,\iota}_{z,z}(t)=0$ otherwise. The off-diagonals elements are given by 
\begin{align}
l^{j,\iota}_{z,z'}(t)=-l^{j,\iota}_{z,z}(t)g^{j,\iota}_{z,z'},
\label{off-diag-rate}
\end{align}
if $\nu_j^{\iota}(z)> \xi_j^{\iota}(z)$ and $z\neq z'$, and $l^{j,\iota}_{z,z'}(t)=0$ if $\nu_j^{\iota}(z)\leq\xi_j^{\iota}(z)$, for any $z'\in\mathcal{Z}$. Using these definitions of the rate matrices $L_{j,c}(t)$ and $L_{j,p}(t)$, and the properties of the mass transport coefficients $g^{j,c}_{z,z'}$ and $g^{j,p}_{z,z'}$, it is easy to prove  that  for $1\leq j\leq r$ and $t\in[0,T]$ (see \cite[p. 360]{Bork+Sund2012} ),
\begin{align*}
\dot{\mu}_j^c(t)&=L_{j,c}(t)^*\mu_j^c(t),\\
\dot{\mu}_j^p(t)&=L_{j,p}(t)^*\mu_j^p(t).
\end{align*}

 Let us now evaluate the difficulty of the passage $S_{[0,T]}(\mu|\nu)$ at this constant velocity path $\mu$. Theorem $\ref{large-dev-emp-proc}$ tells us that, if  $S_{[0,T]}(\mu|\nu)<\infty$, then $S_{[0,T]}(\mu|\nu)$ is given by $(\ref{rate-emp-proc-2})$ and this is in particular true if the $2r$ integral terms in $(\ref{rate-emp-proc-2})$ are finite. Notice that if, for a given $j$, $\mu^{\iota}_j(T)=\xi_j^{\iota}(z)=0$ and $\nu_j^{\iota}(z)>0$ for $z\neq z'$, then $l^{j,\iota}_{z,z'}(T)$ is not bounded.
Therefore, we need to provide a bound for $(\ref{rate-emp-proc-2})$ in the case of the constant velocity path. We start by upper bounding the sums in $(\ref{rate-emp-proc-2})$ for which the rates $l^{j,\iota}_{z,z'}(t)$ are strictly positive on $[0,T]$. To this end we introduce the set
\begin{align*}
\Upsilon_j^{\iota}= \{(z,z') | z\neq z', \nu_j^{\iota}(z) > \xi_j^{\iota}(z), \nu_j^{\iota}(z') < \xi_j^{\iota}(z'),g^{j,\iota}_{z,z'}>0\}.
\end{align*}
Thus, from the definition of the Legendre transform $\tau^*$, the right hand side of $(\ref{rate-emp-proc-2})$ can be written as
\begin{equation}
\begin{split}
\sum_{j=1}^r\bigg[&\alpha_jp_j^c \int_0^T \bigg( \sum_{(z,z')\in\Upsilon_j^{c}}\bigg((\mu_j^c(t)(z))l_{z,z'}^{j,c}(t) \log\bigg(\frac{l_{z,z'}^{j,c}(t)}{\lambda^c_{z,z'}\big(\mu_j^c(t),\mu_j^p(t)\big)}\bigg)-(\mu_j^c(t)(z))l_{z,z'}^{j,c}(t)\\ 
&\qquad\qquad\qquad\qquad\qquad+(\mu_j^c(t)(z))\lambda^c_{z,z'}\big(\mu_j^c(t),\mu_j^p(t)\big) \bigg)\bigg)dt \\
 &+\alpha_jp_j^p\int_0^T  \bigg( \sum_{(z,z')\in\Upsilon_j^{p}}\bigg((\mu_j^p(t)(z))l_{z,z'}^{j,p}(t) \log\bigg(\frac{l_{z,z'}^{j,p}(t)}{\lambda^p_{z,z'}\big(\mu_j^c(t),\mu_1^p(t),\ldots,\mu_r^p(t)\big)}\bigg)-(\mu_j^p(t)(z))l_{z,z'}^{j,p}(t)\\
 &\qquad\qquad\qquad\qquad\qquad+(\mu_j^p(t)(z))\lambda^p_{z,z'}\left(\mu_j^c(t),\mu_1^p(t),\ldots,\mu_r^p(t)\right)\bigg) \bigg) dt \bigg].
\label{rate-emp-proc-2-loc}
\end{split}
\end{equation}

 From $(\ref{diag-rate})$ and $(\ref{off-diag-rate})$ we find that $(\mu_j^{\iota}(t)(z))l_{z,z'}^{j,\iota}(t)=T^{-1}(\nu_j^{\iota}(z)-\xi_j^{\iota}(z))g^{j,\iota}_{z,z'}$. Moreover, from Assumption $\ref{ass-prin}$ we obtain $|\log\lambda^{\iota}_{z,z'}(\cdot,\cdot)|\leq |\log c|+|\log C|$. Therefore, $(\ref{rate-emp-proc-2-loc})$ is upper bounded by 
 \begin{equation}
\begin{split}
\sum_{j=1}^r\bigg[&\alpha_jp_j^c \int_0^T \bigg( \sum_{(z,z')\in\Upsilon_j^{c}}\bigg((T^{-1}(\nu_j^{c}(z)-\xi_j^{c}(z))g^{j,c}_{z,z'} \log\bigg(\frac{(\nu_j^{c}(z)-\xi_j^{c}(z))g^{j,c}_{z,z'}}{T(\mu_j^{c}(t)(z))}\bigg)\\ 
&\qquad\qquad\qquad\qquad\qquad+T^{-1}(\nu_j^{c}(z)-\xi_j^{c}(z))(|\log c|+|\log C|+1\big)+C \bigg)\bigg)dt \\
 &+\alpha_jp_j^p\int_0^T \bigg( \sum_{(z,z')\in\Upsilon_j^{p}}\bigg((T^{-1}(\nu_j^{p}(z)-\xi_j^{p}(z)) g^{j,p}_{z,z'}\log\bigg(\frac{(\nu_j^{p}(z)-\xi_j^{p}(z))g^{j,p}_{z,z'}}{T(\mu_j^{p}(t)(z))}\bigg)\\ 
&\qquad\qquad\qquad\qquad\qquad+T^{-1}(\nu_j^{p}(z)-\xi_j^{p}(z))(|\log c|+|\log C|+1\big)+C \bigg)\bigg)dt  \bigg],
\label{rate-emp-proc-2-loc-2}
\end{split}
\end{equation}
 which in turn is bounded by 
 \begin{equation}
\begin{split}
\sum_{j=1}^r\bigg[&\alpha_jp_j^c  \bigg( \sum_{(z,z')\in\Upsilon_j^{c}}(\nu_j^{c}(z)-\xi_j^{c}(z)) g^{j,c}_{z,z'}\big|\log\big((\nu_j^{c}(z)-\xi_j^{c}(z))g^{j,c}_{z,z'}\big)\big|\\ 
&-\sum_{(z,z')\in\Upsilon_j^{c}}T^{-1}(\nu_j^{c}(z)-\xi_j^{c}(z)) g^{j,c}_{z,z'}\int_0^T\log(\mu_j^{c}(t)(z))dt\\
&\qquad\qquad\qquad+\|\nu_j^{c}-\xi_j^{c}\||\log T|+\|\nu_j^{c}-\xi_j^{c}\|(|\log c|+|\log C|+1\big)+CTK^2 \bigg) \\
 &+\alpha_jp_j^p \bigg( \sum_{(z,z')\in\Upsilon_j^{p}}(\nu_j^{p}(z)-\xi_j^{p}(z)) g^{j,p}_{z,z'}\big|\log\big((\nu_j^{p}(z)-\xi_j^{p}(z))g^{j,p}_{z,z'}\big)\big|\\ 
&-\sum_{(z,z')\in\Upsilon_j^{p}}T^{-1}(\nu_j^{p}(z)-\xi_j^{p}(z)) g^{j,p}_{z,z'}\int_0^T\log(\mu_j^{p}(t)(z))dt\\
&\qquad\qquad\qquad+\|\nu_j^{p}-\xi_j^{p}\||\log T|+\|\nu_j^{p}-\xi_j^{p}\|(|\log c|+|\log C|+1\big)+CTK^2\bigg)  \bigg],
\label{rate-emp-proc-2-loc-2}
\end{split}
\end{equation} 
with $\|\cdot\|$ referring to the total variation distance. Observe that for $x\in [0,1]$, the function $x|\log x|$ is upper bounded, say by $\digamma>0$. Moreover, for a fixed $z$, $\sum_{z':z'\neq z}g_{zz'}^{j,\iota}=1$. Therefore the terms $\sum_{(z,z')\in\Upsilon_j^{\iota}}(\nu_j^{\iota}(z)-\xi_j^{\iota}(z)) g^{j,\iota}_{z,z'}\big|\log\big((\nu_j^{\iota}(z)-\xi_j^{\iota}(z))g^{j,\iota}_{z,z'}\big)\big| $ are bounded by 
\begin{equation}
\begin{split}
&\sum_{z:\nu_j^{\iota}(z)>\xi_j^{\iota}(z)}(\nu_j^{\iota}(z)-\xi_j^{\iota}(z)) \big|\log\big(\nu_j^{\iota}(z)-\xi_j^{\iota}(z)\big)\big|\bigg(\sum_{z':(z,z')\in\Upsilon_j^{\iota}}g^{j,\iota}_{z,z'}\bigg)\\
&\qquad\qquad+\sum_{(z,z')\in\Upsilon_j^{\iota}}(\nu_j^{\iota}(z)-\xi_j^{\iota}(z)) g^{j,\iota}_{z,z'}\big|\log g^{j,\iota}_{z,z'}\big|\\
&\leq \sum_{z:\nu_j^{\iota}(z)>\xi_j^{\iota}(z)}(\nu_j^{\iota}(z)-\xi_j^{\iota}(z)) \big|\log\big(\nu_j^{\iota}(z)-\xi_j^{\iota}(z)\big)\big|+\|\nu_j^{\iota}-\xi_j^{\iota}\| \digamma.
\end{split}
\end{equation}  

On the other hand, using the change of variable $u=\mu_j^{\iota}(t)(z)$ and then $(\ref{deriv-path})$, one obtains 
\begin{equation}
\begin{split}
\int_0^T\log(\mu_j^{\iota}(t)(z))dt&=T(\xi_j^{\iota}(z)-\nu_j^{\iota}(z))^{-1}\int_{\nu_j^{\iota}(z)}^{\xi_j^{\iota}(z)}\log udu\\
                 &=T(\xi_j^{\iota}(z)-\nu_j^{\iota}(z))^{-1}[u\log u-u]_{\nu_j^{\iota}(z)}^{\xi_j^{\iota}(z)}.
\end{split}
\end{equation} 

Therefore $(\ref{rate-emp-proc-2-loc-2})$ is upper bounded by 
 \begin{equation}
\begin{split}
\sum_{j=1}^r\bigg[&\alpha_jp_j^c  \bigg( \sum_{z:\nu_j^{c}(z)>\xi_j^{c}(z)}(\nu_j^{c}(z)-\xi_j^{c}(z)) \big|\log\big(\nu_j^{c}(z)-\xi_j^{c}(z)\big)\big|+\|\nu_j^{c}-\xi_j^{c}\| \digamma\\ 
&+\sum_{ (z,z')\in\Upsilon_j^{c}}|\xi_j^c(z)\log \xi_j^c(z)-\nu_j^c(z)\log \nu_j^c(z)| +\|\nu_j^c-\xi_j^c\| \\
&\qquad\qquad\qquad+\|\nu_j^{c}-\xi_j^{c}\||\log T|+\|\nu_j^{c}-\xi_j^{c}\|(|\log c|+|\log C|+1\big)+CTK^2 \bigg) \\
 &+\alpha_jp_j^p \bigg( \sum_{z:\nu_j^{p}(z)>\xi_j^{p}(z)}(\nu_j^{p}(z)-\xi_j^{p}(z)) \big|\log\big(\nu_j^{p}(z)-\xi_j^{p}(z)\big)\big|+\|\nu_j^{p}-\xi_j^{p}\| \digamma\\ 
&+\sum_{ (z,z')\in\Upsilon_j^{p}}|\xi_j^p(z)\log \xi_j^p(z)-\nu_j^p(z)\log \nu_j^p(z)| +\|\nu_j^p-\xi_j^p\| \\
&\qquad\qquad\qquad+\|\nu_j^{p}-\xi_j^{p}\||\log T|+\|\nu_j^{p}-\xi_j^{p}\|(|\log c|+|\log C|+1\big)+CTK^2\bigg)  \bigg]\\
&\leq \sum_{j=1}^r(\alpha_jp_j^cC'(T)+\alpha_jp_j^pC"(T)),
\label{rate-emp-proc-2-loc-3}
\end{split}
\end{equation}   
where $C'(T)$ and $C"(T)$ are two constants that do not depend on $\nu_j^{\iota}$ and $\xi_j^{\iota}$ for any $\iota\in\{c,p\}$ and $1\leq j\leq r$. 

Let us consider now the case where the pairs $(z,z')\notin\Upsilon_j^{\iota}$, i.e. for which $l_{z,z'}^{j,\iota}=0$. Using $\tau^{\star}(-1)=1$, the corresponding integral terms in $(\ref{rate-emp-proc-2})$ become 
\begin{equation}
\begin{split}
&\int_0^T  \sum_{(z,z')\notin\Upsilon_j^{c}}\mu_j^{c}(t)(z))\lambda^{c}_{z,z'}(\mu_j^{c}(t),\mu_j^{p}(t))dt\leq K^2CT,\\
&\int_0^T  \sum_{(z,z')\notin\Upsilon_j^{c}}\mu_j^{p}(t)(z))\lambda^{p}_{z,z'}\left(\mu_j^{c}(t),\mu_1^{p}(t),\ldots,\mu_r^{p}(t)\right)dt\leq K^2CT.
\label{zero-term-bound}
\end{split}
\end{equation}

Combining this with $(\ref{rate-emp-proc-2-loc-3})$ gives  $S_{[0,T]}(\mu|\nu)\leq \sum_{j=1}^r(\alpha_jp_j^c\bar{C}'(T)+\alpha_jp_j^p\bar{C}"(T))$ for some positive constants $\bar{C}'$ and $\bar{C}"$ that do not depend on $\nu_j^{\iota}$ and $\xi_j^{\iota}$.  
    
 Finally, let us consider the case where only transitions in $\mathcal{E}$ are allowed. Since the directed graph $(\mathcal{Z},\mathcal{E})$ is irreducible, the Markov chain $(\mu(t),t\geq 0)$ is also irreducible and thus there exist a finite
sequence of intermediate points through which one can move from $\nu$ to $\xi$ in $m=m(|\mathcal{Z}|,\mathcal{E})<+\infty$ steps
\begin{align*}
\nu= \nu^{(0)} \rightarrow\nu^{(1)}\rightarrow\cdots\rightarrow\nu^{(m)}=\xi.
\end{align*}

Therefore, one can construct a piecewise linear path $\mu$ with constant velocities on each of the $m$ segments such that each segment is covered in time duration $T/m$. Indeed, define for each $k=0,\ldots,m-1$, the constant velocity path
\begin{align*}
\mu_{k+1}(t)=(1-\frac{tm}{T})\nu^{(k)}+\frac{tm}{T}\nu^{(k+1)},\quad\forall t\in[kT/m,(k+1)T/m],
\end{align*}
and take $\mu(t)=\mu_{k+1}(t)$ for $t\in[kT/m,(k+1)T/m]$. Hence 
\begin{align*}
S_{[0,T]}(\mu|\nu)\leq C_1(T)= m\sum_{j=1}^r(\alpha_jp_j^c\bar{C}'(T/m)+\alpha_jp_j^p\bar{C}"(T/m)),
\end{align*}
 which completes the proof of the first assertion.  

\item This follows immediately from the previous result and $(\ref{T-rate})$.

\item Fix $\varepsilon>0$ and take $T=\varepsilon$ in $(\ref{rate-emp-proc-2-loc-3})$. We can always find a $\delta\in(0,\varepsilon)$ such that, if  $\rho^{2r}_0(\nu,\xi)<\delta$, then $(\ref{rate-emp-proc-2-loc-3})$ is bounded by 
 \begin{equation}
\begin{split}
\sum_{j=1}^r\bigg[&\alpha_jp_j^c  \bigg( 6\varepsilon+CK^2\varepsilon \bigg) +\alpha_jp_j^p\bigg( 6\varepsilon+CK^2\varepsilon\bigg)  \bigg].
\end{split}
\end{equation} 

Combining this with $(\ref{zero-term-bound})$ gives that, for any $\nu,\xi$ such that $\rho^{2r}_0(\nu,\xi)<\delta$,
 \begin{equation}
\begin{split}
S_{[0,\varepsilon]}(\xi|\nu)\leq \sum_{j=1}^r\bigg[&\alpha_jp_j^c  \bigg( 6\varepsilon+2CK^2\varepsilon \bigg) +\alpha_jp_j^p\bigg( 6\varepsilon+2CK^2\varepsilon\bigg)  \bigg].
\end{split}
\end{equation} 

The result then follows from $(\ref{T-rate})$.
\end{enumerate}
\carre

\subsection{Proof of Lemma \ref{unif-cont-lem-1}} 
\label{unif-cont-lem-1-proof}
We generalize \cite[Lem. 7.1]{Bork+Sund2012}. From Theorem $\ref{large-dev-emp-proc}$, since $S_{[0,T]} (\mu|\nu)<+\infty$, the rate function $S_{[0,T]} (\mu|\nu)$ is given by $(\ref{rate-emp-proc-2})$. Moreover, we can verify that $\tau^*(u-1)= u \log u - u + 1 \geq  u- e + 1$ for all $u \geq 0$. Using this together with $(\ref{rate-emp-proc-2})$ gives to us
\begin{align*}
S_{[0,T ]} (\mu|\nu)&\geq \sum_{j=1}^r\bigg[\alpha_jp_j^c \int_0^T \bigg( \sum_{(z,z')\in\mathcal{E}}(\mu_j^c(t)(z))\lambda^c_{z,z'}(\mu_j^c(t),\mu_j^p(t)) \bigg(\frac{l_{z,z'}^{j,c}(t)}{\lambda^c_{z,z'}(\mu_j^c(t),\mu_j^p(t))}-e+1\bigg) \bigg)dt \\
 &\quad+\alpha_jp_j^p\int_0^T  \bigg( \sum_{(z,z')\in\mathcal{E}}(\mu_j^p(t)(z))\lambda^p_{z,z'}(\mu_j^c(t),\mu_1^p(t),\ldots,\mu_r^p(t)) \bigg(\frac{l_{z,z'}^{j,p}(t)}{\lambda^p_{z,z'}(\mu_j^c(t),\mu_1^p(t),\ldots,\mu_r^p(t))}-e+1\bigg) \bigg) dt \bigg]\\
 &\geq \sum_{j=1}^r\bigg[\alpha_jp_j^c \int_0^T \bigg( \sum_{(z,z')\in\mathcal{E}}(\mu_j^c(t)(z))l_{z,z'}^{j,c}(t)  \bigg)dt+\alpha_jp_j^p\int_0^T  \bigg( \sum_{(z,z')\in\mathcal{E}}(\mu_j^p(t)(z)) l_{z,z'}^{j,p}(t) \bigg) dt\bigg]\\
 &\quad-\sum_{j=1}^r\bigg[\alpha_jp_j^c (e-1)CT|\mathcal{E}|\bigg]+\sum_{j=1}^r\bigg[\alpha_jp_j^p (e-1)CT|\mathcal{E}| \bigg],\\
\end{align*}
which completes the proof. \carre

\subsection{Proof of Lemma \ref{unif-cont-lem-2}}
\label{unif-cont-lem-2-proof}
 We generalize \cite[Lem. 7.2]{Bork+Sund2012}. By the construction we have that $\tilde{\mu}= \mu(0) = \nu$ and $\tilde{\mu}(T')= \mu(\beta T')= \mu(T)=\xi$. Fix $1\leq j\leq r$. Using $\dot{\mu}_j^c(t)=L_{j,c}(t)^*\mu_j^c(t)$ we find that, for all $t\in[0,T']$,
\begin{align*}
\dot{\tilde{\mu}}_j^c(t)=\frac{d\tilde{\mu}_j^c( t)}{dt}=\frac{d \mu_j^c( \beta t)}{dt}=\beta \dot{\mu}_j^c(\beta t)=\beta L_{j,c}(\beta t)^*\mu_j^c(\beta t)=\beta L_{j,c}(\beta t)^*\dot{\tilde{\mu}}_j^c(t),
\end{align*}
from which we deduce that $\tilde{L}_{j,c}(t) = \beta L_{j,c}(\beta t)$. Using the exact same steps leads to $\tilde{L}_{j,p}(t) = \beta L_{j,p}(\beta t)$. Let us now evaluate the rate function $S_{[0,T']} (\tilde{\mu}|\nu)$ associated with the path $\tilde{\mu}:[0,T]\rightarrow(\mathcal{M}_1(\mathcal{Z}))^{2r}$. From $(\ref{rate-emp-proc-2})$ we have, 
\begin{equation}
\begin{split}
S_{[0,T']} (\tilde{\mu}|\nu)&=\sum_{j=1}^r\bigg[\alpha_jp_j^c \int_0^{T'} \bigg( \sum_{(z,z')\in\mathcal{E}}(\tilde{\mu}_j^c(t)(z))\lambda^c_{z,z'}(\tilde{\mu}_j^c(t),\tilde{\mu}_j^p(t)) \tau^*\bigg(\frac{\tilde{l}_{z,z'}^{j,c}(t)}{\lambda^c_{z,z'}(\tilde{\mu}_j^c(t),\tilde{\mu}_j^p(t))}-1\bigg) \bigg)dt \\
 &+\alpha_jp_j^p\int_0^{T'}  \bigg( \sum_{(z,z')\in\mathcal{E}}(\tilde{\mu}_j^p(t)(z))\lambda^p_{z,z'}(\tilde{\mu}_j^c(t),\tilde{\mu}_1^p(t),\ldots,\tilde{\mu}_r^p(t)) \tau^*\bigg(\frac{\tilde{l}_{z,z'}^{j,p}(t)}{\lambda^p_{z,z'}(\tilde{\mu}_j^c(t),\tilde{\mu}_1^p(t),\ldots,\tilde{\mu}_r^p(t))}-1\bigg) \bigg) dt \bigg]\\
 &=\sum_{j=1}^r\bigg[\alpha_jp_j^c \int_0^{T'} \bigg( \sum_{(z,z')\in\mathcal{E}}(\mu_j^c(\beta t)(z))\lambda^c_{z,z'}(\mu_j^c(\beta t),\mu_j^p(\beta t))\\
 & \qquad\qquad\qquad\qquad\qquad\qquad\times\tau^*\bigg(\frac{\beta l_{z,z'}^{j,c}(\beta t)}{\lambda^c_{z,z'}(\mu_j^c(\beta t),\mu_j^p(\beta t))}-1\bigg) \bigg)dt \\
 &\qquad\quad+\alpha_jp_j^p\int_0^{T'}  \bigg( \sum_{(z,z')\in\mathcal{E}}(\mu_j^p(\beta t)(z))\lambda^p_{z,z'}(\mu_j^c(\beta t),\mu_1^p(\beta t),\ldots,\mu_r^p(\beta t))\\
 &\qquad\qquad\qquad\qquad\qquad\qquad\times \tau^*\bigg(\frac{\beta l_{z,z'}^{j,p}(\beta t)}{\lambda^p_{z,z'}(\mu_j^c(\beta t),\mu_1^p(\beta t),\ldots,\mu_r^p(\beta t))}-1\bigg) \bigg) dt \bigg].
\end{split}
\end{equation}
From $(\ref{Legen-tran})$ and using the properties of the $\log$ function we can easily verify that
\begin{align*}
\tau^{\star}(\beta u-1)=\beta\bigg(u\log \beta+\tau^{*}(u-1)+\frac{1-\beta}{\beta}\bigg),\quad\forall u\geq 0.
\end{align*}
Therefore we find 
\begin{equation}
\begin{split}
S_{[0,T']} (\tilde{\mu}|\nu) &=\sum_{j=1}^r\bigg[\alpha_jp_j^c \int_0^{T'} \bigg( \sum_{(z,z')\in\mathcal{E}}(\mu_j^c(\beta t)(z))\lambda^c_{z,z'}(\mu_j^c(\beta t),\mu_j^p(\beta t))\times\\
&\beta\left\{\frac{ l_{z,z'}^{j,c}(\beta t)}{\lambda^c_{z,z'}(\mu_j^c(\beta t),\mu_j^p(\beta t))}(\log\beta)+\tau^{*}\left(\frac{ l_{z,z'}^{j,c}(\beta t)}{\lambda^c_{z,z'}(\mu_j^c(\beta t),\mu_j^p(\beta t))}-1\right)+\frac{1-\beta}{\beta}\right\}\bigg)dt \\
 &\qquad\quad+\alpha_jp_j^p\int_0^{T'}  \bigg( \sum_{(z,z')\in\mathcal{E}}(\mu_j^p(\beta t)(z))\lambda^p_{z,z'}(\mu_j^c(\beta t),\mu_1^p(\beta t),\ldots,\mu_r^p(\beta t))\times\\
 &\beta\left\{\frac{ l_{z,z'}^{j,p}(\beta t)}{\lambda^p_{z,z'}(\mu_j^c(\beta t),\mu_1^p(\beta t),\ldots,\mu_r^p(\beta t))}(\log\beta)+\tau^*\bigg(\frac{ l_{z,z'}^{j,p}(\beta t)}{\lambda^p_{z,z'}(\mu_j^c(\beta t),\mu_1^p(\beta t),\ldots,\mu_r^p(\beta t))}-1\bigg)+\frac{1-\beta}{\beta}\right\}\bigg) dt \bigg].
\end{split}
\end{equation}

Introducing the change of variables $\beta t\rightarrow t$ we obtain
\begin{equation}
\begin{split}
S_{[0,T']} (\tilde{\mu}|\nu) =S_{[0,T]} (\mu|\nu)+ \sum_{j=1}^r\bigg[&\frac{1-\beta}{\beta}\alpha_jp_j^c \int_0^{T} \bigg( \sum_{(z,z')\in\mathcal{E}}(\mu_j^c( t)(z))\lambda^c_{z,z'}(\mu_j^c( t),\mu_j^p( t))\bigg)dt\\
&+ (\log\beta)\alpha_jp_j^c \int_0^{T} \bigg( \sum_{(z,z')\in\mathcal{E}}(\mu_j^c( t)(z))l_{z,z'}^{j,c}(t)\bigg)dt \\
 &+\frac{1-\beta}{\beta}\alpha_jp_j^p\int_0^{T}  \bigg( \sum_{(z,z')\in\mathcal{E}}(\mu_j^p( t)(z))\lambda^p_{z,z'}(\mu_j^c( t),\mu_1^p(t),\ldots,\mu_r^p( t))\bigg)dt\\
 &+(\log\beta)\alpha_jp_j^p \int_0^{T} \bigg( \sum_{(z,z')\in\mathcal{E}}(\mu_j^p( t)(z))l_{z,z'}^{j,p}(t)\bigg)dt  \bigg].
\end{split}
\end{equation}

Finally, using Assumption $\ref{ass-prin}$ we find $(\ref{bound-time-scal})$. \carre

\end{appendices}

\section*{Acknowledgment}
This research was supported by the Natural Sciences and Engineering Research Council of Canada Discovery Grants and by Carleton University.

\bibliographystyle{livre} 
\bibliography{biblio}

\end{document}